\documentclass{amsart}

\usepackage{amsopn}
\usepackage{amsmath}
\usepackage{amssymb}
\usepackage{amsfonts}
\usepackage{latexsym}
\usepackage{graphics}
\usepackage{graphicx}
\usepackage{hyperref}
\usepackage{color}
\usepackage{amscd}

\newtheorem{theorem}{Theorem}[section]
\newtheorem{lemma}[theorem]{Lemma}
\newtheorem{definition}[theorem]{Definition}
\newtheorem{corollary}[theorem]{Corollary}

\newtheorem{example}[theorem]{Example}

\newtheorem{observation}[theorem]{Observation}

\def\Proof{\noindent{\bf Proof.\ \ }}

\newcommand{\comment}[1]{}

\begin{document}

\title[Generalizing the Splits Equivalence Theorem]{Generalizing the Splits Equivalence Theorem and Four Gamete Condition: Perfect Phylogeny on Three State Characters}
\date{\today}
\author{Fumei Lam} 
\address{Department of Computer Science, University of California, Davis}
\email{flam@cs.ucdavis.edu}
\author{Dan Gusfield}
\address{Department of Computer Science, University of California, Davis}
\email{gusfield@cs.ucdavis.edu}
\author{Srinath Sridhar}
\address{Department of Computer Science, Carnegie Mellon University}
\email{srinath@cs.cmu.edu}

\maketitle

%\keywords{Four gametes, splits equivalence, perfect phylogeny, multiple state}

\begin{abstract} We study the three state perfect phylogeny problem and establish a generalization of the four gamete condition (also called the Splits Equivalence Theorem) for sequences over three state characters.  Our main result is that a set of input sequences over three state characters allows a perfect phylogeny if and only if every subset of three characters allows a perfect phylogeny.   In establishing these results, we prove fundamental structural features of the perfect phylogeny problem on three state characters and completely characterize the minimal obstruction sets that must occur in three state input sequences that do not have a perfect phylogeny.  We further give a proof for a stated lower bound involved in the conjectured generalization of our main result to any number of states.  

Until this work, the notion of a conflict, or incompatibility, graph has been defined for two state characters only.  Our generalization of the four gamete condition allows us to generalize the notion of incompatibility to three state characters. The resulting incompatibility structure is a hypergraph, which can be used to solve algorithmic and theoretical problems for three state characters. 
\end{abstract}

\section{Introduction}

% end comment

One of the fundamental problems in biology is the construction of phylogenies, or evolutionary trees, to describe ancestral relationships between a set of observed taxa.  Each taxon is represented by a sequence and the evolutionary tree provides an explanation of branching patterns of mutation events transforming one sequence into another.  There have been many elegant theoretical and algorithmic results on the problems of reconstructing a plausible history of mutations that generate a given set of observed sequences and determining the minimum number of such events needed to explain the sequences.  

A widely used model in phylogeny construction and population genetics is the {\it infinite sites} model, in which the mutation of any character can occur at most once in the phylogeny. This implies that the data must be binary (a character can take on at most two states), and that without recombination the phylogeny must be a tree, called a (binary) Perfect Phylogeny. The problem of determining if a set of binary sequences fits the infinite sites model without recombination, corresponds to determining if the data can be derived on a binary Perfect Phylogeny. A
generalization of the infinite sites model is the {\it infinite alleles} model, where any character can mutate many times but each mutation of the character must lead to a different allele (state). Again, without recombination, the phylogeny is tree, called a {\it multi-state} Perfect Phylogeny. Correspondingly, the problem of determining if multi-state data fits the infinite-alleles model without recombination corresponds to determining if the data can be derived on a multi-state perfect phylogeny.

In the case of binary sequences, the well-known Splits Equivalence Theorem (also known as the \emph{four gamete condition}) gives a necessary and sucient condition for the existence of a (binary) perfect phylogeny.

\begin{theorem}[Splits Equivalence Theorem, Four Gamete Condition \\ \cite{EJM76,Gu91,Me83}] \label{thm:fourgamete} A perfect phylogeny exists for binary input sequences if and only if no pair of characters contains all four possible binary pairs 00, 01, 10, 11. \end{theorem}

It follows from this theorem that for binary input, it is possible to either construct a perfect phylogeny, or output a pair of characters containing all four gametes as an obstruction set witnessing the nonexistence of a perfect phylogeny.  This test is the building block for many theoretical results and practical algorithms.  Among the many applications of this theorem, Gusfield et al.~\cite{GB05,GBBS07} and Huson et al.~\cite{H+05} apply the theorem to achieve decomposition theorems for phylogenies, Gusfield, Hickerson, and Eddhu~\cite{GHE07} Bafna and Bansal~\cite{BB05,BB04}, and Hudson and Kaplan~\cite{HK85} use it to obtain lower bounds for recombination events, Gusfield et al.~\cite{Gu05,GEL04} use it to obtain algorithms for constructing networks with constrained recombination, Sridhar et al.~\cite{BDH06,SDBHRS06,S+07} and Satya et al.~\cite{SMAPB06} use it to achieve a faster near-perfect phylogeny reconstruction algorithm, Gusfield~\cite{Gu02} uses it to infer phase inference (with subsequent papers by Gusfield et al. \cite{BGHY04,BGLY03,DFG06,GB05}, Eskin, Halperin, and Karp~\cite{EHK03,HE04}, Satya and Mukherjee~\cite{SM06} and Bonizzoni~\cite{B07}), and Sridhar~\cite{SBRS06} et al. use it to obtain phylogenies from genotypes.   

%While this list is not comprehensive, it demonstrates the power of the theorem. 

The focus of this work is to extend results for the binary perfect phylogeny problem to the multiple state character case, addressing the following natural questions arising from the Splits Equivalence Theorem.   Given a set of sequences on $r$ states ($r \geq 3$), is there a necessary and sufficient condition for the existence of a perfect phylogeny analogous to the Splits Equivalence Theorem?  If no perfect phylogeny exists, what is the size of the smallest witnessing obstruction set? 

 In 1975, Fitch gave an example of input $S$ over three states such that every \emph{pair} of characters in $S$ allows a perfect phylogeny while the entire set of characters $S$ does not~\cite{Fe04,F75,F77,SS03}.  In 1983, Meacham generalized these results to characters over $r$ states ($r \geq 3$)\cite{Me83}, constructing a class of sequences called \emph{Fitch-Meacham examples}, which we examine in detail in Section \ref{sec:fitchmeacham}.   Meacham writes:

\begin{quote}
``The Fitch examples show that any algorithm to determine whether a set of characters is compatible must consider the set as a whole and cannot take the shortcut of only checking pairs of characters." \cite{Me83}
\end{quote}

However, while the Fitch-Meacham construction does show that checking pairs of characters is not sufficient for the existence of a perfect phylogeny, our main result will show that for three state input, there is a sufficient condition which does \emph{not} need to consider the entire set of characters simultaneously.  In particular, we give a complete answer to the questions posed above for three state characters, by 

\begin{enumerate}
\item showing the existence of a necessary and sufficient condition analogous to the Splits Equivalence Theorem (Sections \ref{sec:mequals3} and \ref{sec:generalm}),
\item in the case no perfect phylogeny exists, proving the existence of a small obstruction set as a witness (Section \ref{sec:generalm}),
\item giving a complete characterization of all minimal obstruction sets (Section \ref{sec:forbidden}), and
\item proving a stated lower bound involved in the conjectured generalization of our main result to any number of states (Section \ref{sec:fitchmeacham}).
\end{enumerate}

In establishing these results, we prove fundamental structural features of the perfect phylogeny problem on three state characters.

\section{Perfect Phylogenies and Partition Intersection Graphs}

The input to our problem is a set of $n$ sequences (representing taxa), where each sequence is a string of length $m$ over $r$ states.  Throughout this paper, the states under consideration will be the set $\{ 0, 1, 2, \ldots r-1\}$ (in particular, in the case $r=2$, the input are binary sequences over $\{0, 1\}$).  The input can be considered as a matrix of size $n \times m$, where each row corresponds to a sequence and each column corresponds to a character (or site).  We denote characters by $\mathcal{C} = \{ \chi^1, \chi^2, \chi^3, \ldots\ \chi^m \}$ and the states of character $\chi^i$ by $\chi^i_j$ for $0 \leq j \leq r-1$.  A \emph{species} is a sequence $s_1, s_2, \ldots s_m \in \chi^1_{j_1} \times \chi^2_{j_2} \times \cdots \chi^m_{j_m}$, where $s_i$ is the \emph{state} of character $\chi^i$ for $s$.

The \emph{perfect phylogeny problem} is to determine whether an input set $S$ can be displayed on a tree such that

\begin{enumerate}
\item each sequence in input set $S$ labels exactly one leaf in $T$
\item each vertex of $T$ is labeled by a species
\item for every character $\chi^i$ and for every state $\chi^i_j$ of character $\chi^i$, the set of all vertices in $T$ such that the state of character $\chi^i$ is $\chi^i_j$ forms a connected subtree of $T$. 
\end{enumerate}

The general perfect phylogeny problem (with no constraints on $r$, $n$, and $m$) is NP-complete \cite{BFW92,St92}.  However, the perfect phylogeny problem becomes polynomially solvable (in $n$ and $m$) when $r$ is fixed.  For $r=2$, this follows from the Splits Equivalence Theorem \ref{thm:fourgamete}.  For larger values of $r$, this was shown by Dress and Steel for $r=3$ \cite{DS93}, by Kannan and Warnow for $r=3$ or $4$ \cite{KW94}, and by Agarwala and Fern\'andez-Baca for all fixed $r$ \cite{AF94} (with an improved algorithm by Kannan and Warnow \cite{KW97}). 

\begin{definition}[\cite{B74,SS03}] For a set of input sequences $S$, the \emph{partition intersection graph} $G(S)$ is obtained by associating a vertex for each character state and an edge between two vertices $\chi^i_j$ and $\chi^k_l$ if there exists a sequence $s$ with state $j$ in character $\chi^i \in \mathcal{C}$ and state $l$ in character $\chi^k \in \mathcal{C}$.  We say $s$ is a row that \emph{witnesses} edge $(\chi^i_j, \chi^k_l)$.  For a subset of characters $\Phi = \{ \chi^{i_1}, \chi^{i_2}, \ldots \chi^{i_k}\}$, let $G(\Phi)$ denote the partition intersection graph $G(S)$ restricted to the characters in $\Phi$.  \end{definition}

Note that by definition, there are no edges in the partition intersection graph between states of the same character.

\begin{definition} A graph $H$ is \emph{chordal}, or \emph{triangulated}, if there are no induced chordless cycles of length four or greater in $H$.
\end{definition}

Consider coloring the vertices of the partition intersection graph $G(S)$ in the following way.  For each character $\chi^i$, assign a single color to the vertices $\chi^i_0, \chi^i_1, \ldots \chi^i_{r-1}$.  A \emph{proper triangulation} of the partition intersection graph $G(S)$ is a chordal supergraph of $G(S)$ such that every edge has endpoints with different colors.  In \cite{B74}, Buneman established the following fundamental connection between the perfect phylogeny problem and triangulations of the corresponding partition intersection graph.

\begin{theorem} \label{thm:buneman} \cite{B74,SS03} A set of taxa $S$ 
admits a perfect phylogeny if and only if the corresponding partition 
intersection graph $G(S)$ has a proper triangulation. \end{theorem}

We will use Theorem \ref{thm:buneman} to extend the Splits Equivalence Theorem to a test for the existence of a perfect phylogeny on trinary state characters.  In a different direction, Theorem \ref{thm:buneman} and triangulation were also recently used to obtain an algorithm to handle perfect phylogeny problems with missing data \cite{Gu09}.

To outline our approach, suppose a perfect phylogeny exists for $S$ and consider every subset of three characters. Then each of these $\binom{m}{3}$ characters also has a perfect phylogeny.  We show that this necessary condition is also sufficient and moreover, we can systematically piece together the proper triangulations for each triple of characters to obtain a triangulation for the entire set of characters.  On the other hand, if no perfect phylogeny exists, then we show there exists a witness set of three characters for which no perfect phylogeny exists.  This extends the Splits Equivalence Theorem to show that for binary and trinary state input, the number of characters needed for a witness obstruction set is equal to the number of character states.   The following is the main theorem of the paper.

\begin{theorem} \label{thm:main} Given an input set $S$ on $m$ characters with at most three states per character ($r \leq 3$), $S$ admits a perfect phylogeny if and only if every subset of three characters of $S$ admits a perfect phylogeny. \end{theorem}

By this theorem, in order to verify that a trinary state input matrix $S$ has a perfect phylogeny, it suffices to verify that partition intersection graphs $G[\chi^i, \chi^j, \chi^k]$ have proper triangulations for all triples $\chi^i, \chi^j, \chi^k \in \mathcal{C}$.    In Section \ref{sec:fitchmeacham}, we will show that the Fitch-Meacham examples \cite{F75,Me83} demonstrate that the size of the witness set in Theorem \ref{thm:main} is best possible.

\comment{ 

Note that the existence of a unique triangulation does \emph{not} 
imply the existence of a unique perfect phylogeny.  Unlike the binary instance, it 
is possible for a set of sequences over three state characters to have several 
perfect phylogenies. 

A proper triangulation is \emph{minimal} if removing any edge results in 
a graph with an induced chordless cycle.

We note that a colored graph can have a unique triangulation but still 
be extended into a $k$-tree in multiple ways.  Therefore, a unique 
triangulation does not necessarily imply a unique perfect phylogeny, as 
the following example shows.

However, if it is the case that in the triangulation, all minimal 
separators of the graph have size $k$, then the graph must be a $k$-tree 
and the perfect phylogeny must be unique.

} % end comment

\section{Structure of Partition Intersection Graphs for Three Characters} \label{sec:mequals3}

We begin by studying the structure of partition intersection graphs on three characters with at most three states per character ($m \leq 3$, $r \leq 3$).  For convenience, we will denote the three characters by the letters $a,b,c$ (interchangeably referring to them as characters and colors) and denote the states of these characters by $a_i, b_i, c_i$ ($i \in \{0,1,2\}$).  

\comment{
For a given triple of characters $a, b,c$, it is possible to construct graph $G[a,b,c]$ as follows.  Since there are three characters over three states, each row is one of  $3^3=27$ possible sequences.   Scan the rows and count the number of times each sequence appears.  Now, create the nine vertices $a_i, b_i, c_i, i \in \{0,1,2\}$ in graph $G[a,b,c]$.  For each sequence $a_j, b_k, c_l$ appearing in the input, add the edges $(a_j, b_k), (b_k, c_l), (a_j, c_l)$ to $E(G[a,b,c])$ if not already present.  

For each row, consider the entires in the row restricted to columns $a,b$, and $c$.  There are $3^3 = 27$ possibilities for this 3-tuple and we can scan all $n$ rows in$O(n)$ time to count the number of times eachpossible $3$-tuple appears.  Then, for each $3$-tuple appearing at least once, add the edges between the corresponding states in the partition intersection graph if not already present. 
}

The problem of finding proper triangulations for graphs on at most three colors and arbitrary number of states ($m=3$, $r$ arbitrary) has been studied in a series of papers \cite{BK93,IS93,KW91}. However, it will be unnecessary in our problem to employ these triangulation algorithms, as our instances will be restricted to those arising from character data on at most three states ($m=3, r \leq 3$).  In such instances, we will show that if a proper triangulation exists, then the structure of the triangulation is very simple.  We begin by proving a sequence of lemmas characterizing the possible cycles contained in the partition intersection graph.  

\begin{lemma}\label{lemma:nothreestates} Let $S$ be a set of input species on three characters $a, b,$ and $c$ with at most three states per character.  Suppose every pair of characters induces a properly triangulatable partition intersection graph (i.e., $G[a,b]$, $G[b,c]$ and $G[a,c]$ are properly triangulatable) and let $C$ be a chordless cycle in $G[a,b,c]$.   Then $C$ cannot contain all three states of any character.  \end{lemma}

\noindent {\bf Proof.}  Suppose there is a color, say $a$, such that all three states $a_0, a_1$ and $a_2$ appear in $C$.   Note that $C$ must contain all three colors $a,b,$ and $c$ (since any pair of colors induces a properly triangulatable graph and any cycle on two colors cannot be properly triangulated).  We have the following cases.

\begin{figure}
        	\begin{center}
        	\includegraphics[scale=0.35]{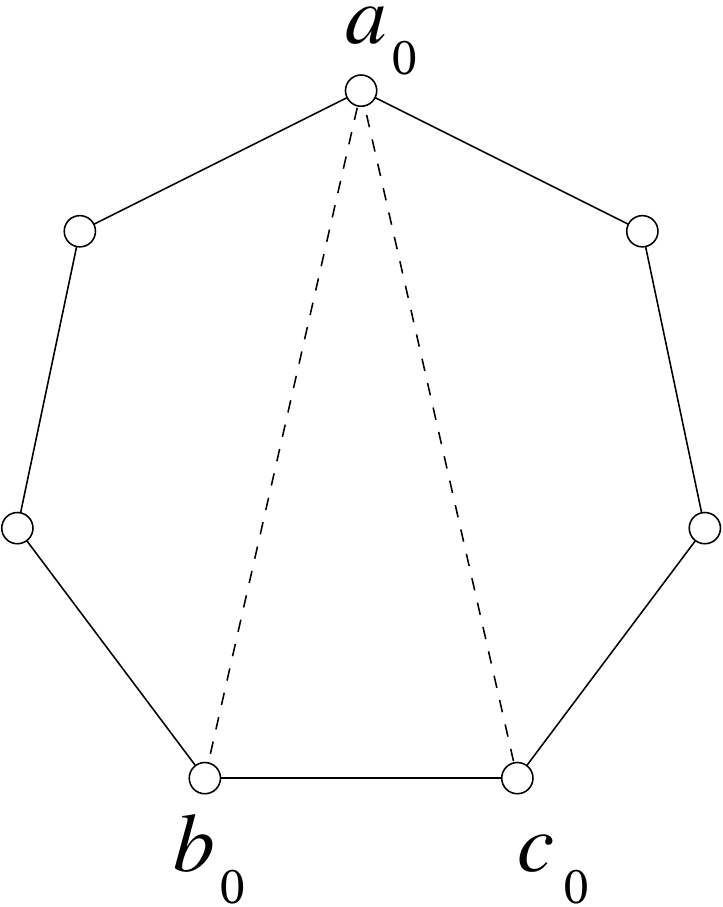}
        	\hspace*{0.5cm}
        	\caption{The row witnessing edge $(b_0, c_0)$ must contain a state in character $a$ . \label{fig:caseI} }
	\end{center}
\end{figure}

\noindent {\bf Case I.}   Suppose there is an edge $e$ in $C$ neither of whose endpoints have color $a$ (without loss of generality, let $e = (b_0, c_0)$).  The row that witnesses this edge must contain some state in $a$, say $a_0$.  This implies that the vertices $a_0, b_0,$ and $c_0$ form a triangle in $G[a,b,c]$, a contradiction since $C$ is assumed to be chordless (see Figure \ref{fig:caseI}).

\begin{figure}[h!]
        \centering
        \includegraphics[scale=0.35]{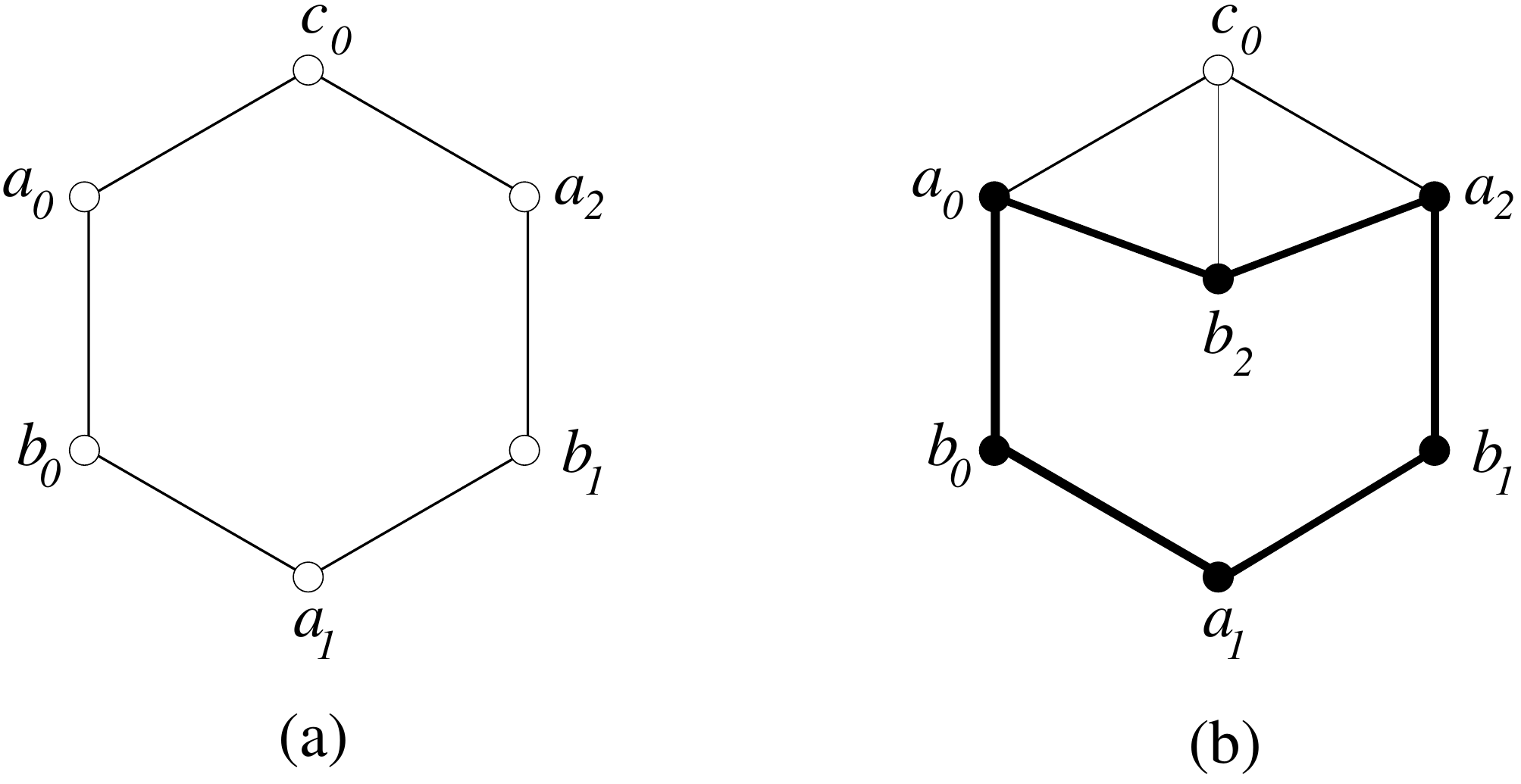}
        \caption{The row witnesses for edges $(a_0, c_0)$ and $(c_0, a_2)$ must share the same state of $b$.\label{fig:nothreestates}}
\end{figure}

\noindent {\bf Case II.}  Otherwise, every edge has an endpoint of color $a$, implying each edge has color pattern either $(a,b)$ or $(a,c)$.  Since all three states of $a$ appear, the color pattern up to relabeling must be as shown in Figure \ref{fig:nothreestates}(a) (in the figure, color $b$ appears twice and color $c$ appears once). 

In this case, the row witness for edge $(a_0, c_0)$ must contain the final state $b_2$ of $b$ (otherwise there would be an edge between $c_0$ and either $b_0$ or $b_1$, a contradiction since $C$ is chordless).  Similarly, the row witness for edge $(c_0, a_2)$ must also be state $b_2$.  As shown in Figure \ref{fig:nothreestates}(b), this gives a cycle $(a_0, b_2)$, $(b_2, a_2)$, $(a_2, b_1)$, $(b_1, a_1)$, $(a_1, b_0)$, $(b_0, a_0)$ on two colors.  Such a cycle is not properly triangulatable, and therefore $G[a,b]$ is not properly triangulatable, a contradiction.  

Since Case I and Case II cannot occur, it follows that $a_0, a_1$ and $a_2$ cannot all appear in $C$, proving the lemma. \qed

Before stating the next lemma, we give the following definition.

\begin{definition} Suppose the endpoints of edge $e$ have colors $\chi^i$ and $\chi^j$.  Then any other edge whose endpoints also have colors $\chi^i$ and $\chi^j$ is called \emph{color equivalent} to $e$.  
Two edges are called \emph{nonadjacent} if they do not share a common endpoint. \end{definition}

For example, the edges $(c_1, a_2)$ and $(c_0, a_1)$ in Figure \ref{fig:caseI} are color equivalent and nonadjacent.  

\begin{lemma}\label{lemma:unique}  Let $S$ be a set of input species on three characters $a, b,$ and $c$ with at most three states per character.  If the partition intersection graph $G[a,b,c]$ is properly triangulatable, then for every chordless cycle $C$ in $G[a,b,c]$, there exists a color ($a,b$, or $c$) that appears exactly once in $C$.  \end{lemma}

\Proof Consider any chordless cycle $C$ of $G[a,b,c]$.  By Lemma \ref{lemma:nothreestates}, no color appears in all three states in $C$. To obtain a contradiction, suppose each color $a,b$, and $c$ appears exactly twice in $C$ and without loss of generality, relabel the states so that the vertices appearing on the cycle are $a_0, a_1, b_0, b_1, c_0,$ and $c_1$.   We first show that $C$ has a pair of nonadjacent edges that are color equivalent.  Up to symmetry and relabelling of colors, there are two cases for the color pattern of $C$ as follows.

\begin{figure}[h!]
        \centering
        \includegraphics[scale=0.35]{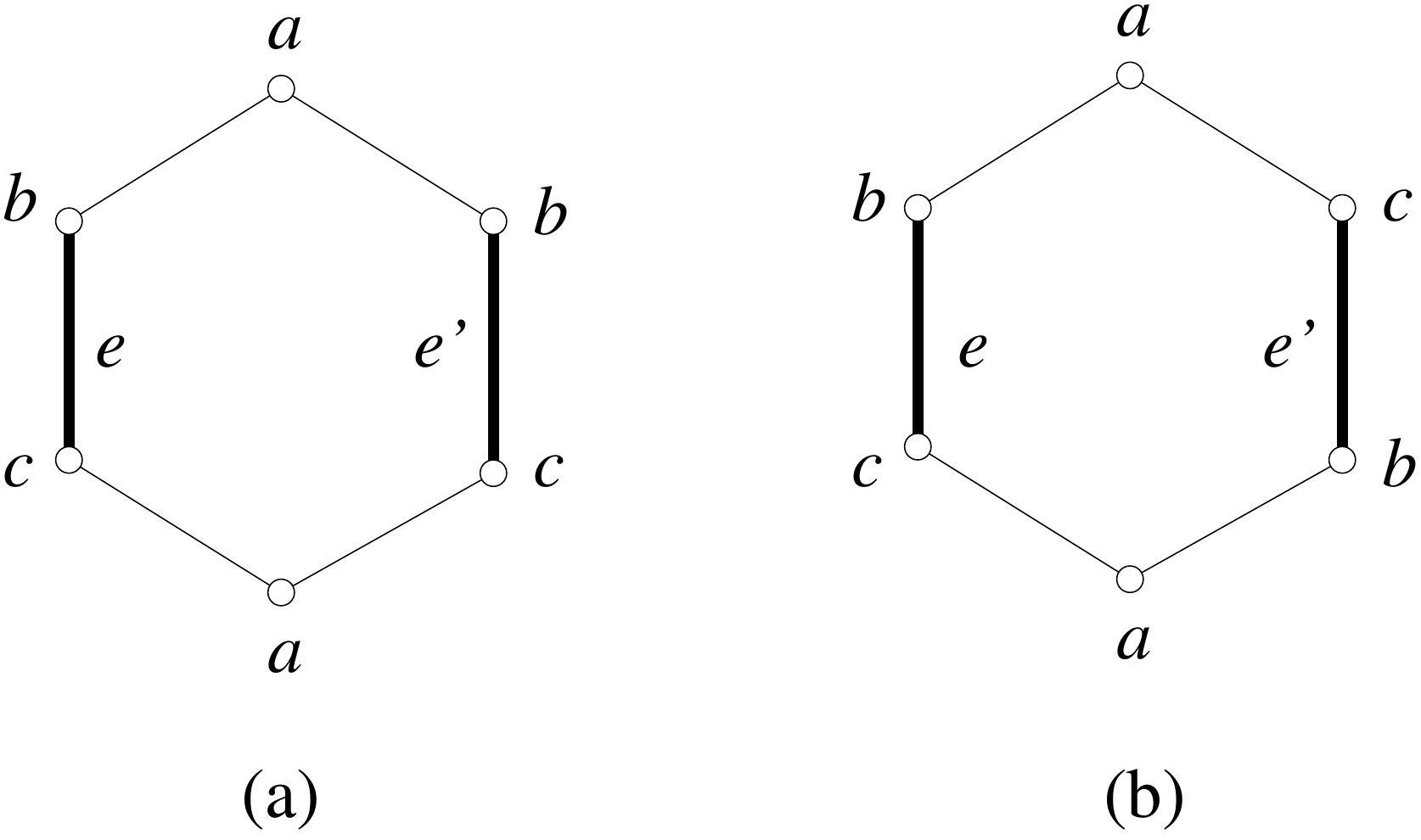}
        \hspace*{0.5cm}
        \caption{Color Patterns and nonadjacent color equivalent edges $e$ and $e'$. \label{fig:colorpatterns}}
\end{figure}

\begin{enumerate}
\item[]  Case 1.  There is a vertex in the cycle whose neighbors in the cycle have the same color.  Up to relabeling, we can assume this vertex has color $a$ (say in state $a_0$) and the two adjacent vertices have color $b$.  The states for the remaining vertices of the cycle are $a_1$, $c_0$, and $c_1$.  Now, consider the vertices adjacent to $b_0$ and $b_1$ other than $a_0$.  These vertices must be $c_0$ and $c_1$ (otherwise, the two states of $c$ would be adjacent in the cycle).  This color pattern is shown in Figure \ref{fig:colorpatterns}(a). 
\item[]
\item[] Case 2.  No vertex in the cycle is adjacent to two vertices of the same color.  Then the two neighbors of a vertex with color $a$ must have colors $b$ and $c$.  Then the vertex following $b$ in the cycle must have color $c$ (otherwise vertex $b$ is adjacent to two vertices of the same color).  By working this way around the cycle, the only color pattern possible is as shown in Figure \ref{fig:colorpatterns}(b).  
\end{enumerate}

Note that both color patterns contain a pair of nonadjacent and color equivalent edges (edges $e$ and $e'$ in Figure \ref{fig:colorpatterns}).

\begin{figure}[h!]
        \centering
        \includegraphics[scale=0.35]{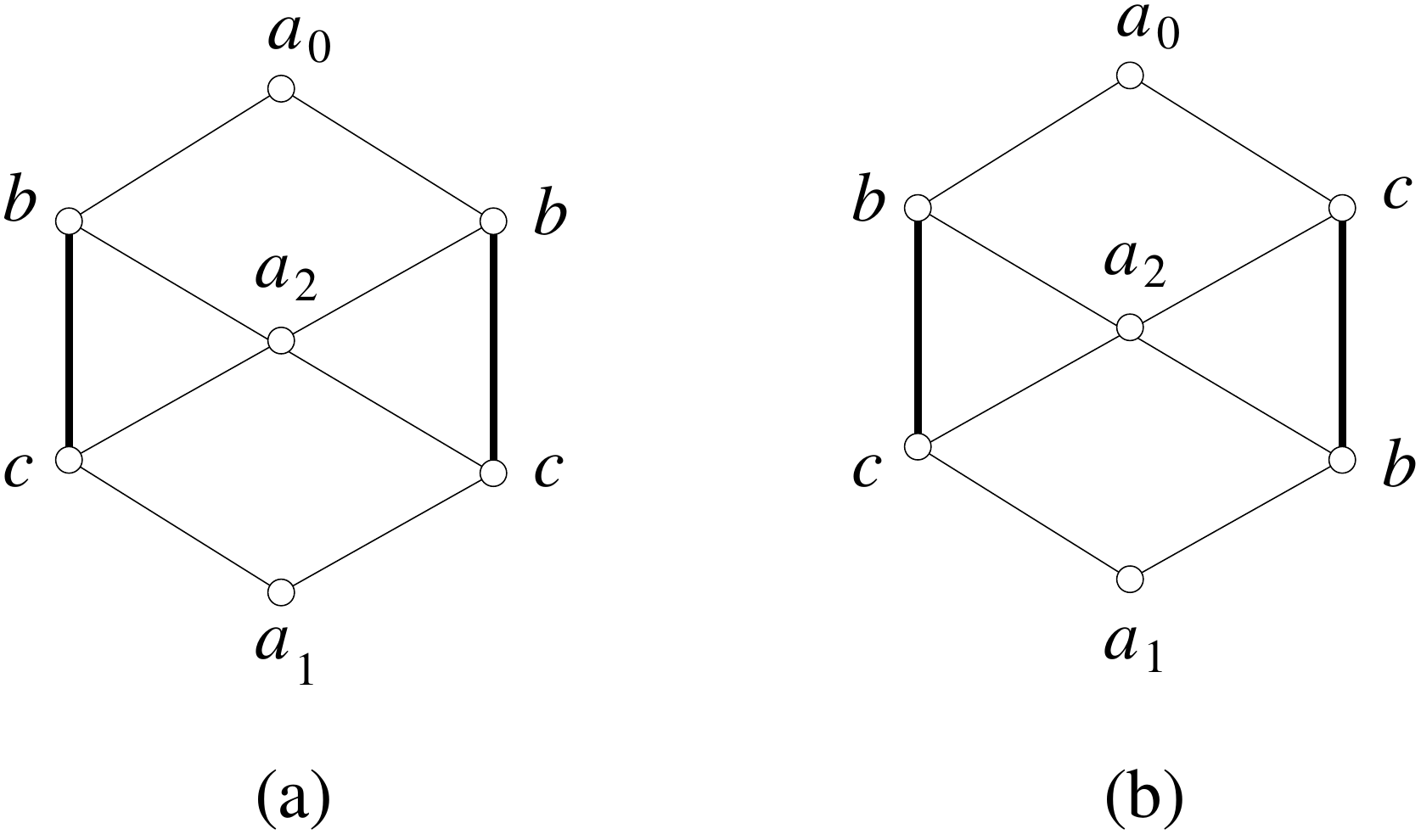}
        \hspace*{0.5cm}
        \caption{Induced Color Patterns \label{fig:inducedgraph}}
\end{figure}

Consider this pair of nonadjacent and color equivalent edges $e$ and $e'$. Without loss of generality, assume that the endpoints of these edges have colors $b$ and $c$.  Let $s$ be the row witness for $e$ and $s'$ be the row witness for $e'$.  Since cycle $C$ is chordless, the state in character $a$ of row $s$ cannot be $a_0$ or $a_1$.  
Similarly, the state in character $a$ of row $s'$ cannot be $a_0$ or $a_1$.  Since $a_2$ is the only remaining state of character $a$, both $s$ and $s'$ must contain $a_2$.  This implies that the partition intersection graph $G[a,b,c]$ must induce one of the two color patterns in Figure \ref{fig:inducedgraph}.

\begin{figure}[h!]
        \centering
        \includegraphics[scale=0.35]{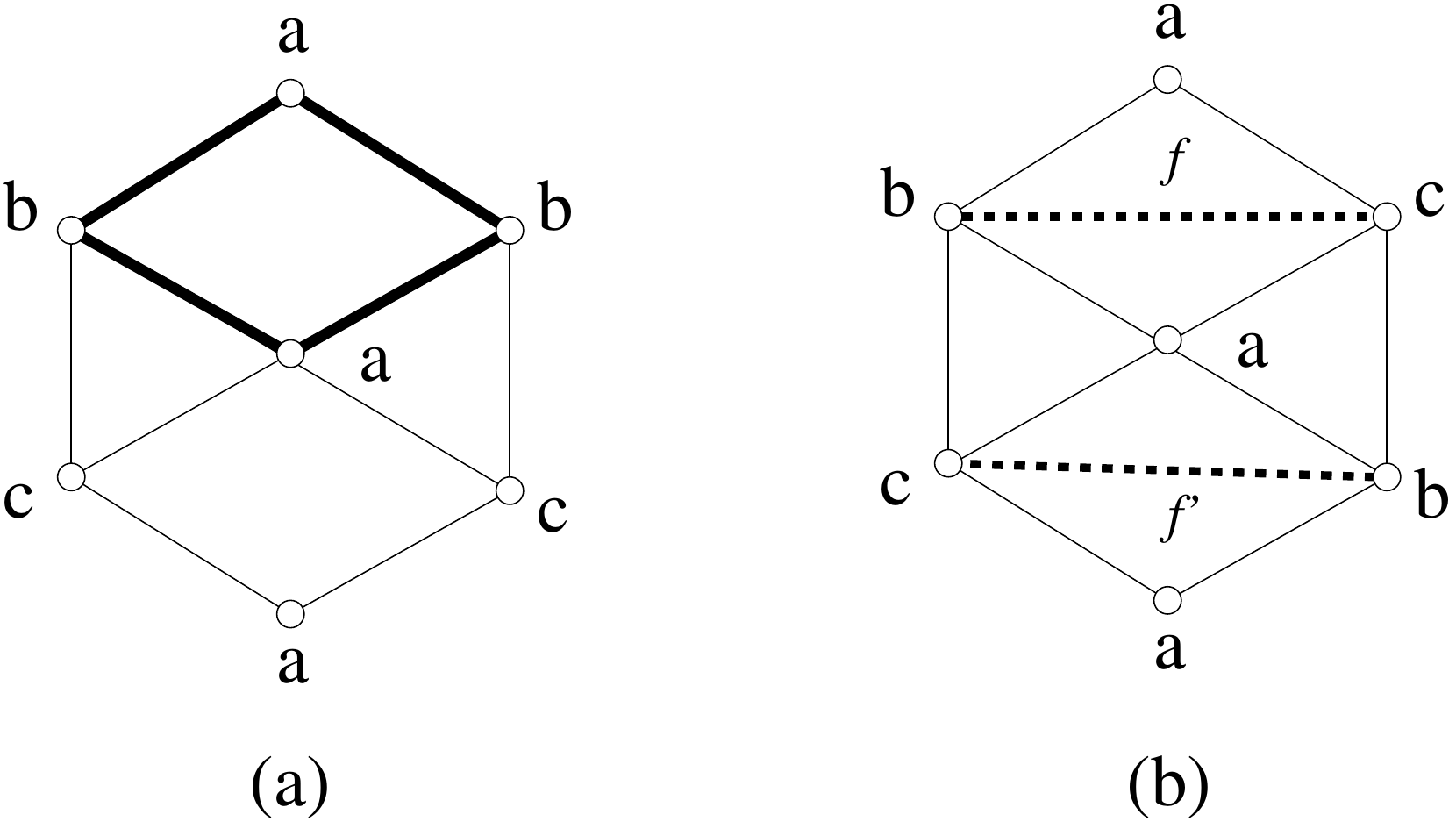}
        \hspace*{0.5cm}
        \caption{(a) Induced cycle of length four on two colors; (b) Forced Edges $f$ and $f'$ \label{fig:forcededges}}
\end{figure}

In the case illustrated in Figure \ref{fig:inducedgraph}(a), there is a cycle 
on four vertices induced by the two characters $a$ and $b$ (see Figure 
\ref{fig:forcededges}(a)), implying $G[a,b,c]$ is not properly triangulatable.  In the 
case illustrated in Figure \ref{fig:inducedgraph}(b), there are two edge-disjoint 
cycles of length four with color pattern $a, b, a, c$.  Since edges in a proper 
triangulation cannot connect vertices of the same color, any proper triangulation of 
$G$ must contain the two edges $f$ and $f'$ connecting vertices of color $b$ and $c$ 
(see Figure \ref{fig:forcededges}(b)).  However, this induces a cycle of length four on 
the states of $b$ and $c$, which does not have a proper triangulation.  This again shows that $G[a,b,c]$ is not properly triangulatable.

Since all of these cases result in contradictions, it follows that there exists 
a color that appears exactly once in $C$.  \qed

\begin{figure}[h!]
        \centering
        \includegraphics[scale=0.35]{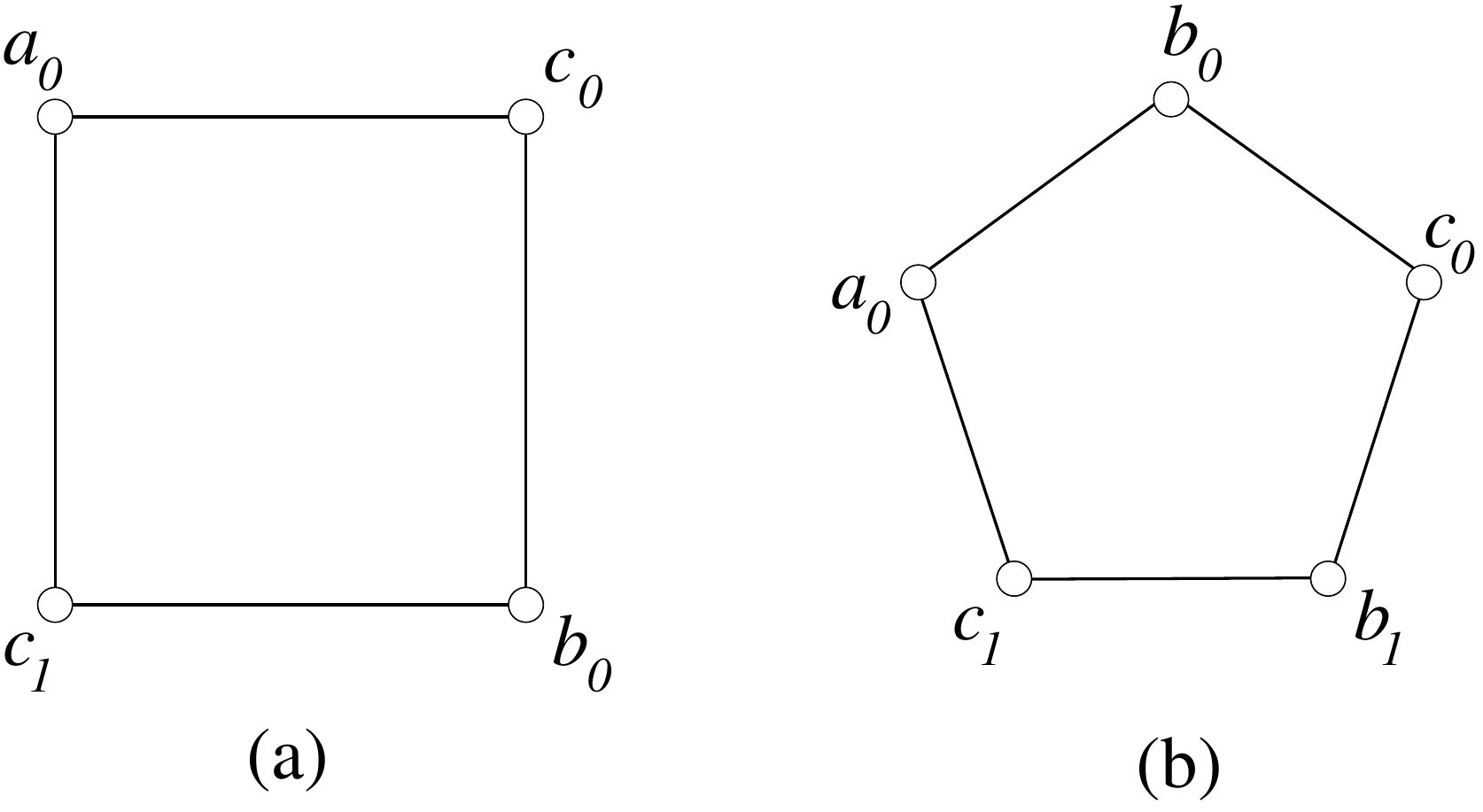}
        \caption{The only possible chordless cycles in $G[a,b,c]$: (a) characters $a$ and $b$ appear uniquely while character $c$ appears twice; (b) character $a$ appears uniquely while characters $b$ and $c$ each appear twice.\label{fig:possible_chordless}}
\end{figure}

Lemmas \ref{lemma:nothreestates} and \ref{lemma:unique} show that if $C$ is a chordless cycle in a properly triangulatable graph $G[a,b,c]$, then no color can appear in all three states and one color appears uniquely.  This leaves two possibilities for chordless cycles in $G[a,b,c]$ (see Figure \ref{fig:possible_chordless}):

\begin{quote}
\begin{itemize}
\item[$\bullet$] a chordless four cycle, with two colors appearing uniquely and the remaining color appearing twice 
\item[$\bullet$] a chordless five cycle, with one color appearing uniquely and the other two colors each appearing twice
\end{itemize}
\end{quote}

In the next lemma, we show that if $G[a,b,c]$ is properly triangulatable, the second case cannot occur, i.e., $G[a,b,c]$ cannot contain a chordless five cycle.

\begin{lemma}\label{lemma:nofivecycle} Let $S$ be a set of input species on three characters $a, b,$ and $c$ with at most three states per character.   If the partition intersection graph $G[a,b,c]$ is properly triangulatable, then $G[a,b,c]$ cannot contain chordless cycles of length five or greater. 
\end{lemma}

\noindent {\bf Proof.} Lemmas \ref{lemma:nothreestates} and \ref{lemma:unique} together show that $G[a,b,c]$ cannot contain chordless cycles of length six or greater, so it remains to show that $G[a,b,c]$ cannot contain chordless cycles of length equal to five.  

Suppose $C$ is a chordless cycle in $G[a,b,c]$ of length five; without loss of generality, let $a$ be the color appearing exactly once in $C$ (say in state $a_0$), let $b_0, b_1$ be the two states of $b$ in $C$, and let $c_0, c_1$ be the two states of $c$ in $C$.  Up to relabeling of the states, the cycle is as shown in Figure \ref{fig:possible_chordless}(b). 

\begin{figure}[h!]
        \centering
        \includegraphics[scale=0.35]{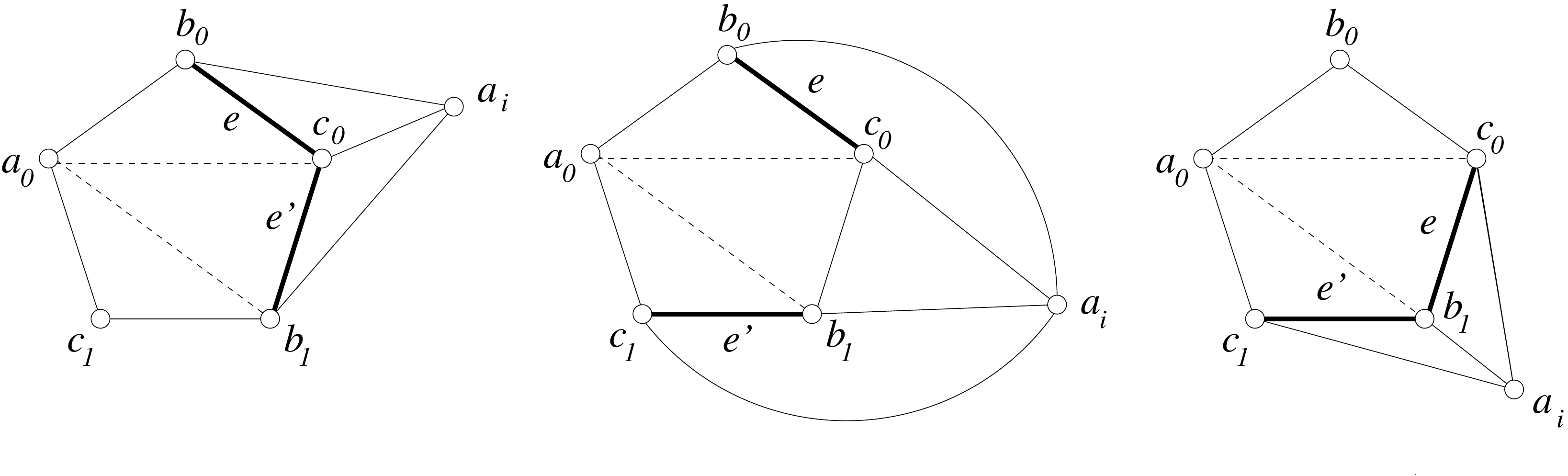}
        \caption{Edges $e$ and $e'$ are both witnessed by state $a_i$.\label{fig:forbidden_5cycle2}}
\end{figure}

Now, any proper triangulation of $G[a,b,c]$ must triangulate cycle $C$ by edges $(a_0, c_0)$ and $(a_0, b_1)$ shown in Figure \ref{fig:forbidden_5cycle2} (since the only other edge between nonadjacent vertices of different colors is $(b_0, c_1)$, which would create a non-triangulatable four cycle on the two colors $b$ and $c$. 
%$(b_0, c_0)$, $(c_0, b_1)$, $(b_1, c_1)$, $(c_1, b_0)$).

The row witnesses for edges $(b_0, c_0)$, $(c_0, b_1)$, and $(c_1, b_1)$ must contain a state in color $a$ that is one of $a_1$ or  $a_2$ (otherwise, $a_0$ would have an edge to a non-adjacent vertex in cycle $C$, implying $C$ is not  chordless).  Since there are three edges and two possible witness states in color $a$, there are two edges among $(b_0, c_0)$, $(c_0, b_1)$, $(c_1, b_1)$ that share a witness $a_i$.  We denote these two edges by $e$ and $e'$; as shown in Figure \ref{fig:forbidden_5cycle2}, there are three ways to choose $e$ and $e'$.  

%For any choice of the two edges sharing a witness, either 
%\begin{enumerate}
%\item both $b_0$ and $b_1$ are adjacent to $a_i$ or 
%\item both $c_0$ and $c_1$ are adjacent to $a_i$. 
%\end{enumerate}

\begin{figure}[h!]
        \centering
        \includegraphics[scale=0.35]{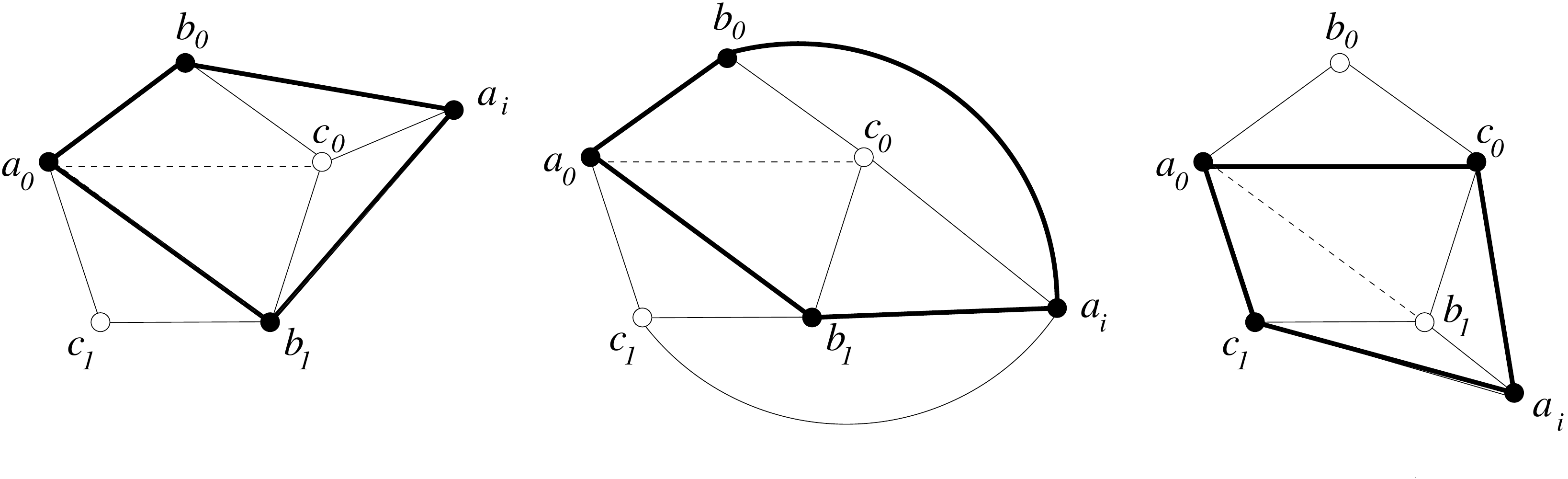}
        \caption{Forced cycles of length four on two colors.\label{fig:forbidden_5cycle3}}
\end{figure}

Figure \ref{fig:forbidden_5cycle3} shows that all three cases induce a four cycle on two colors, a contradiction since $G[a,b,c]$ is properly triangulatable.  Therefore, $G[a,b,c]$ cannot contain a chordless 5-cycle.  \qed

\begin{lemma}\label{lemma:unique2} Let $S$ be a set of input species on three 
characters $a, b,$ and $c$ with at most three states per character.   If the 
partition intersection graph $G[a,b,c]$ is properly triangulatable, then every chordless cycle in $G[a,b,c]$ is \emph{uniquely} triangulatable. \end{lemma}

\begin{figure}[h!]
        \centering
        \includegraphics[scale=0.35]{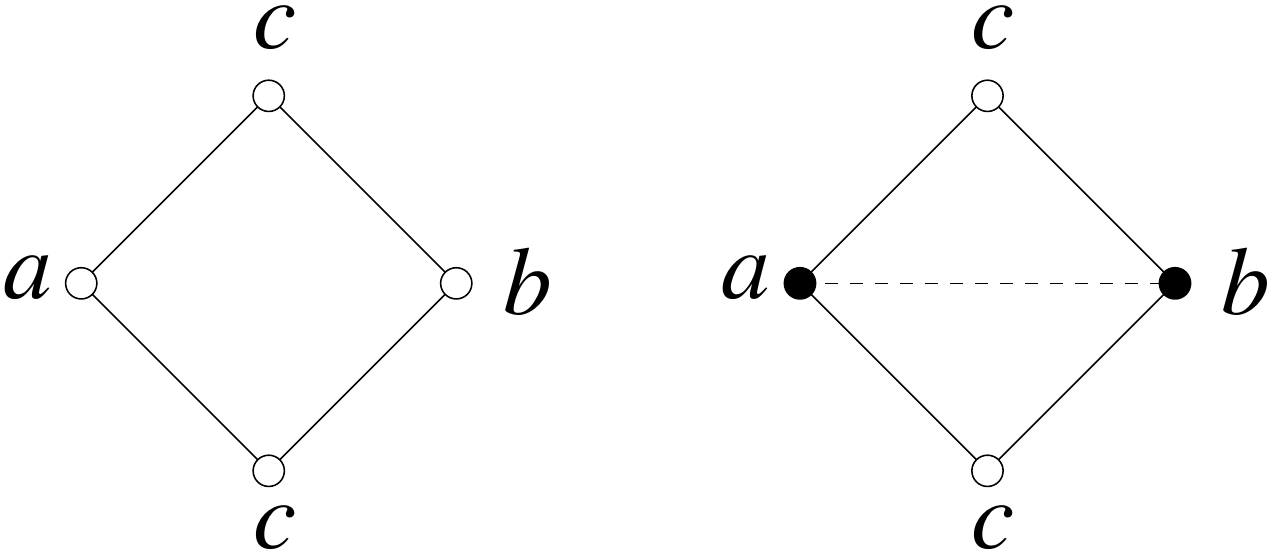}
        \caption{Color pattern for chordless cycle $C$. \label{fig:twounique}}
\end{figure}

\Proof By Lemma \ref{lemma:nofivecycle}, if $C$ is a chordless cycle in $G[a,b,c]$, then $C$ must be a four cycle with the color pattern shown in Figure \ref{fig:twounique} (up to relabeling of the colors).  Then $C$ is uniquely triangulatable by adding the edge between the two colors appearing uniquely (in Figure \ref{fig:twounique}, these are colors $a$ and $b$). 

\qed

For any three colors $a,b,c$, Lemma \ref{lemma:unique2} gives a simple algorithm to properly triangulate $G[a,b,c]$: for each chordless cycle $C$ in $G[a,b,c]$, check that $C$ is a four cycle with two nonadjacent vertices having colors that appear exactly once in $C$ and add an edge between these two vertices.

\section{The 3-character test} \label{sec:generalm}

\subsection{Triangulating Triples of Characters}

We now consider the case of trinary input sequences $S$ on $m$ characters (for $m$ greater or equal to 4).  Our goal is to prove that the existence of proper triangulations for all subsets of three characters at a time is a sufficient condition to guarantee existence of a proper triangulation for \emph{all} $m$ characters.

By Lemma \ref{lemma:unique}, if a set of three characters $\chi^i, \chi^j, \chi^k$ is properly triangulatable, then there is a unique set of edges $F(\chi^i, \chi^j, \chi^k)$ that must be added to triangulate the chordless cycles in $G[\chi^i, \chi^j, \chi^k]$. Construct a new graph $G'(S)$ on the same vertices as $G(S)$ with edge set $E(G(S)) \cup \{ \cup_{1 \leq i < j < k \leq m} F(\chi^i, \chi^j, \chi^k) \}$.  $G'(S)$ is the partition intersection graph $G(S)$ together with all of the additional edges used to properly triangulate chordless cycles in $G[\chi^i, \chi^j, \chi^k]$ ($1 \leq i < j < k \leq m$).  In $G'(S)$, edges from the partition intersection graph $G(S)$ are called $E$-edges and edges that have been added as triangulation edges for some triple of columns are called $F$-edges.  We call a cycle consisting only of $E$-edges an \emph{$E$-cycle}.  

\begin{example} \label{ex:chordless4} \emph{Consider input set $S$ and the corresponding partition intersection graph $G(S)$ in Figure \ref{fig:chordless4cycle}.  Each triple of characters in $S$ induces a chordal graph while the entire partition intersection graph $G(S)$ contains a chordless cycle of length four.  Since each triple of characters induces a chordal graph, no $F$-edges are added and $G(S) = G'(S)$. }

\begin{figure}[h!]
        \centering
        \includegraphics[scale=0.35]{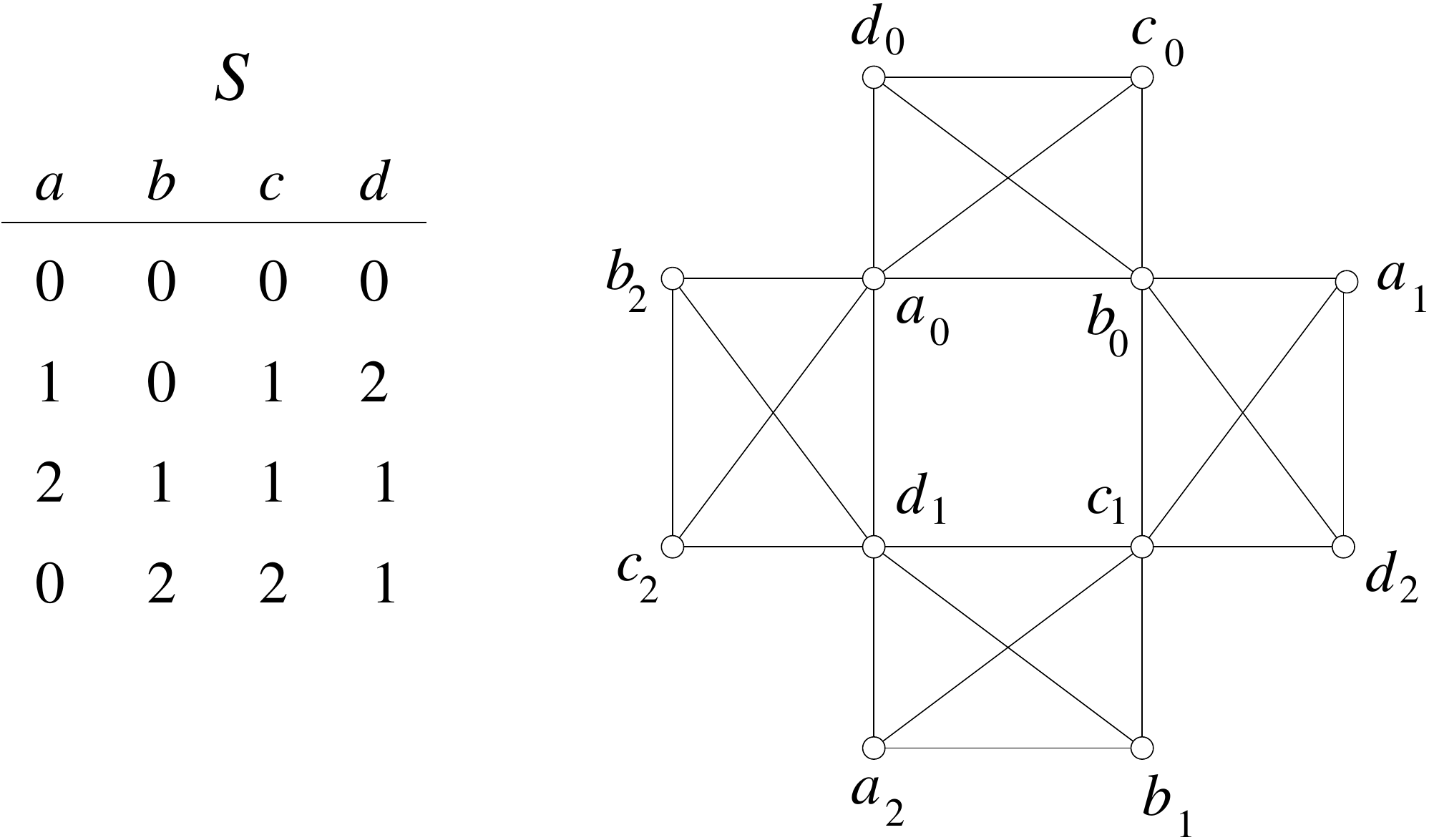}
        \caption{Partition intersection graph $G'(S)$ contains a chordless four cycle.\label{fig:chordless4cycle}}
\end{figure}

\end{example}

As Example \ref{ex:chordless4} illustrates, the addition of $F$-edges alone may not be sufficient to triangulate the entire partition intersection graph.  We now turn to the problem of triangulating the remaining chordless cycles in $G'(S)$.

Consider any $E$-cycle $C$ that is chordless in $G'(S)$ satisfying the properties

\begin{quote}
\begin{enumerate}
\item $C$ has length equal to four 
\item all colors of $C$ are distinct
\end{enumerate}
\end{quote}

For every such chordless cycle, add the chords between the two pairs of nonadjacent vertices in $C$ (note that these are legal edges).  Call this set of edges $F'$-edges and let $G''(S)$ denote the graph $G'(S)$ with the addition of $F'$-edges.  Note that the sets of $E$-edges, $F$-edges, and $F'$-edges are pairwise disjoint; we call the set of $F$ and $F'$-edges \emph{non-$E$ edges}.

\comment{
 any pair of vertices adjacent by an $F$-edge must be nonadjacent in the original partition intersection graph, implying the vertices cannot be adjacent by an $E$-edge.   
We first study some properties of $E$-edges and $F$-edges. Suppose $C$ is a chordless $E$-cycle in $G'(S)$ (i.e., $C$ is a chordless cycle in the original partition intersection graph $G(S)$ that remains chordless upon adding all $F$-edges $\cup_{1 \leq i < j < k \leq m} F(\chi^i, \chi^j, \chi^k)$).  Since all chordless $E$-cycles on three colors in $G(S)$ have been properly triangulated in $G'(S)$, $C$ must contain four or more colors.
}

We begin by investigating structural properties of cycles in $G'(S)$ and $G''(S)$ containing at least one $F$-edge or $F'$-edge.  Let $C$ be a cycle in $G'(S)$ or $G''(S)$ containing an edge $f$ that is an $F$-edge or $F'$-edge (without loss of generality, let $f=(a_0, b_0)$).  This edge must be added due to an $E$-cycle $D$ containing $a_0, b_0$ and two other vertices $w$ and $z$ as shown in Figure \ref{fig:fchord1}(a) (note that $w$ and $z$ cannot have color $a$ or $b$).  If $f$ is an $F$-edge, then $w$ and $z$ have the same color and therefore cannot be adjacent in $G'(S)$.  If $f$ is an $F'$-edge, then since $D$ is a chordless $E$-cycle in $G'(S)$, $w$ and $z$ are again nonadjacent in $G'(S)$.  The cycle $C$ created by edge $f$ is shown in Figure \ref{fig:fchord1}(b).

\begin{figure}[h!]
        \centering
        \includegraphics[scale=0.35]{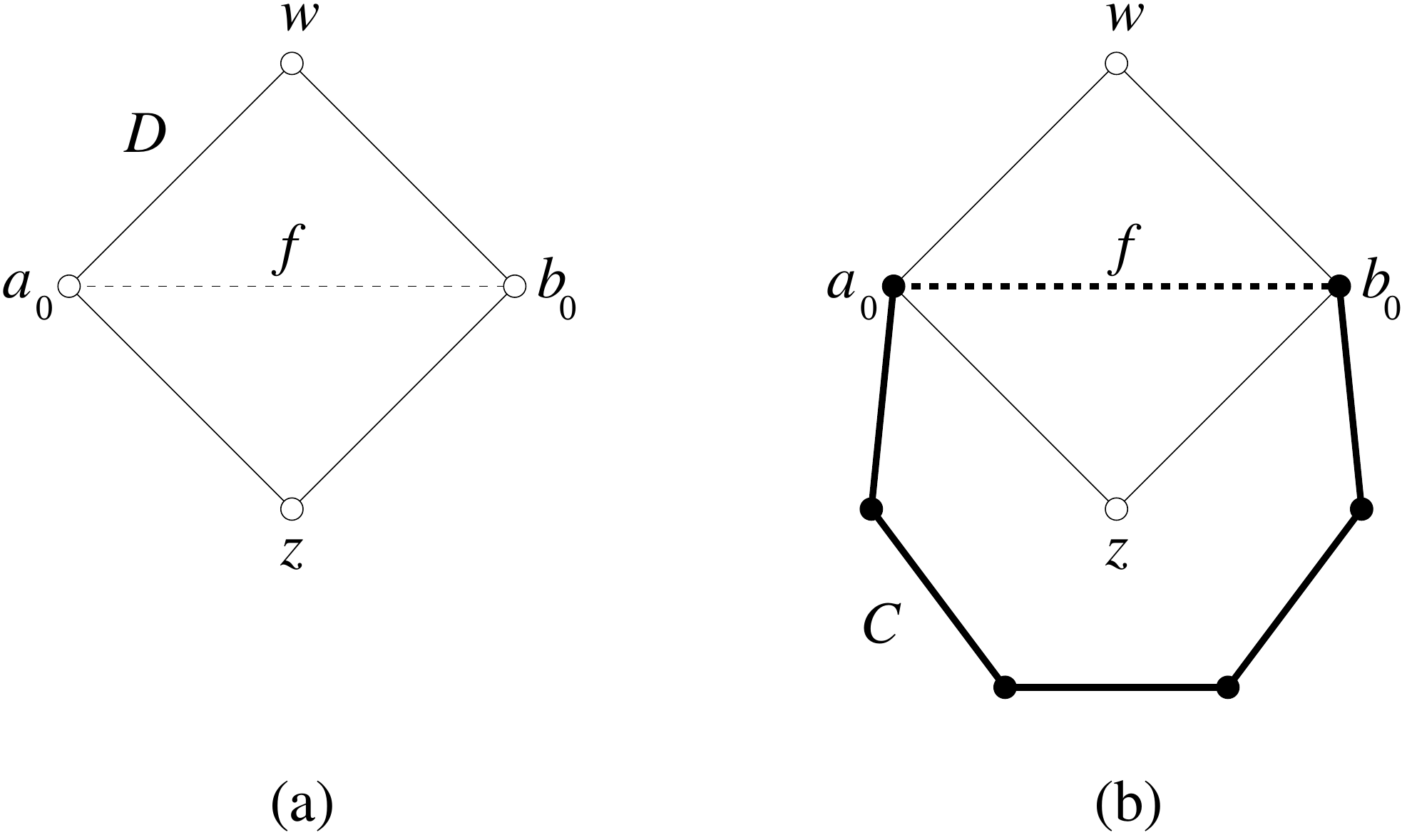}
        \caption{(a) Chordless cycle $D$ (b) edge $f = (a_0, b_0)$ creates cycle $C$ (shown in bold). \label{fig:fchord1}}
\end{figure}

\begin{figure}[h!]
        \centering
        \includegraphics[scale=0.35]{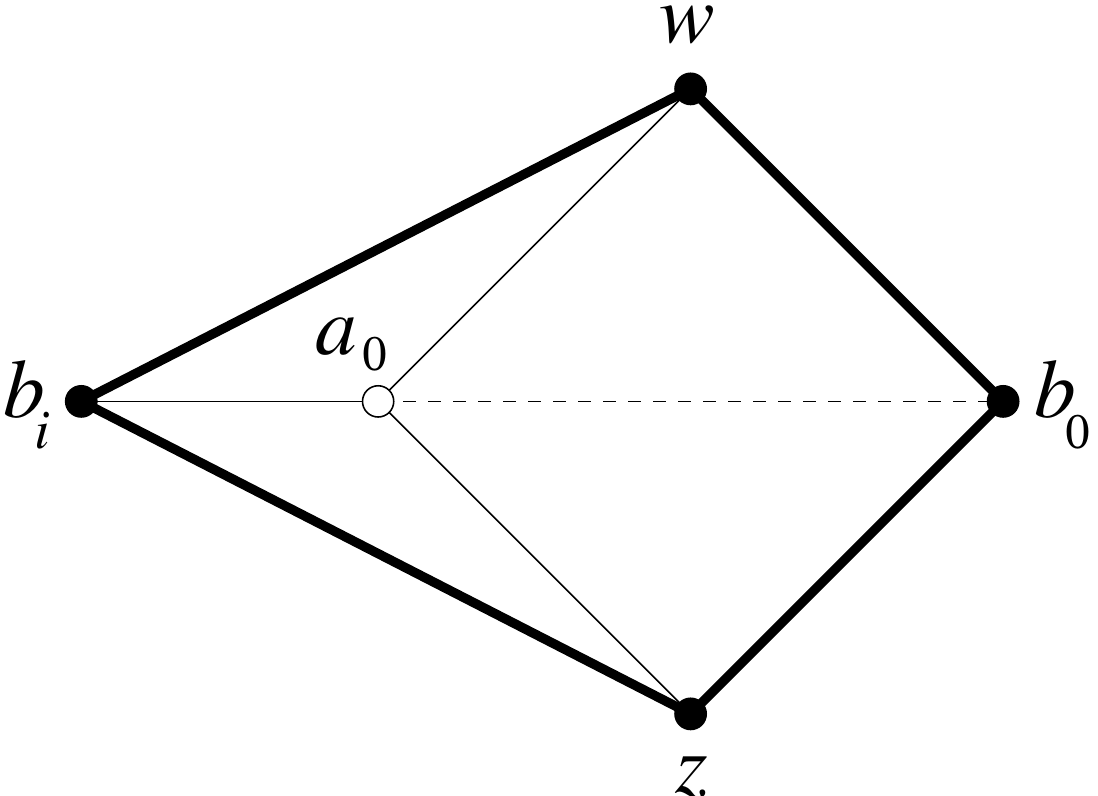}
        \caption{If the row witnesses for $(a_0, w)$ and $(a_0, z)$ share a state of $b$, there is a chordless $E$-cycle of length four on at most three colors.\label{fig:fchord2}}
\end{figure}

Since $D$ is an $E$-cycle, each edge in $D$ has a row witness.  Consider first the row witnesses for edges $(a_0, w)$ and $(a_0, z)$.  These row witnesses must contain a state of $b$ other than $b_0$ (since $a_0$ and $b_0$ are not connected by an $E$-edge).  If both row witnesses share the same state $b_i$ of $b$, then the cycle $(b_i, w)$, $(w, b_0)$, $(b_0, z)$, $(z, b_i)$ is a chordless $E$-cycle on at most three colors in $G'(S)$  as shown in Figure \ref{fig:fchord2} (as argued above, $w$ and $z$ are nonadjacent in $G'(S)$).  However, all chordless $E$-cycles on at most three colors have been triangulated in $G'(S)$, a contradiction.

\begin{figure}[h!]
        \centering
        \includegraphics[scale=0.35]{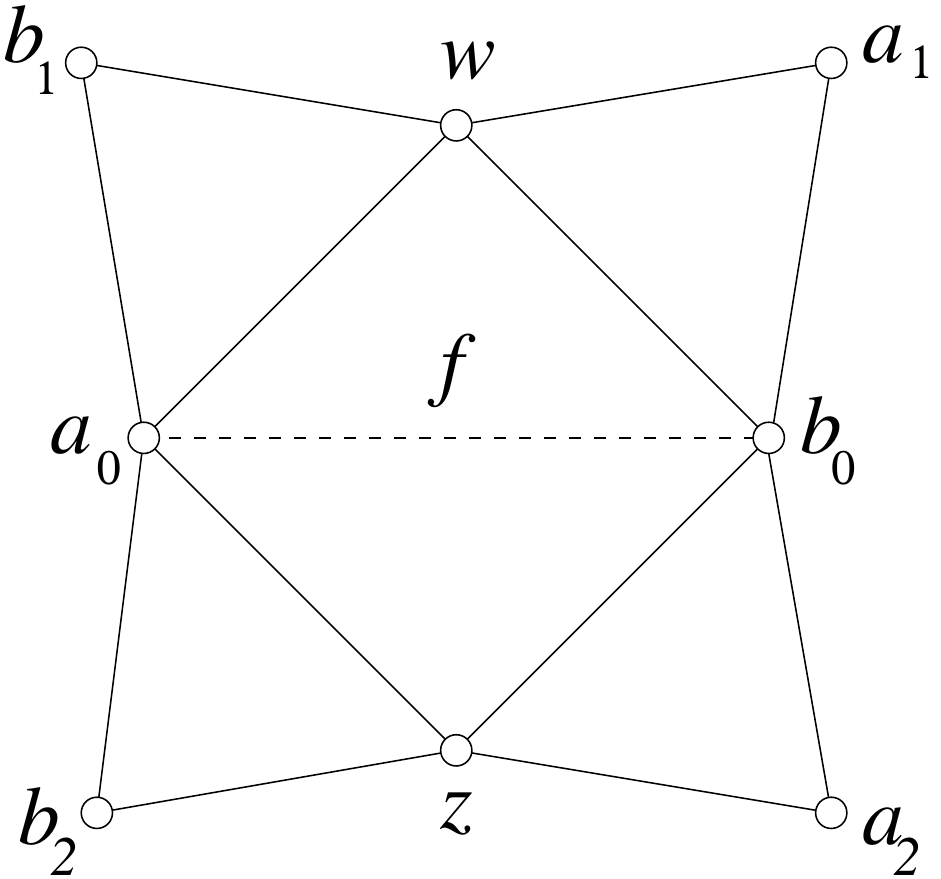}
        \caption{Pattern of forced witnesses for edges in $D$.\label{fig:fchord3}}
\end{figure}

Therefore, the row witnesses for $(a_0, w)$ and $(a_0, z)$ cannot share the same state of $b$.  Similarly, the row witnesses for $(b_0, w)$ and $(b_0, z)$ cannot share the same state of $a$. This implies the following situation, up to relabeling of the states, illustrated in Figure \ref{fig:fchord3}. 

In particular, the following two conditions must be satisfied.

\begin{eqnarray*}
(1) \ a_0 \text{ is adjacent to both } b_1 \text{ and } b_2 &\\
\vspace{-0.75ex}  & \hspace{.5in} (*)\\
(2) \vspace{-0.75ex} \ b_0 \text{ is adjacent to both } a_1 \text{ and } a_2 & 
\end{eqnarray*}

We use this structure to prove a sequence of lemmas eliminating the possibilities for chordless cycles in graph $G'(S)$.  This sequence of lemmas will show $G''(S)$ cannot contain a chordless cycle with exactly one non-$E$ edge (Lemmas \ref{lemma:onefedge} and \ref{lemma:onef'edge}), a chordless cycle with two or more non-$E$ edges (Lemma \ref{lemma:nonE}), or a chordless $E$-cycle (Corollary \ref{cor:ecycle}).  

\begin{lemma}\label{lemma:onefedge} $G'(S)$ cannot contain a chordless cycle with exactly one $F$-edge. \end{lemma}

\begin{figure}[h!]
        \centering
        \includegraphics[scale=0.35]{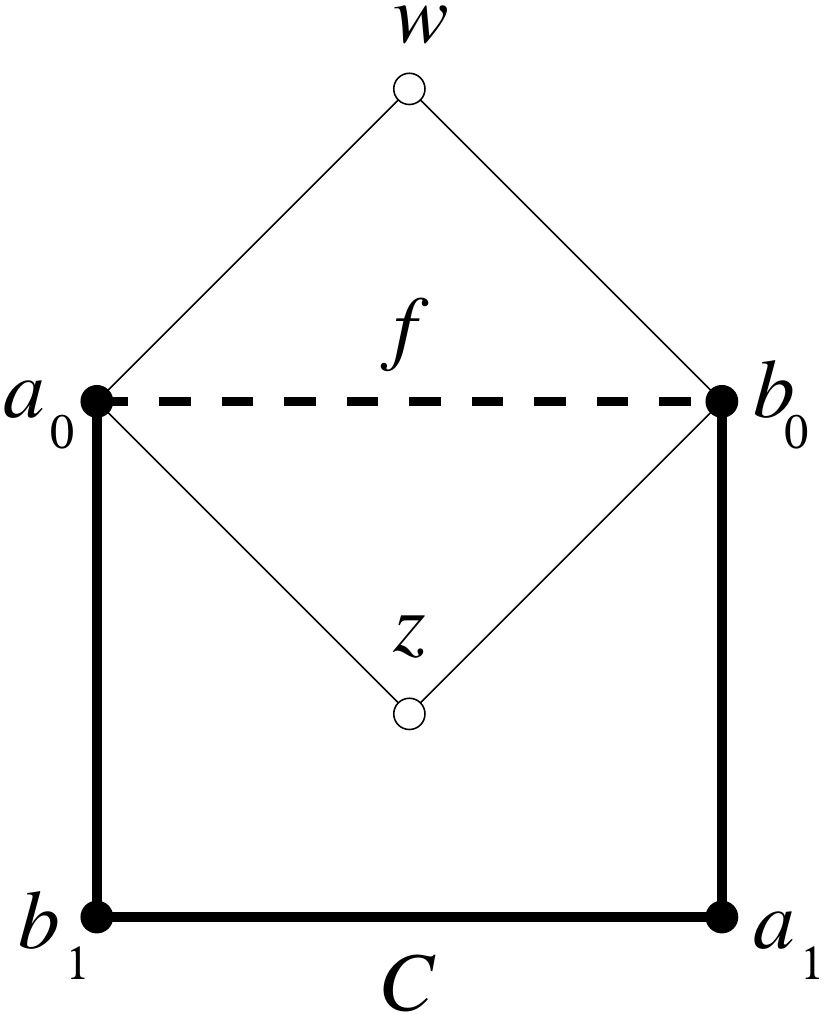}
        \caption{Chordless cycle $C$ on two colors with exactly one $F$-edge. \label{fig:oneedge1}}
\end{figure}

\Proof Suppose that $C$ is a chordless cycle in $G'(S)$ with exactly one $F$-edge, say $f=(a_0, b_0)$.  Edge $(a_0, b_0)$ must have been added due to a chordless $E$-cycle $D$ on three colors as shown in Figure \ref{fig:fchord3}, where $w$ and $z$ are states of the same color.  Note that edge $(a_0, b_0)$ is a forced $F$-edge that creates cycle $C$ (see Figure \ref{fig:oneedge1}). If $C$ contains only the two colors $a$ and $b$, the partition intersection graph on the three colors $a, b$, and the shared color of $w$ and $z$ is not properly triangulatable, a contradiction.  

\begin{figure}[h!]
        \centering
        \includegraphics[scale=0.30]{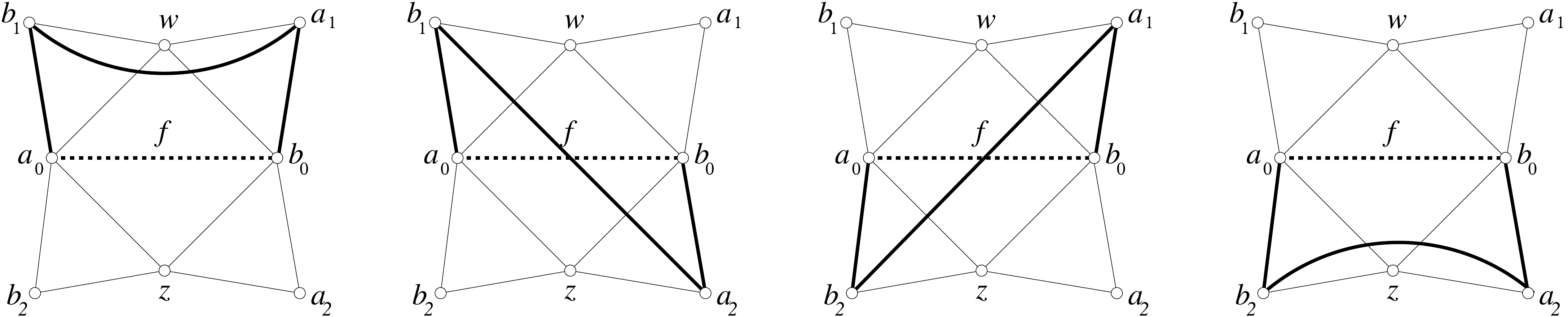}
        \caption{If any of $(b_1, a_1)$, $(b_1, a_2)$, $(b_2, a_1)$, or $(b_2, a_2)$ are $E$-edges, there is a chordless four cycle in $G'(S)$ on two colors with exactly one $F$-edge.\label{fig:fchord3i}}
\end{figure}

This implies any cycle $C$ in $G'(S)$ with exactly one $F$-edge must contain three or more colors.  As shown in Figure \ref{fig:fchord3i}, if any of the edges $(b_1, a_1)$, $(b_1, a_2)$, $(b_2, a_1)$, and $(b_2, a_2)$ are present, there would be a chordless cycle on two colors with exactly one $F$-edge, which we have argued cannot occur.  It follows that $a_1$ and $a_2$ are nonadjacent to $b_1$ and $b_2$ by $E$-edges.
 
Since $a_1$ is nonadjacent to $b_1$ or $b_2$ by $E$-edges, any row that contains $a_1$ must contain state $b_0$ in character $b$.   We call this condition (A1).  By a similar argument, the following conditions must be satisfied:

\medskip

\indent (A2) any row that contains $a_2$ must contain state $b_0$ in character $b$.\\
\indent (B1) any row that contains $b_1$ must contain state $a_0$ in character $a$.\\
\indent (B2) any row that contains $b_2$ must contain state $a_0$ in character $a$.

\medskip

\begin{figure}[h!]
        \centering
        \includegraphics[scale=0.35]{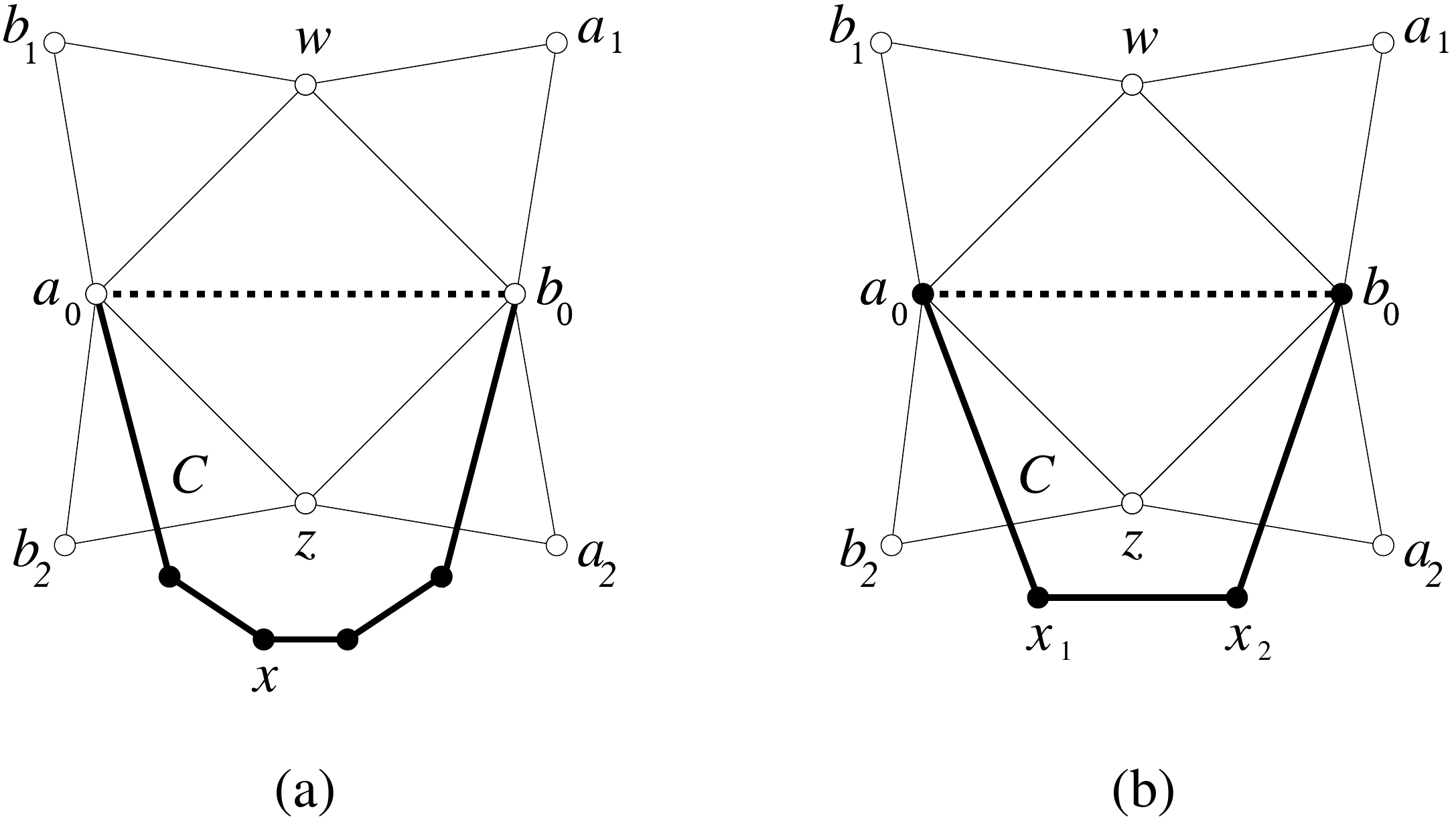}
        \caption{Vertices on cycle $C$. \label{fig:fchord7i}}
\end{figure}

Now, let $x$ be a vertex in $C \backslash \{a_0, b_0 \}$ and consider the state of character $a$ in any row that witnesses $x$ (see Figure \ref{fig:fchord7i}(a)).  If this state is $a_0$, then $x$ is adjacent to $a_0$ by an $E$-edge.  Otherwise, if this state is either $a_1$ or $a_2$, then this row witness for $x$ must contain state $b_0$ by (A1) and (A2).  Since $C$ is a chordless cycle, at most one vertex on $C \backslash \{a_0, b_0\}$ can be adjacent to each of $a_0$ and $b_0$.  This shows there can be at most two such vertices $x_1$ and $x_2$ in $C \backslash \{ a_0, b_0\}$, one of which is adjacent to $a_0$ and the other which is adjacent to $b_0$ (moreover, these are adjacencies by $E$-edges).  Therefore, $C$ has length equal to four formed by edges $(a_0, x_1)$, $(x_1, x_2)$, $(x_2, b_0)$, and $(b_0, a_0)$ (see Figure \ref{fig:fchord7i}(b)).

Edge $(x_1, x_2)$ is an $E$-edge since $f=(a_0, b_0)$ is the unique $F$-edge in $C$ by assumption.  Furthermore, at least one of $x_1$ or $x_2$ has color different from $a$ and $b$ since $C$ has three or more colors.   Without loss of generality, assume vertex $x_1$ has color different from $a$ and $b$. The color of $x_2$ is different from $b$ since $x_2$ and $b_0$ are adjacent, implying edge $(x_1, x_2)$ must have a witness in character $b$.   If this witness is $b_0$, then $x_1$ and $b_0$ are adjacent by an $E$-edge, a contradiction to the chordlessness of cycle $C$.  If this witness is either $b_1$ or $b_2$, then (B1) or (B2) imply that $x_2$ and $a_0$ are adjacent by an $E$-edge, again a contradiction to the chordlessness of cycle $C$.  

This  concludes the proof of Lemma \ref{lemma:onefedge}. \qed

A similar proof shows that the lemma can be extended to the graph $G''(S)$.

\begin{lemma}\label{lemma:onef'edge} $G''(S)$ cannot conatin a chordless cycle with exactly one non-$E$ edge.   \end{lemma}

\Proof  Suppose that $C$ is a chordless cycle in $G''(S)$ with exactly one non-$E$ edge, say $f=(a_0, b_0)$.  If $f$ is an $F$-edge, then $C$ would be a chordless cycle in $G'(S)$ with exactly one $F$-edge, contradicting Lemma \ref{lemma:onefedge}.  Therefore $f$ is an $F'$-edge that is added due to chordless cycle $D$ as shown in Figure \ref{fig:fchord3} (with $w$ and $z$ different colors). 

\noindent {\bf Case I.} $C$ contains only the two colors $a$ and $b$.  Since the graph in Figure \ref{fig:fchord3} contains all the states of characters $a$ and $b$, $C$ must also contain one of the edges $(b_1, a_1)$, $(b_1, a_2)$, $(b_2, a_1)$, or $(b_2, a_2)$ as an $E$-edge and we have the following cases. 

\begin{itemize}
\item[] Case I(i). $C$ contains edge $(b_2, a_2)$.  This results in an $E$-cycle of length five on at most three colors as shown in Figure \ref{fig:oneedge2}(a).  Such a cycle cannot be chordless in $G[a,b,w]$ by Lemma \ref{lemma:nofivecycle}.  Therefore, vertex $w$ must be adjacent to one of $b_2$ or $a_2$ by an $E$-edge.  This creates a chordless $E$-cycle in $G'(S)$ of length four on three colors; either cycle $(b_2, w)$, $(w, b_0)$, $(b_0, z)$, $(z, b_2)$  shown in Figure \ref{fig:oneedge2}(b) or cycle $(a_0, w)$, $(w, a_2)$, $(a_2, z)$, $(z, a_0)$ shown in Figure \ref{fig:oneedge2}(c)  (note that $w$ and $z$ are nonadjacent in $G'(S)$ since cycle $D$ is chordless).  This is a contradiction since all cycles on at most three colors are triangulated in $G'(S)$.

\begin{figure}[h!]
        \centering
        \includegraphics[scale=0.35]{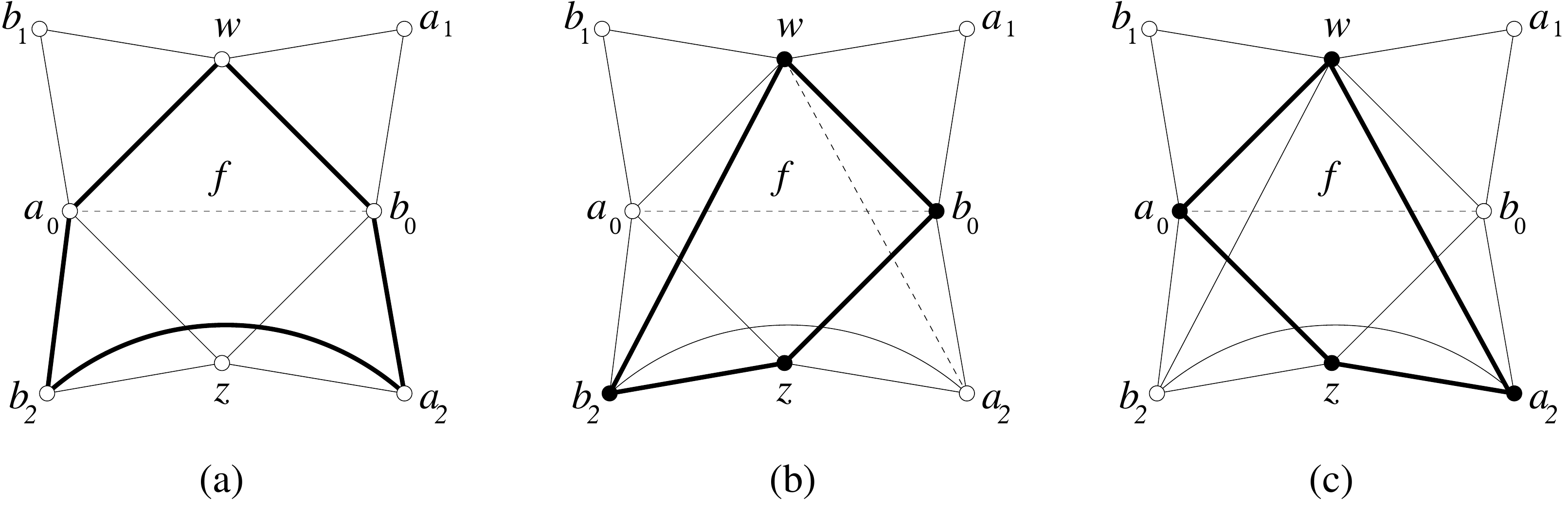}
        \caption{(a) Edge $(b_2, a_2)$ gives a five cycle $C$ on at most three colors $C$ (b),(c) Chordless cycle of length four containing three colors. \label{fig:oneedge2}}
\end{figure}

\item[] Case I(ii).  $C$ contains edge $(b_1, a_1)$.  This case is symmetric to Case I(i).

\item[] Case I(iii). $C$ contains edge $(b_2, a_1)$.  This results in an $E$-cycle of length four on at most three colors as shown in Figure \ref{fig:oneedge3}(a).  Such a cycle is triangulated in $G'(S)$, implying there is either an $E$-edge or an $F$-edge between $b_2$ and $w$.  Then the cycle $(b_2, w)$, $(w, b_0)$, $(b_0, z)$, $(z, b_2)$ is either a $E$-chordless cycle in $G'(S)$ (a contradiction since all $E$-cycles on at most three colors are triangulated in $G'(S)$) or a chordless cycle in $G'(S)$ with exactly one $F$-edge (contradicting Lemma \ref{lemma:onefedge}).

\begin{figure}[h!]
        \centering
        \includegraphics[scale=0.35]{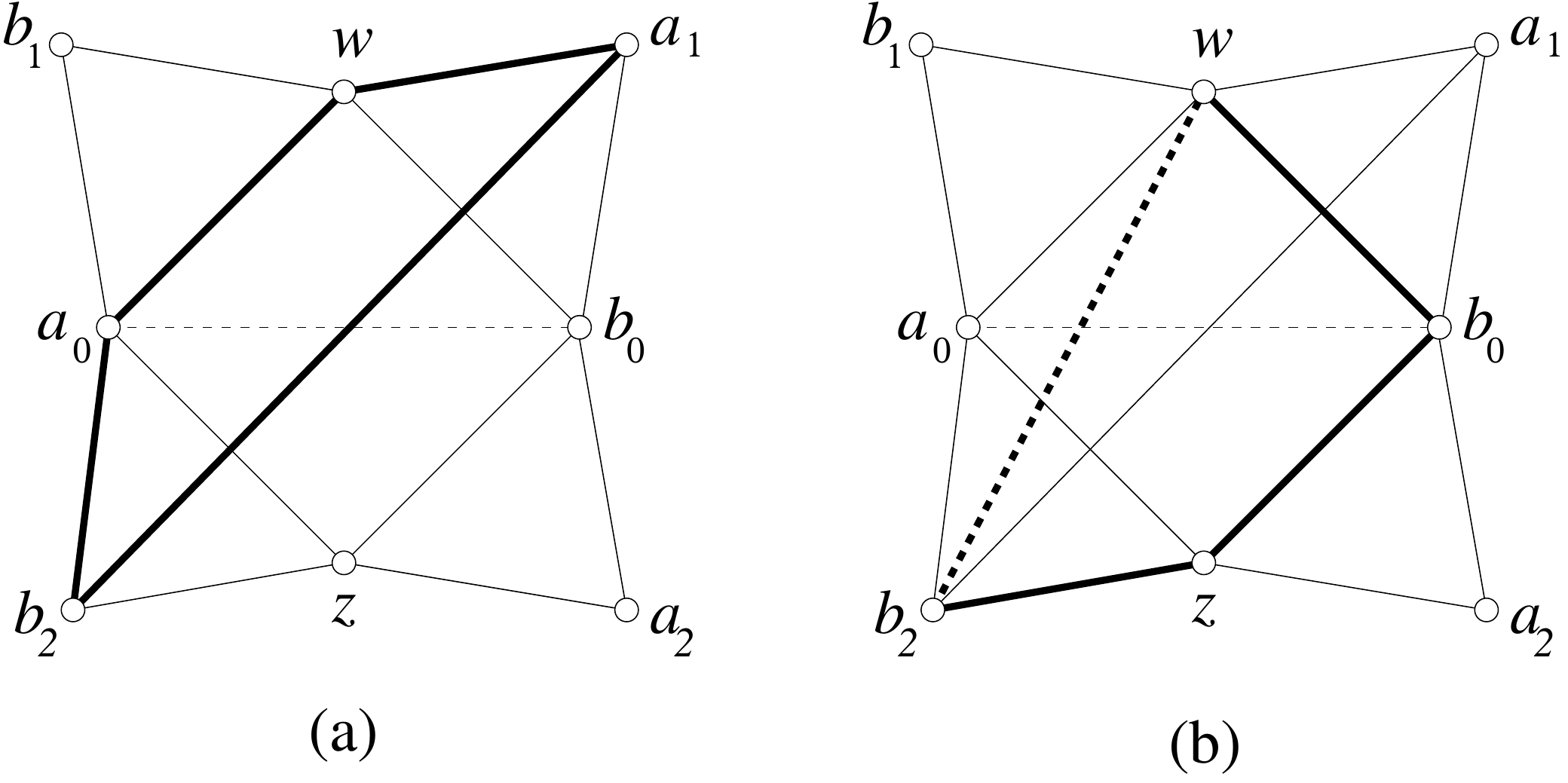}
        \caption{(a) $E$-cycle of length four on three colors (b) Chordless cycle on three colors with exactly one $F'$-edge. \label{fig:oneedge3}}
\end{figure}

\item[] Case I(iv). $C$ contains edge $(b_1, a_2)$.  This case is symmetric to Case I(iii).
\end{itemize}

 It follows that none of the vertex pairs $(b_1, a_1)$, $(b_1, a_2)$, $(b_2, a_1)$, and $(b_2, a_2)$  are adjacent by an $E$-edge and any cycle $C$ in $G''(S)$ with exactly one non-$E$ edge must contain three or more colors. Because of these nonadjacencies, the statements (A1), (A2), (B1), (B2) from Lemma \ref{lemma:onefedge} hold.  
 
 \medskip
 
\indent (A1) any row that contains $a_1$ must contain state $b_0$ in character $b$.\\   
\indent (A2) any row that contains $a_2$ must contain state $b_0$ in character $b$.\\
\indent (B1) any row that contains $b_1$ must contain state $a_0$ in character $a$.\\
\indent (B2) any row that contains $b_2$ must contain state $a_0$ in character $a$.

\medskip

\begin{figure}[h!]
        \centering
        \includegraphics[scale=0.35]{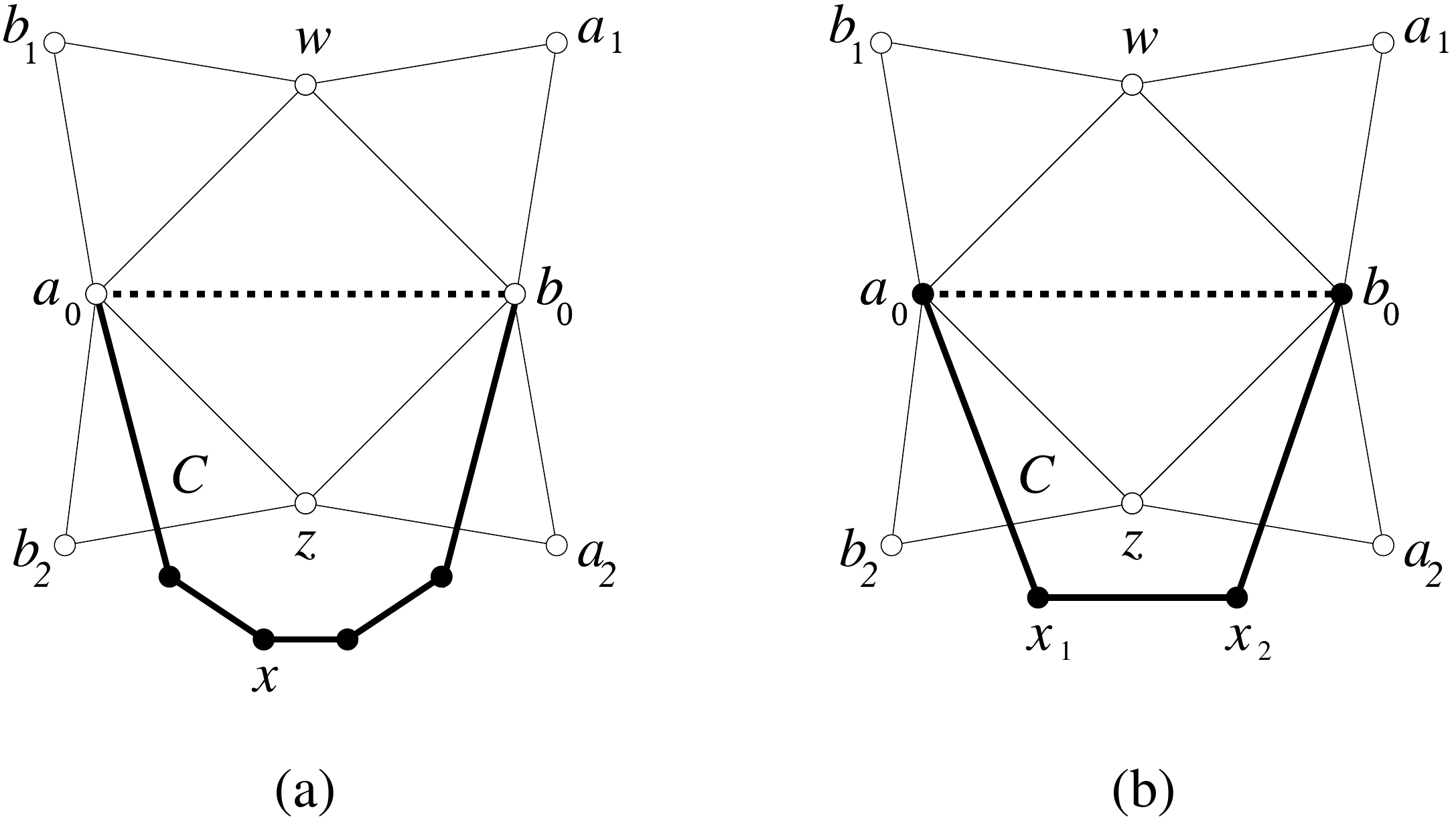}
        \caption{Vertices on cycle $C$. \label{fig:fchord7}}
\end{figure}

As in the proof of Lemma \ref{lemma:onefedge}, it follows that cycle $C$ has length equal to four formed by edges $(a_0, x_1)$, $(x_1, x_2)$, $(x_2, b_0)$, and $(b_0, a_0)$ (see Figure \ref{fig:fchord7}(b)).

Now, edge $(x_1, x_2)$ is an $E$-edge since $f=(a_0, b_0)$ is the unique non-$E$ edge in $C$ by assumption.  The remainder of the proof follows exactly as in the proof of Lemma \ref{lemma:onefedge}.  \qed

We now consider chordless cycles in $G''(S)$ with two or more non-$E$ edges.

\begin{lemma}\label{lemma:nonE} $G''(S)$ cannot contain a chordless cycle with two or more non-$E$ edges.
\end{lemma}

\Proof Suppose otherwise and let $C$ be a chordless cycle in $G''(S)$ with two or more non-$E$ edges.  Let $f$ be one of the $F$ or $F'$-edges in $C$ and without loss of generality, let $f=(a_0, b_0)$.  This edge must have been added due to an $E$-cycle $D$ that is chordless in $G(S)$ on $a_0, b_0$ and two other vertices $w$ and $z$ (see Figure \ref{fig:fchord1}(a)).  If $f$ is an $F$-edge, then $w$ and $z$ have the same color and therefore are not adjacent in $G'(S)$.  If $f$ is an $F'$-edge, then $w$ and $z$ have different colors and are nonadjacent in $G'(S)$ (since they are nonadjacent vertices in chordless cycle $D$). 

%The cycle $C$ created by edge $f$ is shown in Figure \ref{fig:fchord1}(b).

\comment{0

\begin{figure}[h!]
        \centering
        \includegraphics[scale=0.35]{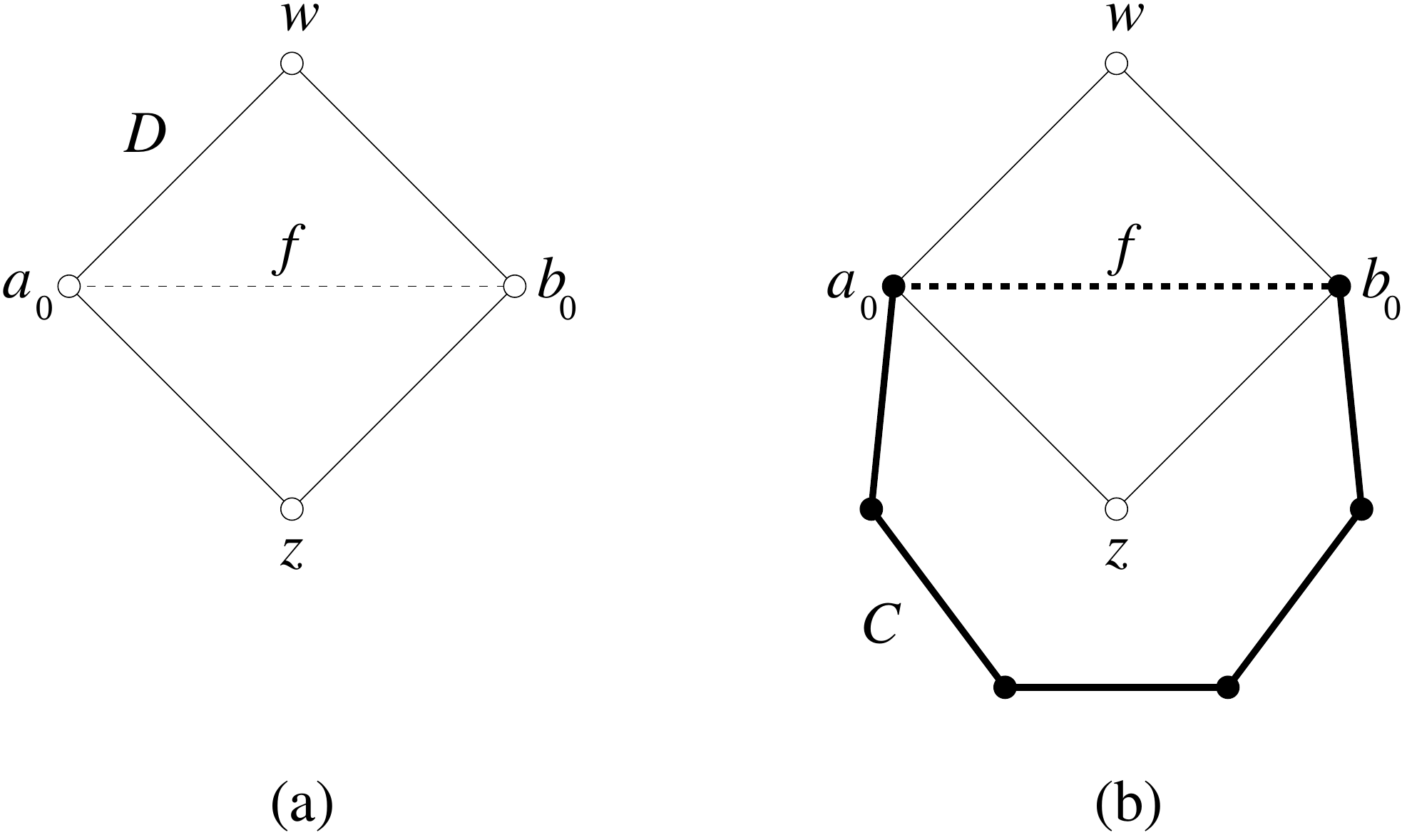}
        \caption{(a) Chordless cycle $D$ (b) cycle $C$ created by adding edge $f = (a_0, b_0)$. \label{fig:fchord1}}
\end{figure}

} % end comment0

As argued previously, the situation up to relabeling of the states is illustrated in Figure \ref{fig:fchord3}.  Furthermore, the proofs of Lemmas \ref{lemma:onefedge} and \ref{lemma:onef'edge} establish conditions (A1), (A2), (B1), and (B2), implying $C$ has length equal to four formed by edges $(a_0, x_1)$, $(x_1, x_2)$, $(x_2, b_0)$, and $(b_0, a_0)$ (see Figure \ref{fig:fchord7}(b)).  Then since $C$ has two or more non-$E$ edges, the edge $(x_1, x_2)$ in $C$ is a non-$E$ edge.  

%\begin{figure}[h!]
   %     \centering
     %   \includegraphics[scale=0.35]{fig/fchord6.pdf}
       % \caption{Impossible $E$-edges in $G''(S)$.\label{fig:fchord6}}
%\end{figure}

We have the following cases for the vertices of $C$.

\begin{figure}[h!]
        \centering
        \includegraphics[scale=0.34]{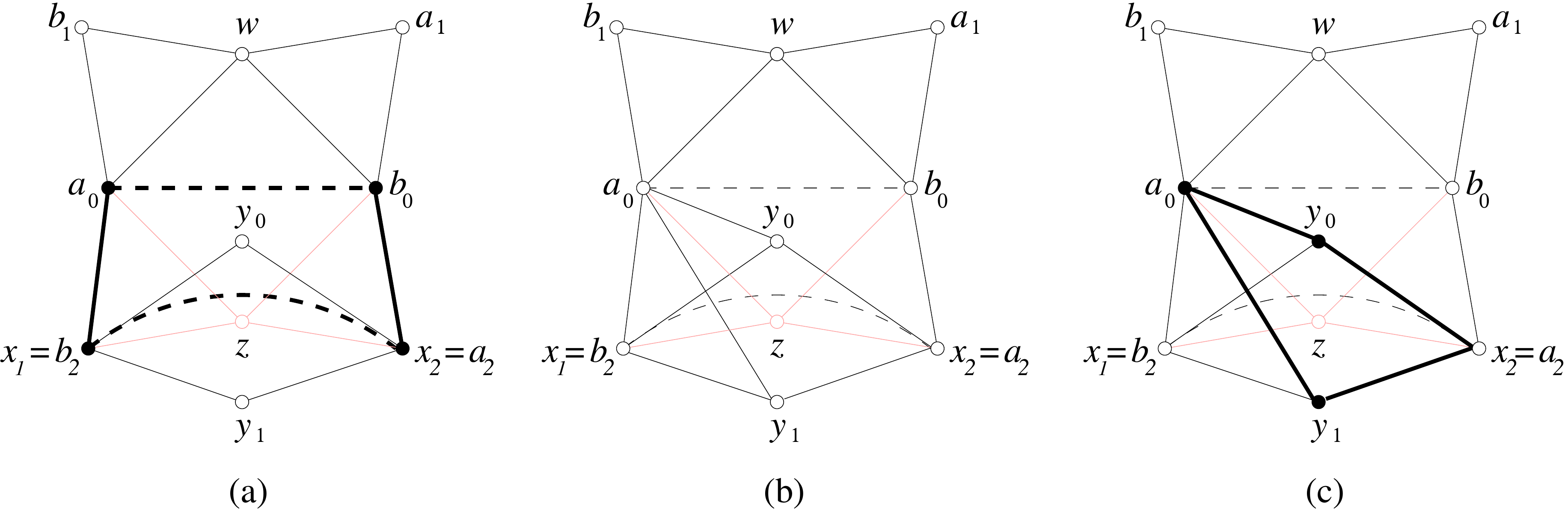}
        \caption{Case I. \label{fig:fcycle1}}
\end{figure}

\noindent {\bf Case I.}  One of $x_1$ and $x_2$ has color $a$ and the other has color $b$.  We can assume without loss of generality that $x_1 = b_2$ and $x_2=a_2$ as illustrated in Figure \ref{fig:fcycle1}(a).  Since edge $(x_1, x_2)$ is either an $F$-edge or an $F'$-edge, it was added because of a chordless $E$-cycle $D'$ containing $a_2, b_2$ and two other vertices $y_0$ and $y_1$ (see Figure \ref{fig:fcycle1}(a)).   By (A2), both $y_0$ and $y_1$ are adjacent to $a_0$, giving an $E$-cycle of length four $(a_0, y_0)$, $(y_0, a_2)$, $(a_2, y_1)$, $(y_1, a_0)$ on at most three colors (see Figure \ref{fig:fcycle1}(c)).  This $E$-cycle must be triangulated in $G'(S)$.  However, this cannot be the case since $D'$ is a chordless cycle in $G'(S)$ and $y_0$ and $y_1$ are nonadjacent vertices in $D'$.

\begin{figure}[h!]
        \centering
        \includegraphics[scale=0.35]{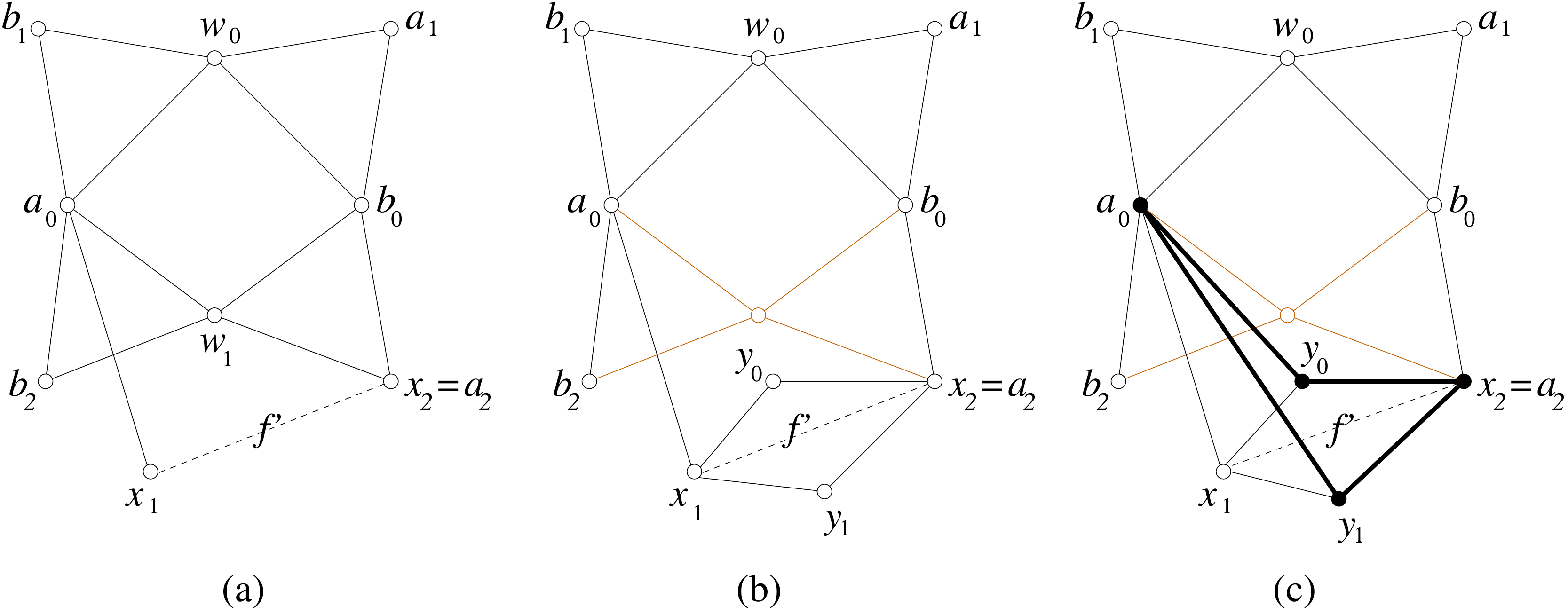}
        \caption{Case II. \label{fig:fcycle2}}
\end{figure}

\noindent {\bf Case II.} The color of $x_1$ is different from $a$ and $b$ and the color of $x_2$ is $a$ or $b$.    Without loss of generality, assume $x_2=a_2$ as illustrated in Figure \ref{fig:fcycle2}(a).  Since edge $(x_1, x_2) ( = (x_1, a_2))$ is an $F$-edge or $F'$-edge, it was added due to a chordless four cycle on $x_1, x_2 ( =a_2)$ and two other vertices $y_0$ and $y_1$.  The row witnesses for edges $(x_1, y_0)$ and $(x_1, y_1)$ must contain state $a_0$ (otherwise, $x_1$ would be adjacent to $b_0$ by (A1) or (A2)).  Then we have the $E$-cycle of length four $(a_0, y_0)$, $(y_0, a_2)$, $(a_2, y_1)$, $(y_1, a_0)$ on at most three colors.  This $E$-cycle must be triangulated in $G'(S)$.  However, this cannot happen since $y_0$ and $y_1$ are nonadjacent vertices in cycle $D'$.

\begin{figure}[h!]
        \centering
        \includegraphics[scale=0.35]{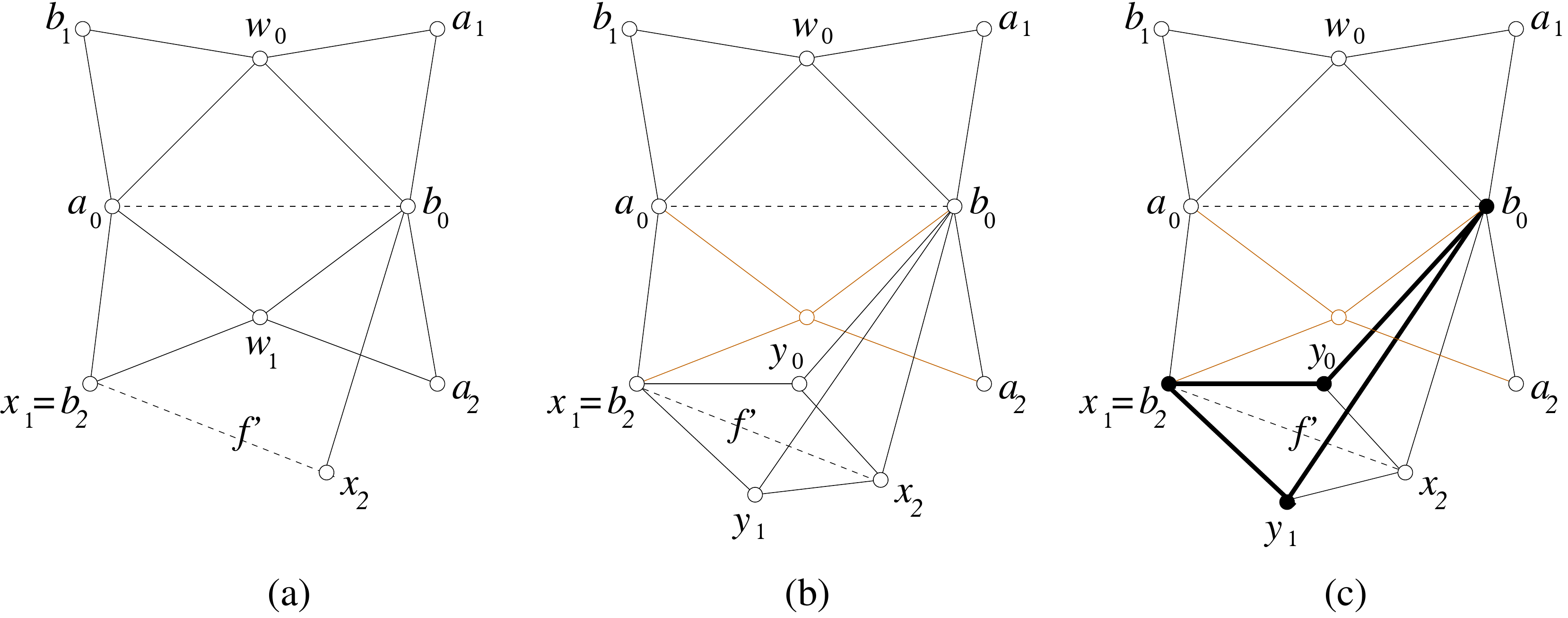}
        \caption{Case II'. \label{fig:fcycle2i}}
\end{figure}

\noindent {\bf Case II'.} The color of $x_1$ is $a$ or $b$ and the color of $x_2$ is different from $a$ and $b$.  This case is symmetric to that in Case II and is shown in Figure \ref{fig:fcycle2i}.

\begin{figure}[h!]
        \centering
        \includegraphics[scale=0.35]{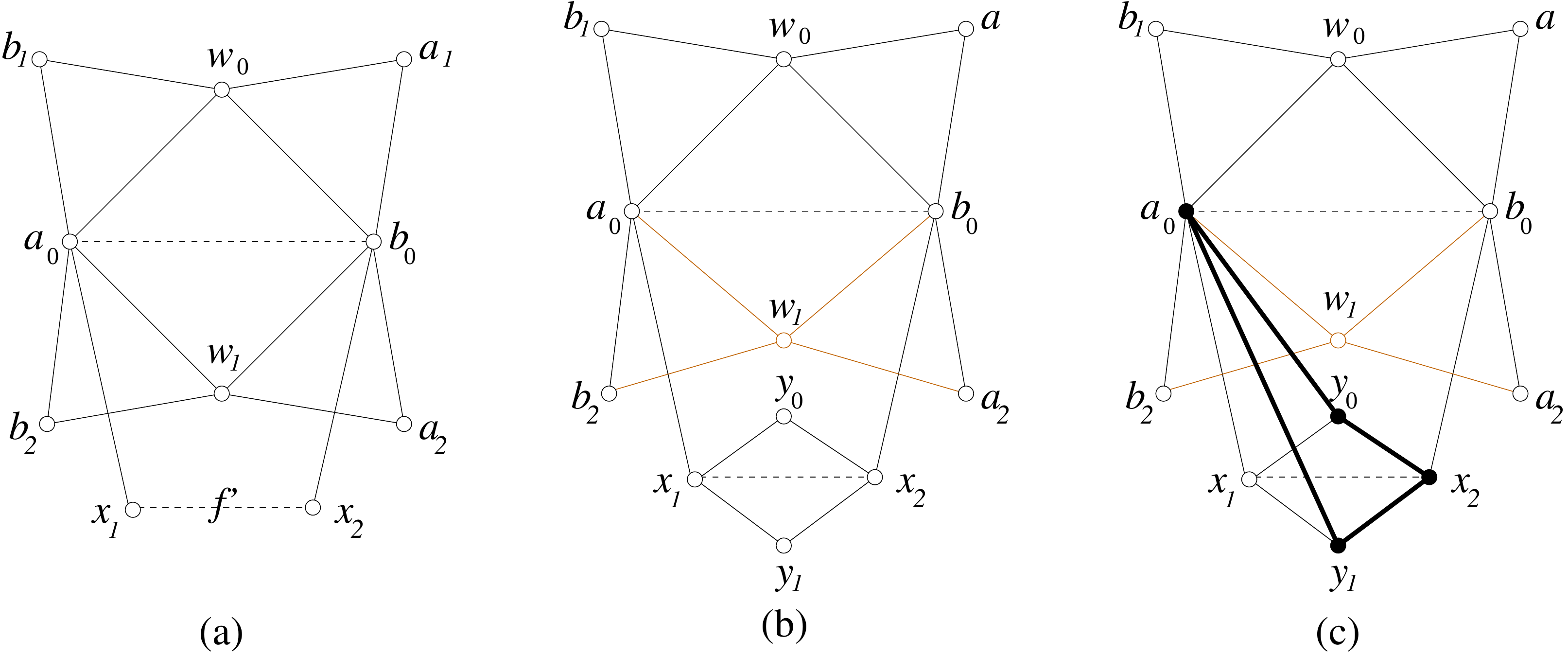}
        \caption{Case III. \label{fig:fcycle3}}
\end{figure}

\noindent {\bf Case III.} Both $x_1$ and $x_2$ have colors different from $a$ and $b$.  Since edge $(x_1, x_2)$ is an $F$-edge or $F'$-edge, it was added due to a chordless four cycle on $x_1, x_2$ and two other vertices $y_0$ and $y_1$.  The row witnesses for edges $(x_1, y_0)$ and $(x_1, y_1)$ must contain state $a_0$ (otherwise, $x_1$ would be adjacent to $b_0$ by (A1) or (A2)).  Then $(a_0, y_0)$, $(y_0, x_2)$, $(x_2, y_1)$, $(y_1, a_0)$ is an $E$-cycle of length four  (see Figure \ref{fig:fcycle3}(c)).  Note that this cycle is chordless in $G'(S)$, since $a_0$ and $x_2$ are nonadjacent vertices in chordless cycle $C$ and $y_0$ and $y_1$ are nonadjacent vertices in chordless cycle $D'$.  If $y_0$ and $y_1$ have the same color, then $C$ has only three colors and this would force edge $(a_0, x_2)$ to be an $F$-edge, a contradiction to the assumption that $C$ is chordless in $G''(S)$.  Therefore, the colors of $a_0, x_2, y_0, y_1$ are all distinct.  This cycle would force edges $(a_0, x_2)$ and $(y_0, y_1)$ to be added as $F'$-edges, a contradiction since cycle $C$ is chordless in $G''(S)$ and $a_0$ and $x_2$ are nonadjacent vertices in $C$.

This proves the lemma. \qed

\begin{figure}[h!]
        \centering
        \includegraphics[scale=0.35]{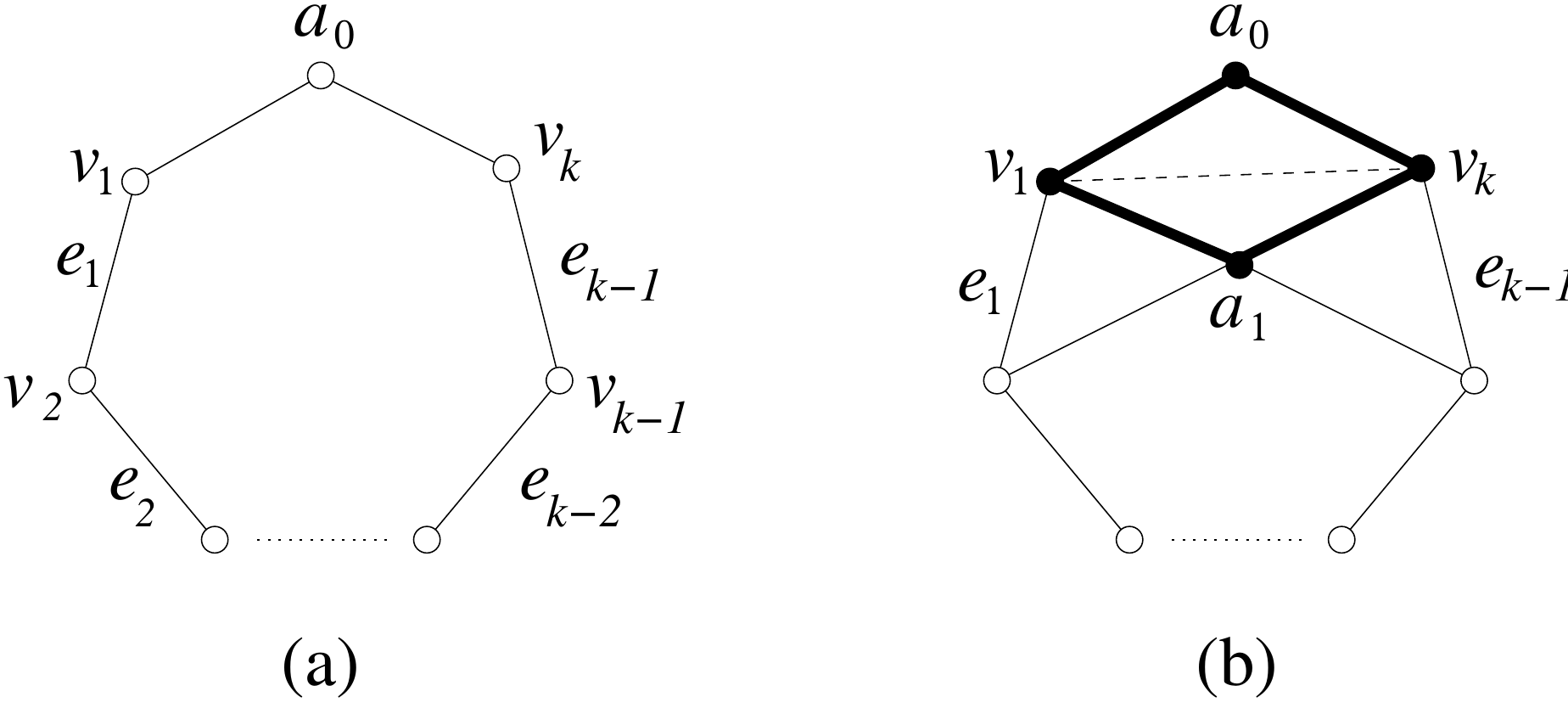}
        \hspace*{0.5cm}
          \caption{Chordless Cycle $C$}\label{fig:forcededges2}
\end{figure}

Lemmas \ref{lemma:onefedge}, \ref{lemma:onef'edge}, and \ref{lemma:nonE} eliminate the possibility of chordless cycles in $G''(S)$ containing non-$E$ edges.  To show that $G''(S)$ is properly triangulated, we proceed to show that $G''(S)$ does not contain chordless $E$-cycles.  Suppose $C$ is an $E$-cycle of length five or greater that is chordless in $G'(S)$ and suppose there is a character $a$ that appears exactly once (say in state $a_0$) in $C$.  Label the edges of the path $C \backslash {a_0}$ in order of appearance by $e_1, e_2, e_3, \ldots e_{k-1}$ with $e_i = (v_i, v_{i+1})$. Since $C$ is chordless and all edges in $C$ are $E$-edges, each edge $e_i$ ($i = 1,2, \ldots k-1$) must be witnessed by a row $s_i$ which contains either state $a_1$ or $a_2$ in color $a$.  Without loss of generality, assume $e_1$ is witnessed by $a_1$ and let $j$ be the largest index such that $e_j$ is witnessed by $a_1$.  If $j$ is equal to $k-1$, then this creates a four cycle $(v_1, a_0), (a_0, v_k), (v_k,a_1), (a_1, v_1)$ on $E$-edges (see Figure \ref{fig:forcededges2}(b)). Since $v_1$ and $v_k$ are nonadjacent (by the chordlessness of $C$ in $G'(S)$), this creates an $E$-cycle on at most three colors that is chordless in $G'(S)$, which cannot occur.

\begin{figure}[h!]
        \centering
        \includegraphics[scale=0.35]{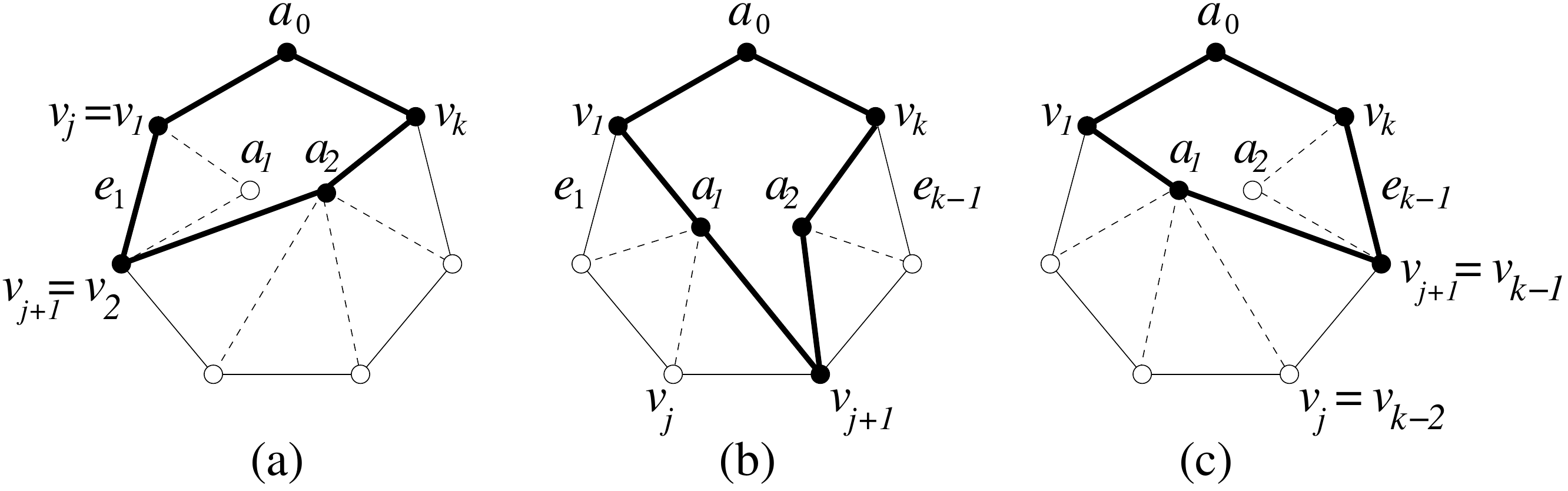}
        \caption{The $a$-complete cycle induced by color $C$ and state
 $a_0$ $I(C,a_0)$ (a) $j = 1$ (b) $1 < j < k-1$ (c) $j = k-2$. \label{fig:completecycle}}
\end{figure}

Therefore, $j$ must be strictly less than $k-1$ and all of the remaining edges $e_{j+1}, \ldots e_{k-1}$ are witnessed by state $a_2$.  Define the \emph{$a$-complete cycle induced by cycle $C$ and state $a_0$} as follows (see Figure \ref{fig:completecycle}): 

$$I(C,a_0) = 
\begin{cases} 
	(a_0, v_1), (v_1, v_2), (v_2, a_2), (a_2, v_k), (v_k, a_0) \text{ if $j=1$}  \\
	(a_0, v_1), (v_1, a_1), (a_1, v_{j+1}), (v_{j+1}, a_2), (a_2, v_k), (v_k, a_0) \text{ if 		$1 < j < k-2$} \\
	(a_0, v_1), (v_1, a_1), (a_1, v_{k-1}), (v_{k-1}, v_k), (v_k, a_0) \text{ if $j=k-2$}
\end{cases}
$$

\begin{observation} \label{obs:completecycle} For an $E$-cycle $C$ such that 
\begin{enumerate}
\item[(i)] $C$ is chordless in $G'(S)$ 
\item[(ii)] $C$ has length five or greater 
\item[(iii)] $C$ contains a character $a$ appearing exactly once in state $a_0$,
\end{enumerate}
the $a$-complete cycle $I(C,a_0)$ exists.  Note that $I(C, a_0)$ contains at least two vertices of color $a$ and has length five or greater.
\end{observation}

We use this construction to prove the following lemma.

\begin{lemma} \label{lemma:forbiddenneighbors} Suppose $C$ is an $E$-cycle of length five or greater that is chordless in $G'(S)$ and suppose there is a character $a$ appearing uniquely in $C$ in state $a_0$.  Then the two vertices adjacent to $a_0$ in $C$ have different colors and $I(C,a_0)$ is an $E$-cycle that is chordless in $G'(S)$. \end{lemma}

\begin{figure}[h!]
        \centering
        \includegraphics[scale=0.35]{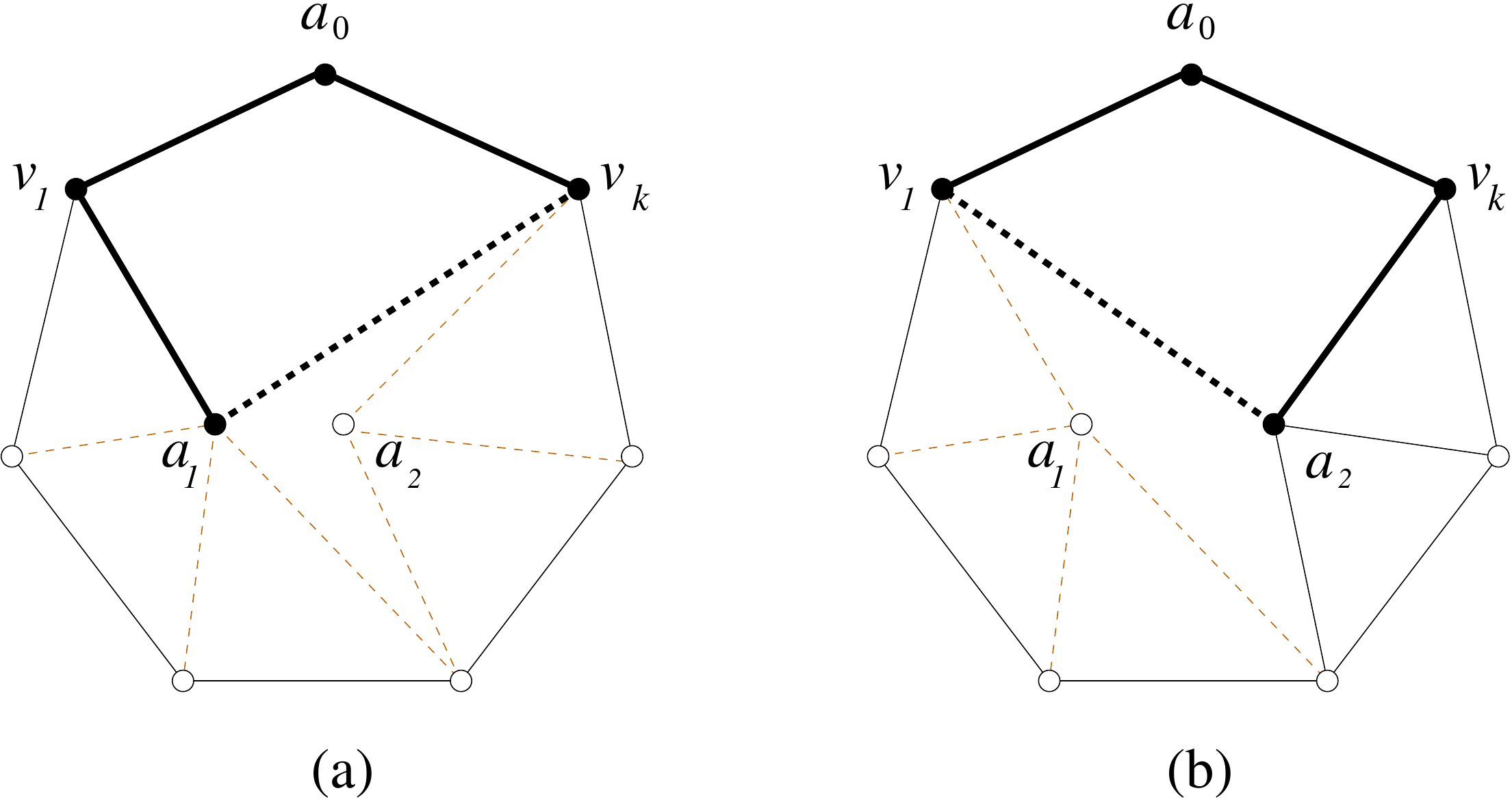}
        \caption{If either ($a_1, v_k$) or $(a_2, v_1)$ are adjacent, there is a cycle on at most three colors with at most one $F$-edge. \label{fig:completecycle2}}
\end{figure}

\Proof Note that  $I(C, a_0)$ exists by Observation \ref{obs:completecycle} and all edges in $I(C, a_0)$ are $E$-edges.  We show $I(C,a_0)$ is chordless in $G'(S)$.  The vertex pairs $(a_1, v_k)$ and $(a_2, v_1)$ are not adjacent in $G'(S)$; otherwise we would obtain a four cycle on at most three colors with at most one $F$-edge that is chordless in $G'(S)$ (see Figure \ref{fig:completecycle2}).  This is a contradiction, since Lemma \ref{lemma:onefedge} implies $G'(S)$ cannot contain a chordless cycle with at most one $F$-edge.  The remaining vertex pairs in $I(C, a_0)$ are in $C$ and are nonadjacent in $G'(S)$ since $C$ is chordless in $G'(S)$.  It follows that $I(C, a_0)$ is chordless in $G'(S)$.

Now suppose for a contradiction that the vertices adjacent to $a_0$ (vertices $v_1$ and $v_k$ in Figure \ref{fig:completecycle}) have the same color.  Then $I(C,a_0)$ is a cycle on at most three colors (color $a$, the color of $v_{j+1}$, and the shared color of vertices $v_1$ and $v_k$).  This is an $E$-cycle that has length five or greater and is chordless in the partition intersection graph on these three colors.  This  is forbidden by Lemma \ref{lemma:nofivecycle}.  Therefore, the two vertices adjacent to $a_0$ are states in two different colors.

This proves the lemma.  \qed

We now use this construction to prove properties of chordless $E$-cycles in $G'(S)$.

\begin{lemma} \label{lemma:fourdistinct} If $C$ is an $E$-cycle that is chordless in $G'(S)$, then $C$ has length exactly four with four distinct colors. \end{lemma}

\noindent {\bf Proof.} Suppose $C$ is a chordless $E$-cycle in $G'(S)$.  Note that $C$ must contain four or more colors since any chordless $E$-cycle on at most three colors is triangulated in $G'(S)$.  We first show every color in $C$ appears uniquely.  Suppose otherwise and let $a$ be the color that appears the most often in $C$ with $f_a$ the number of times $a$ appears.  We consider the following cases.

\noindent {\bf Case I.} $f_a = 3$, i.e., all three states $a_0, a_1,$ and $a_2$ appear in $C$.

If there is an edge $e=(u,v)$ in $C$ that does not have any of $a_0, a_1,$ or $a_2$ 
as endpoints, then consider the row $r$ that witnesses edge $e$; row $r$ must contain 
some state of $a$, say $a_i$. This implies edges $(u,a_i)$ and $(v,a_i)$ are present 
in $G'(S)$ and $C$ is not chordless, a contradiction. Therefore, in this case, every edge $e$ in $C$ must have exactly one endpoint of color $a$.

\begin{figure}[h!]
        \centering
        \includegraphics[scale=0.35]{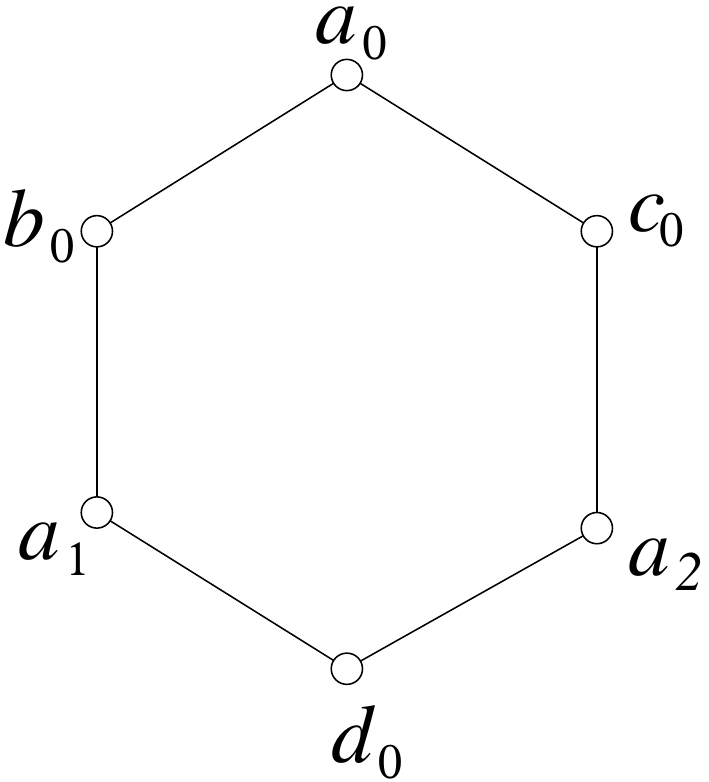}
        \caption{Color pattern in Case I. \label{fig:colorpattern2}}
\end{figure}

Since $C$ contains four or more colors and every edge is adjacent to a state of $a$, by possibly renaming the character states, the color pattern must be as shown in Figure \ref{fig:colorpattern2} (with distinct colors $b$, 
$c$, and $d$).  Now, since $C$ has length at least five and color $b$ appears uniquely in cycle $C$, the $b$-complete graph $I(C,b_0)$ induced by $C$ and $b_0$ exists.  However, the vertices adjacent to $b_0$ in $C$ are the same color (both having color $a$), which is forbidden by Lemma \ref{lemma:forbiddenneighbors}.  It follows that $f_a<3$.

\noindent {\bf Case II.} $f_a=2$: let $a_0$ and $a_1$ be the two states of $a$ appearing in $C$. Since $C$ contains four or more colors, it must be the case that one of the paths from $a_0$ to $a_1$ has three or more edges. 
We have the following cases.

\begin{enumerate}
\item[] Case (IIa) both paths from $a_0$ to $a_1$ have at least three edges \\
\item[] Case (IIb) one path from $a_0$ to $a_1$ has two edges and the other path has three or more edges
\end{enumerate}

\begin{figure}[h!]
        \centering
        \includegraphics[scale=0.35]{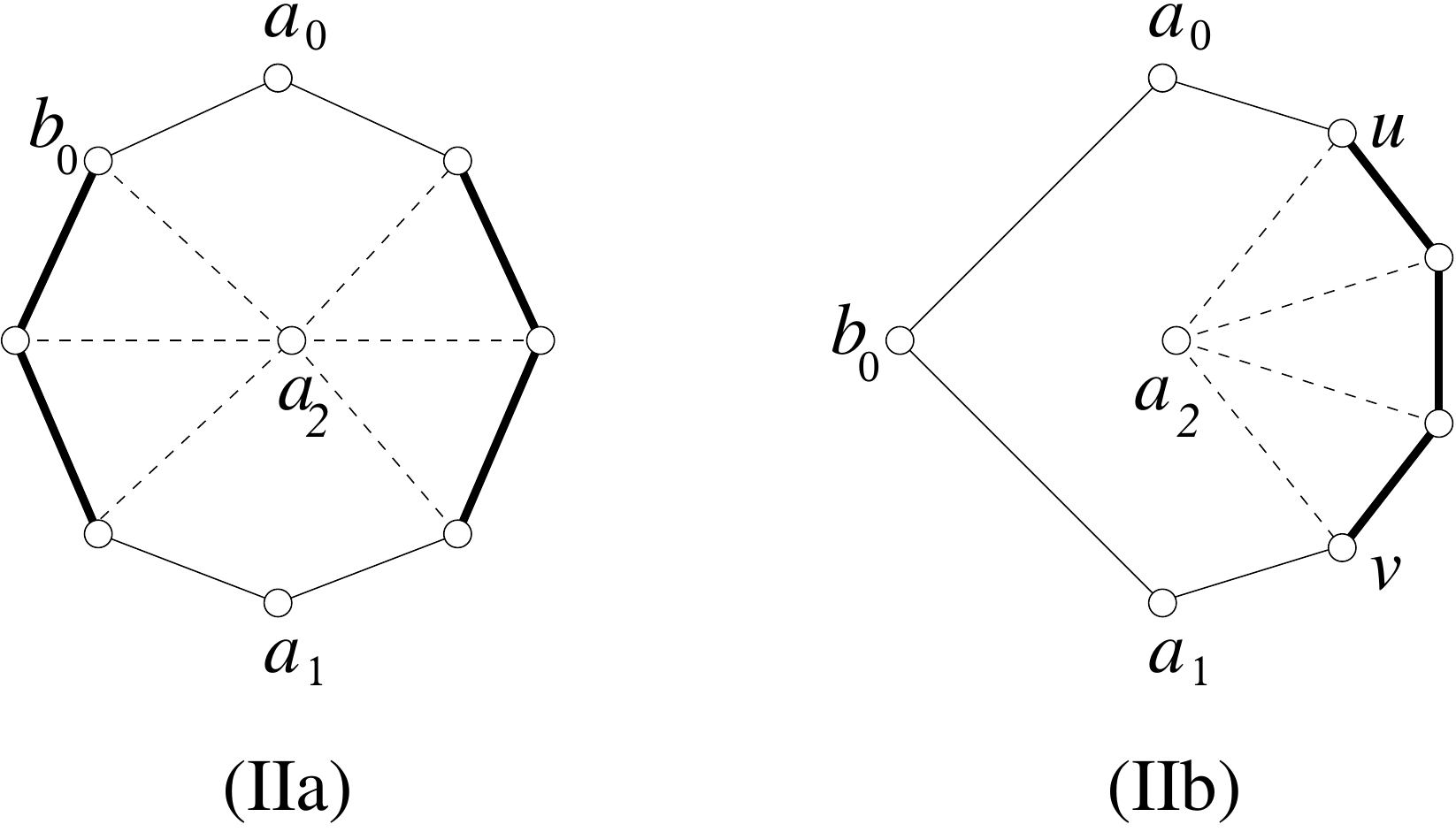}
        \hspace*{0.5cm}
        \caption{Cases (IIa) and (IIb) in the proof of Lemma \ref{lemma:fourdistinct}.  The rows witnessing the edges shown in bold must contain state $a_2$ in character $a$.  \label{fig:caseII}}
\end{figure}

Any edge that does not have color $a$ as one of its endpoints must be witnessed by a row that contains the third state $a_2$, as illustrated in Figure \ref{fig:caseII}.

\begin{figure}[h!]
        \centering
        \includegraphics[scale=0.35]{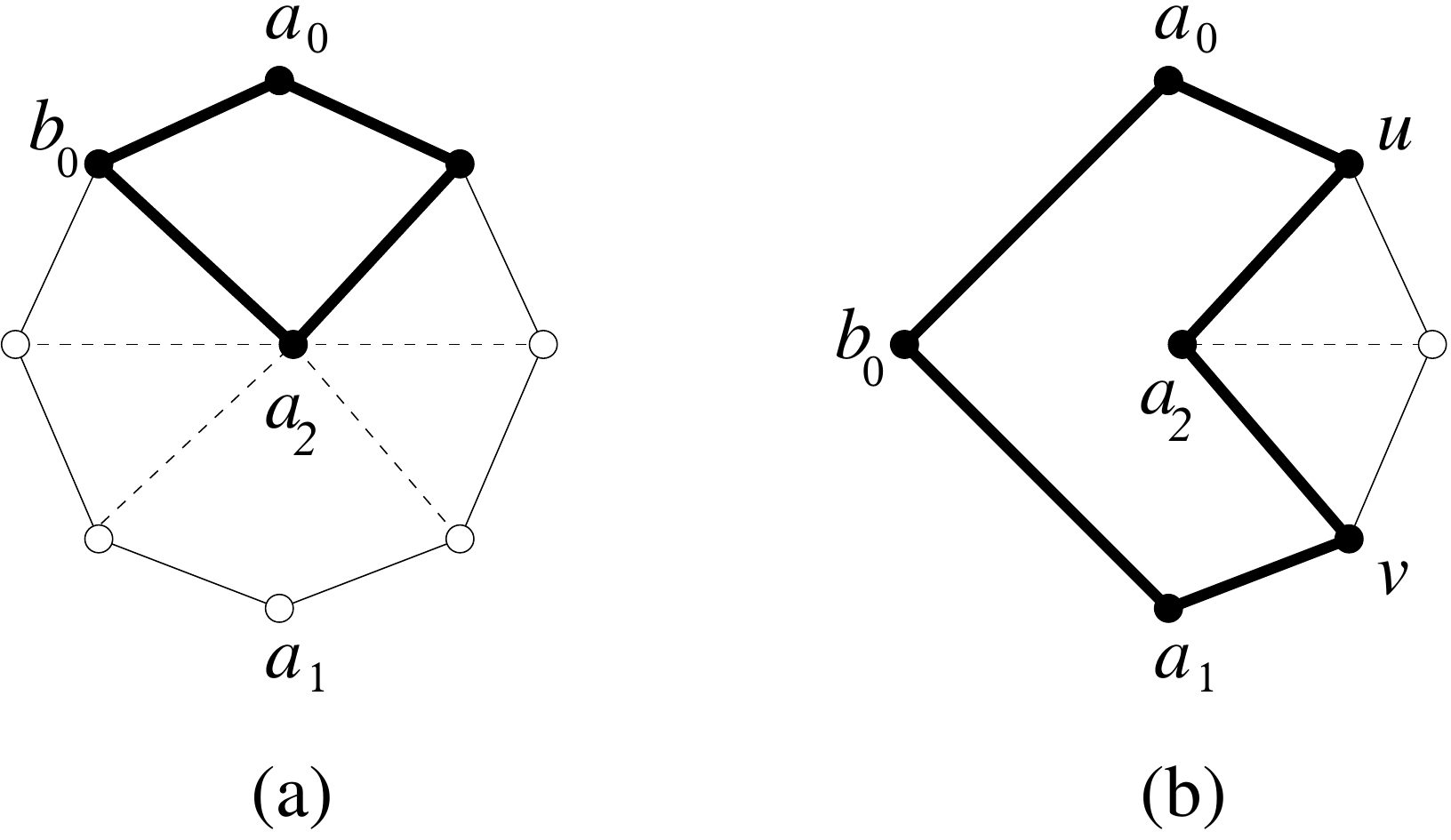}
        \caption{ Case II.\label{fig:caseIIi}}
\end{figure}

In case (IIa), the second edge in both paths from $a_0$ to $a_1$ are witnessed by state $a_2$ and we obtain an $E$-cycle that is chordless in $G'(S)$ on at most three colors (shown in bold in Figure \ref{fig:caseIIi}(a)).  This is a contradiction since all $E$-cycles on at most three colors must be triangulated in $G'(S)$.  In case (IIb), the second and second to last edge on the $a_0$ to $a_1$ path with three or more edges are witnessed by color $a_2$.  Let $D$ denote the $E$-cycle of edges $(b_0, a_0)$, $(a_0, u)$, $(u, a_2)$, $(a_2, v)$, $(v, a_1)$, $(a_1, b_0)$ (shown in Figure \ref{fig:caseIIi}(b)).  Then $a_2$ and $b_0$ are not adjacent in $G'(S)$ (otherwise, we would obtain a cycle of length four on at most three colors with at most one $F$-edge).  This implies $D$ is an $E$-cycle that is chordless in $G(S)$; in this cycle, $b_0$ has two adjacent vertices of color $a$ and therefore cannot be the only state of $b$ appearing in $D$, by Lemma \ref{lemma:forbiddenneighbors}.  This implies one of $u$ or $v$ must also have color $b$ and therefore $D$ is a chordless cycle on at most three colors containing all three states of character $a$, contradicting Lemma \ref{lemma:nothreestates}.

\noindent {\bf Case III.} $f_a=1$, i.e., every color in $C$ appears uniquely.   Suppose for a contradiction that $C$ has length five or greater and let $a_0$ be a state appearing in $C$.  Then $I(C,a_0)$ exists in $G'(S)$ by Observation \ref{obs:completecycle}.  However, this gives a chordless cycle in $G'(S)$ with color $a$ appearing two or more times, which cannot happen by Cases I and II.

It follows that $C$ is a cycle of length four with all colors appearing uniquely in $C$, proving the lemma. \qed
  
Lemma \ref{lemma:fourdistinct} implies all chordless $E$-cycles in $G'(S)$ have length four containing four distinct colors.  We have triangulated all such cycles by $F'$-edges in $G''(S)$, implying the following corollary.

\begin{corollary} \label{cor:ecycle} $G''(S)$ cannot contain a chordless $E$-cycle.  
\end{corollary}

Lemmas \ref{lemma:onefedge}, \ref{lemma:onef'edge}, \ref{lemma:nonE}, and Corollary \ref{cor:ecycle} together imply that $G''(S)$ is properly triangulated, proving the main theorem.

\bigskip

\noindent {\bf Theorem \ref{thm:main}} \emph{Given an input set $S$ on $m$ characters with at most three states per character ($r=3$), $S$ admits a perfect phylogeny if and only if every subset of three characters of $S$ admits a perfect phylogeny.}

\section{Enumerating Obstruction Sets for Three State Characters}\label{sec:forbidden}

We now turn to the problem of enumerating all minimal obstruction sets to perfect phylogenies on three-state character input.  By Theorem \ref{thm:main}, it follows that the minimal obstruction sets contain at most three characters.  We enumerate all instances $S$ on three characters $a,b,$ and $c$ satisfying the following conditions:

\begin{quote}
\begin{enumerate}
\item[(i)] each character $a, b$ and $c$ has at most three states
\item[(ii)] every pair of characters allows a perfect phylogeny
\item[(iii)] the three characters $a$, $b$, and $c$ together do not allow a perfect phylogeny.
\end{enumerate}
\end{quote}

Note that Condition (ii) implies the partition intersection graph $G(S)$ does not contain a cycle on exactly two colors and Condition (iii) implies $G(S)$ contains at least one chordless cycle.  Let $C$ be the largest chordless cycle in $G(S)$, i.e.,

$$C = {\arg\max}_{\text{chordless cycles $D$ in $G(S)$}} \vert D \vert$$

\begin{figure}[h!]
        \centering
        \includegraphics[scale=0.35]{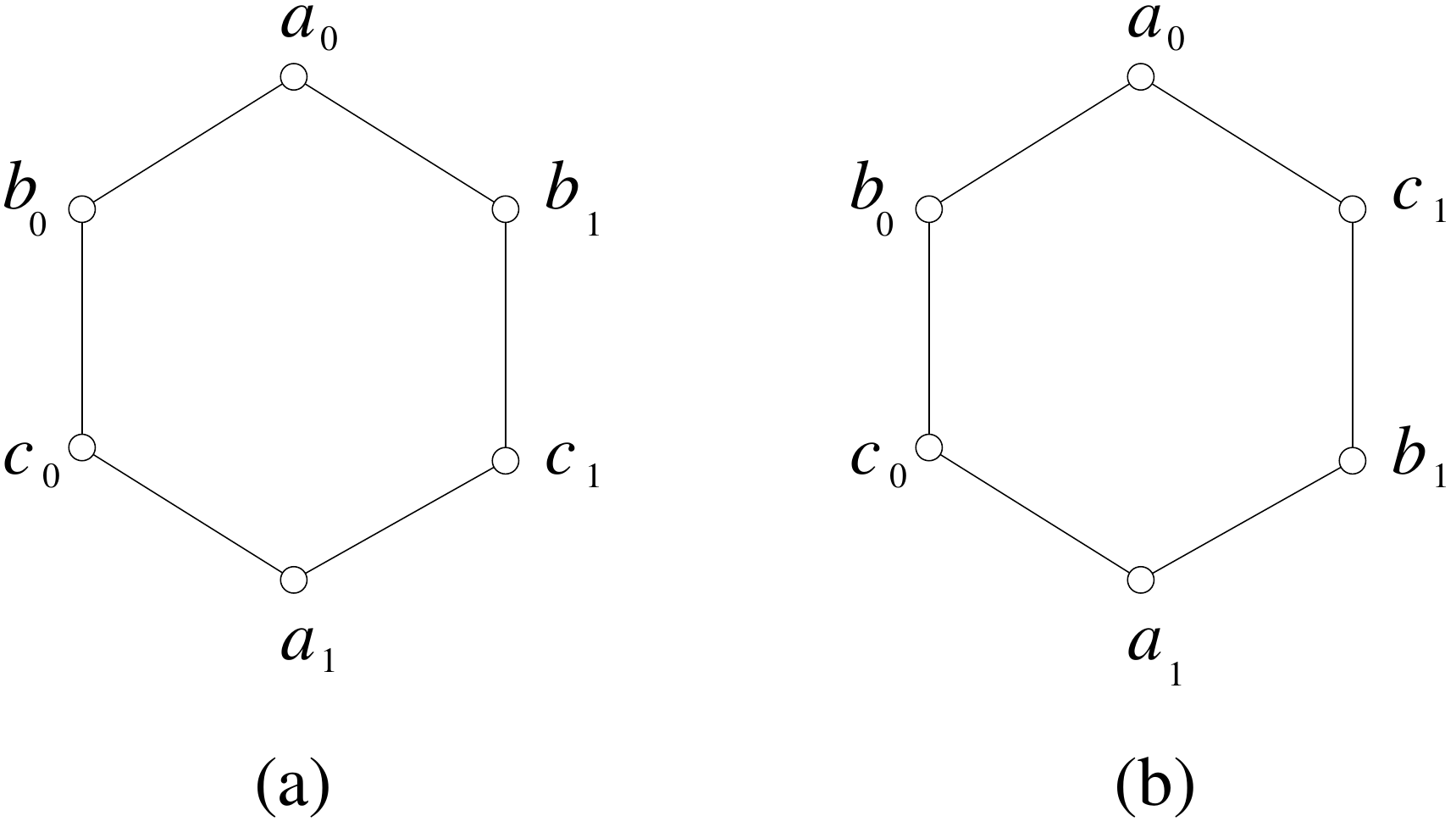}
        \caption{Color patterns for chordless cycle of length six.\label{fig:six_cycle1}}
\end{figure}

Condition (ii) and Lemma \ref{lemma:nothreestates} together imply $C$ cannot contain all three states of any character.  Therefore, $C$ has length at most six.  If $G(S)$ contains a chordless six-cycle $C$, then each color appears exactly twice in $C$ and $C$ must have one of the color patterns (up to relabeling) shown in Figure \ref{fig:six_cycle1}.

\begin{figure}[h!]
        \centering
        \includegraphics[scale=0.35]{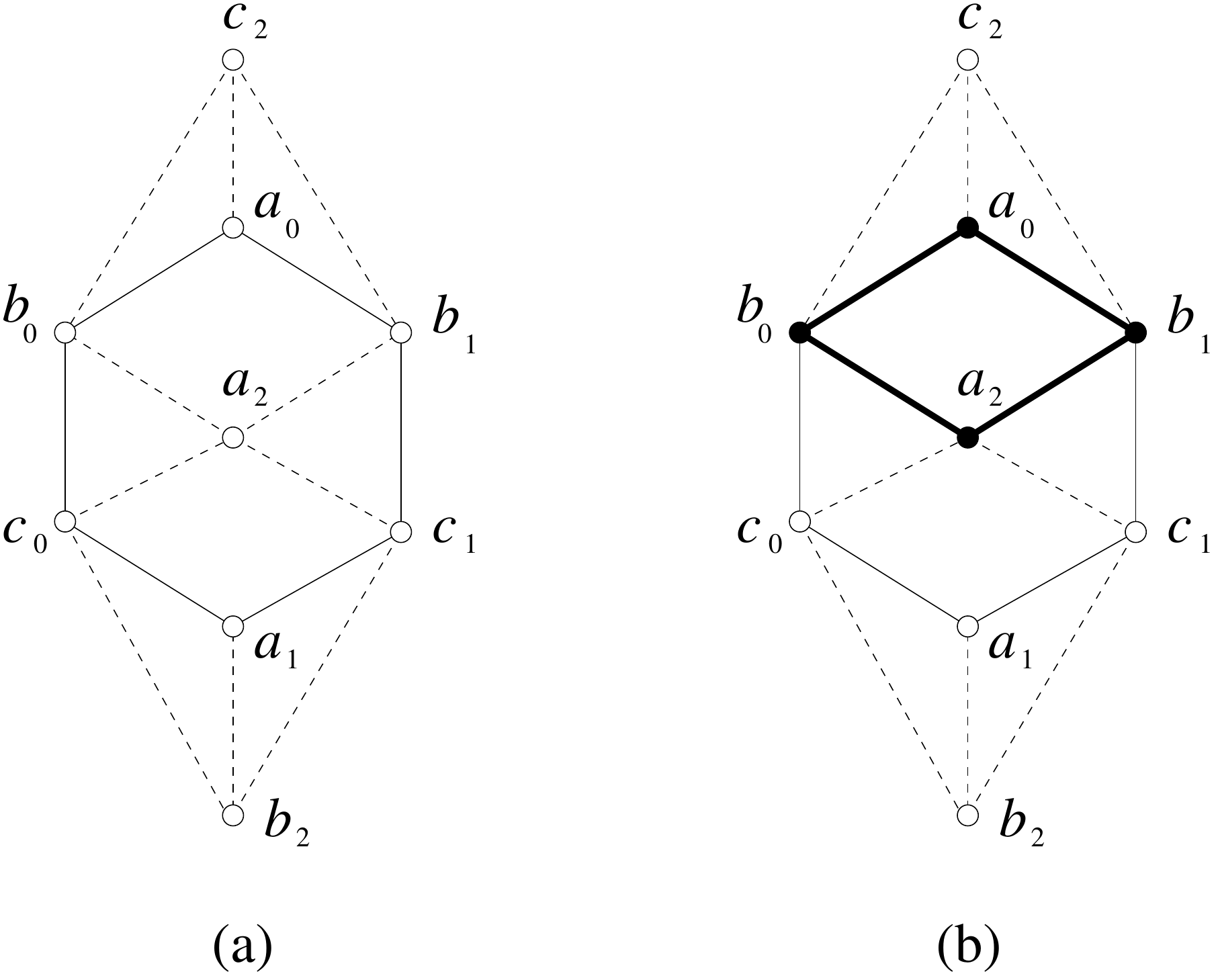}
        \caption{Forced patterns for row witnesses of Figure \ref{fig:six_cycle1}(a).\label{fig:six_cycle2}}
\end{figure}

\begin{figure}[h!]
        \centering
        \includegraphics[scale=0.35]{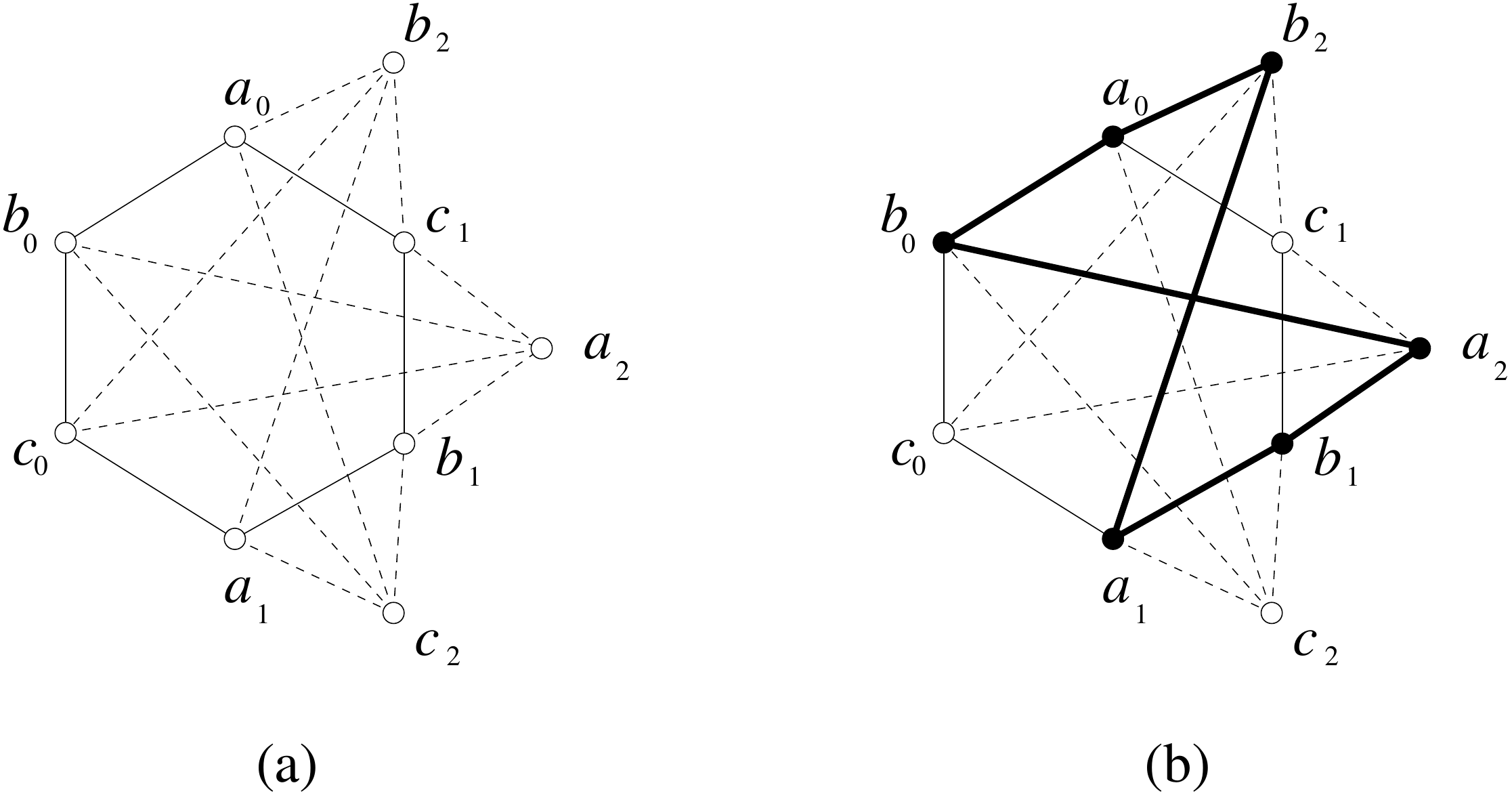}
        \caption{Forced patterns for row witnesses of Figure \ref{fig:six_cycle1}(b). \label{fig:six_cycle3}}
\end{figure}

In Figures \ref{fig:six_cycle1}(a) and \ref{fig:six_cycle1}(b), there is one state in each character that does not appear in $C$ (states $a_2, b_2$, and $c_2$).   Since $C$ is chordless, the witness for each edge is forced to contain the missing state in the third character.  This implies Figure \ref{fig:six_cycle1}(a) must be completed by the edges in Figure \ref{fig:six_cycle2}(a) and Figure \ref{fig:six_cycle1}(b) must be completed by the edges in Figure \ref{fig:six_cycle3}(a).  In both cases, there is a cycle on two characters $a$ and $b$ (see Figures \ref{fig:six_cycle2}(b) and \ref{fig:six_cycle3}(b)).  This implies the pair of characters $a$ and $b$ is not properly triangulatable, a contradiction to condition (ii).  Therefore, $G(S)$ cannot contain chordless cycles of length six.

\clearpage

\begin{figure}[h!]
        \centering
        \includegraphics[scale=0.35]{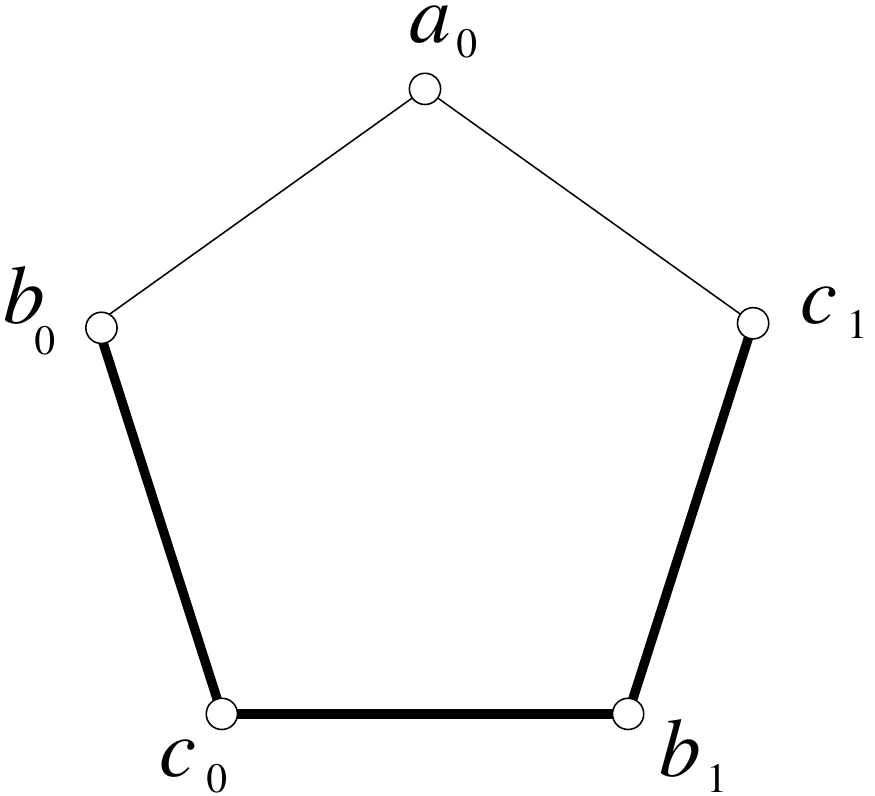}
        \caption{Color pattern for cycle $C$ of length five.\label{fig:five_cycle1}}
\end{figure}

If $C$ is a chordless cycle in $G(S)$ of length five, then $G(S)$ is not properly triangulatable by Lemma \ref{lemma:nofivecycle}, implying Condition (iii) is satisfied.  In this case, there must be two characters (say $b$ and $c$) appearing in two different states and one character appearing once in $C$, as shown in Figure \ref{fig:five_cycle1} (up to relabeling of the states).  Cycle $C$ contains three edges that are not adjacent to character $a$ (edges $(b_0, c_0)$, $(c_0, b_1)$, $(b_1, c_1)$ in Figure \ref{fig:five_cycle1}).  The row witnesses for these edges must contain either state $a_1$ or $a_2$ in character $a$.

\begin{figure}[h!]
        \centering
        \includegraphics[scale=0.35]{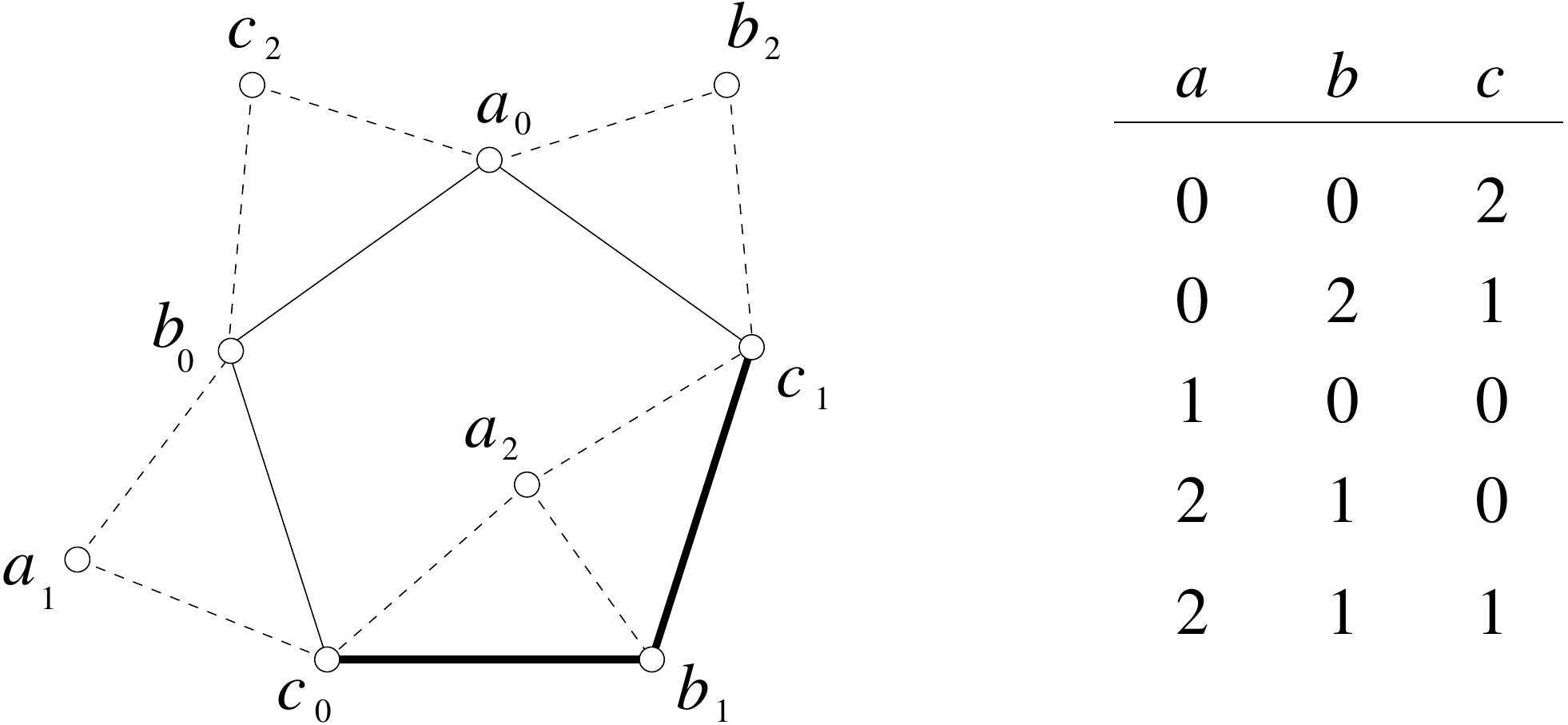}
        \caption{Case I.  Row witnesses for two adjacent edges share the same state of $a$.\label{fig:five_cycle2}}
\end{figure}

\noindent {\bf Case I.}  The row witnesses for two adjacent edges share the same state of $a$ and the row witness for the third edge contains the final state in $a$.  Without loss of generality, assume $(c_0, b_1)$ and $(b_1, c_1)$ are the two adjacent edges sharing the same state of $a$.  In this case, $G(S)$ and the corresponding input sequences $S$ are shown in Figure \ref{fig:five_cycle2} (up to relabeling of the states).
  
\begin{figure}[h!]
        \centering
        \includegraphics[scale=0.35]{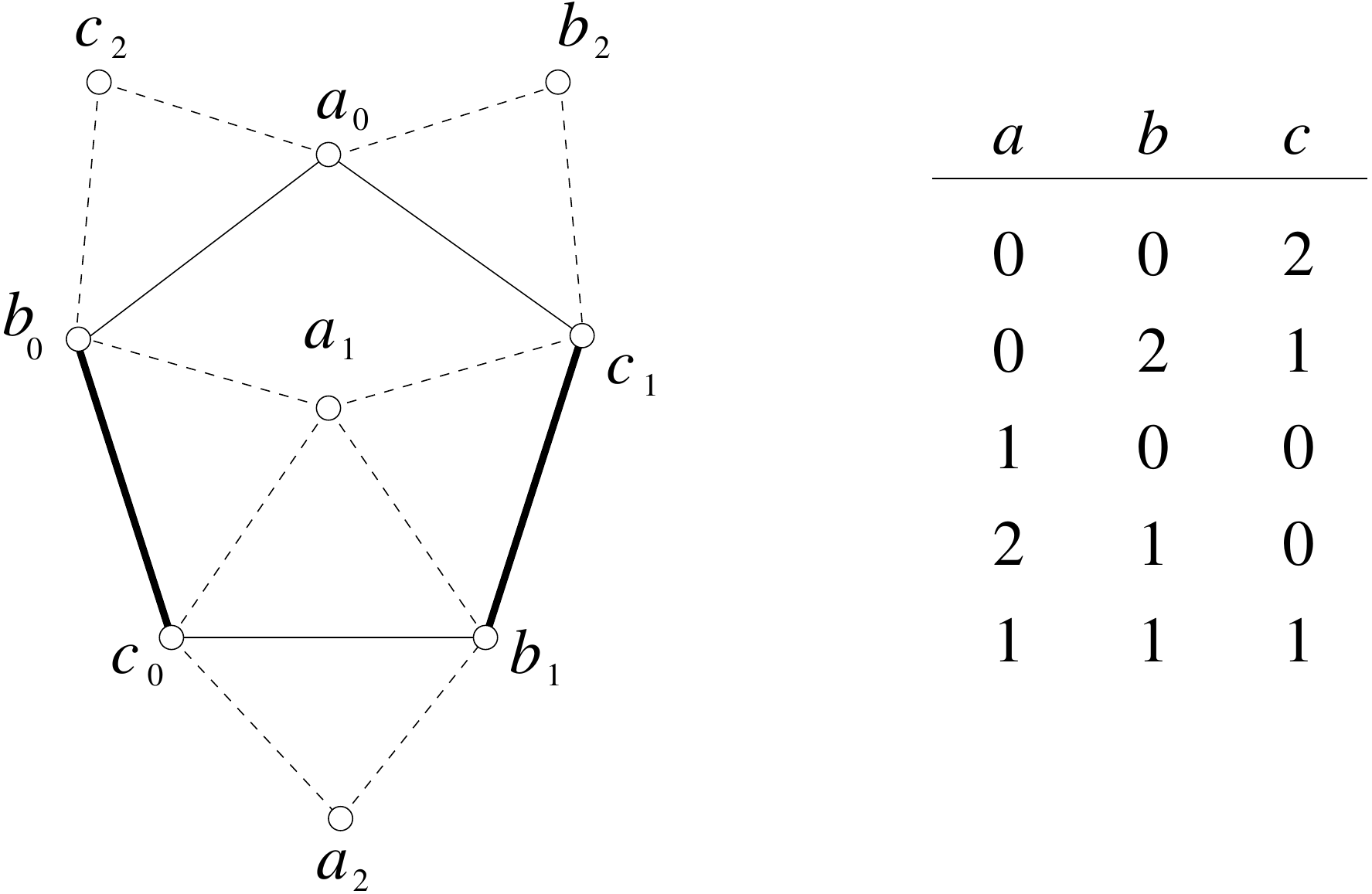}
        \caption{Case II.  Row witnesses for two nonadjacent edges share the same state of $a$.\label{fig:five_cycle3}}
\end{figure}

\noindent {\bf Case II.}  The row witnesses for the two nonadjacent edges share the same state of $a$ and the row witness for the third edge contains the final state in $a$.   In this case, $G(S)$ and the corresponding input sequences $S$ are shown in Figure \ref{fig:five_cycle3} (up to relabeling of the states).

\begin{figure}[h!]
        \centering
        \includegraphics[scale=0.35]{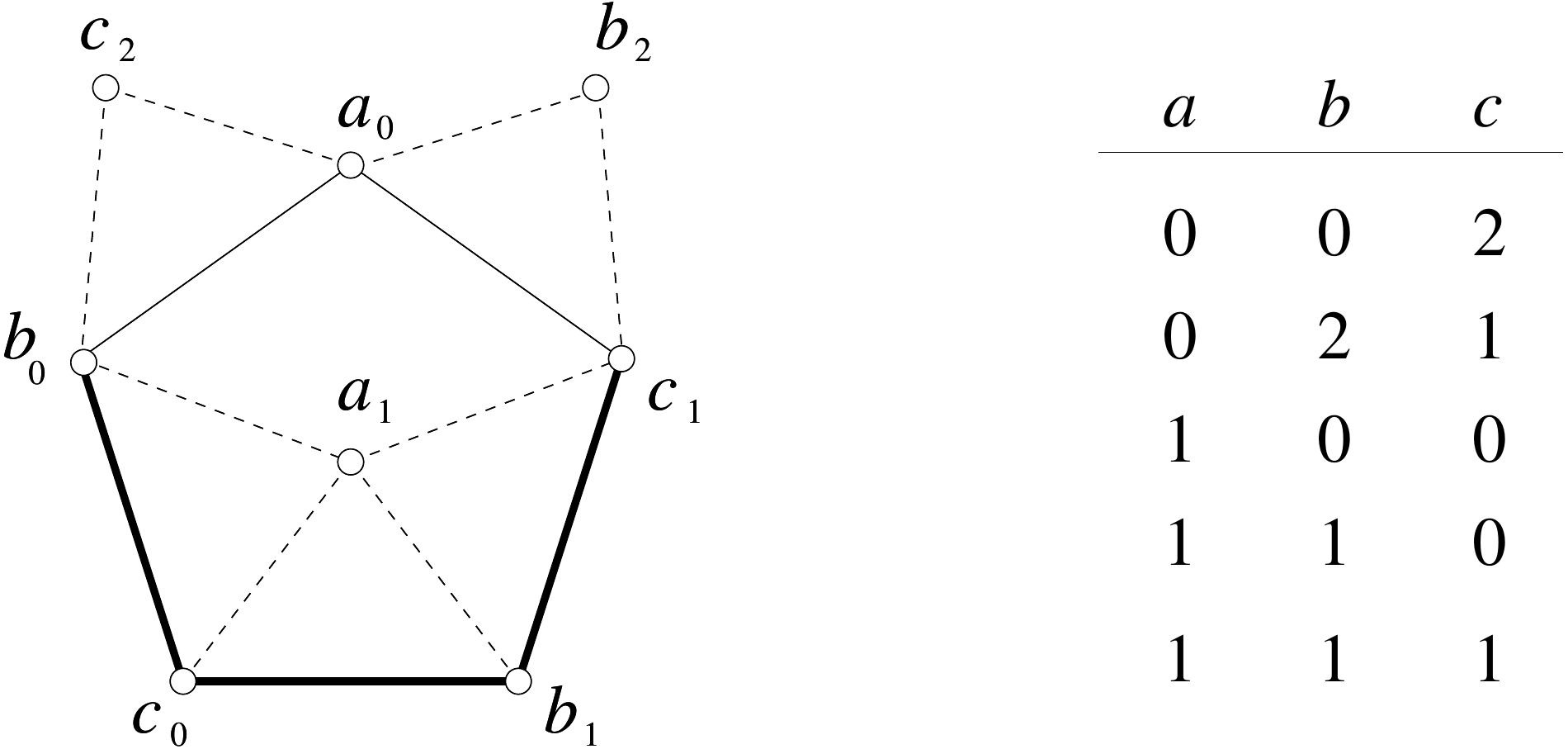}
        \caption{Case III.  Row witnesses for all three edges share the same state of $a$.\label{fig:five_cycle4}}
\end{figure}

\noindent {\bf Case III.}  The row witnesses for all three edges share the same state of $a$.   In this case, $G(S)$ and the corresponding input sequences $S$ are shown in Figure \ref{fig:five_cycle4} (up to relabeling of the states).

\clearpage

\begin{figure}[h!]
        \centering
        \includegraphics[scale=0.35]{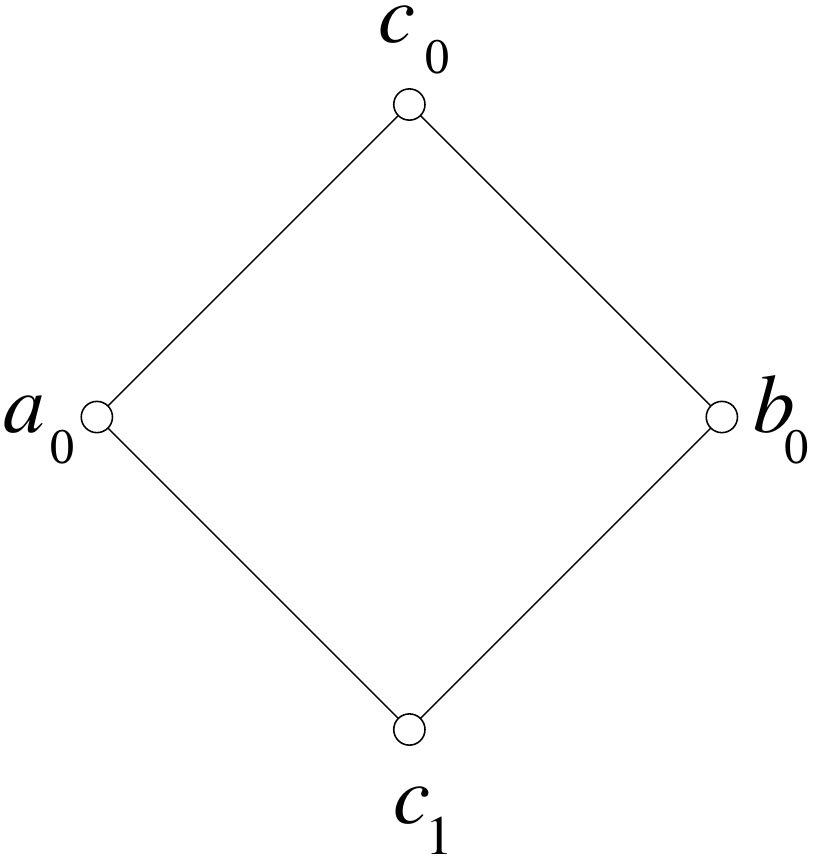}
        \caption{Color Pattern for chordless cycle of length four.\label{fig:four_cycle1}}
\end{figure}

If $C$ is a chordless cycle of length four, then without loss of generality it must have the color pattern shown in Figure \ref{fig:four_cycle1}.

\begin{figure}[h!]
        \centering
        \includegraphics[scale=0.35]{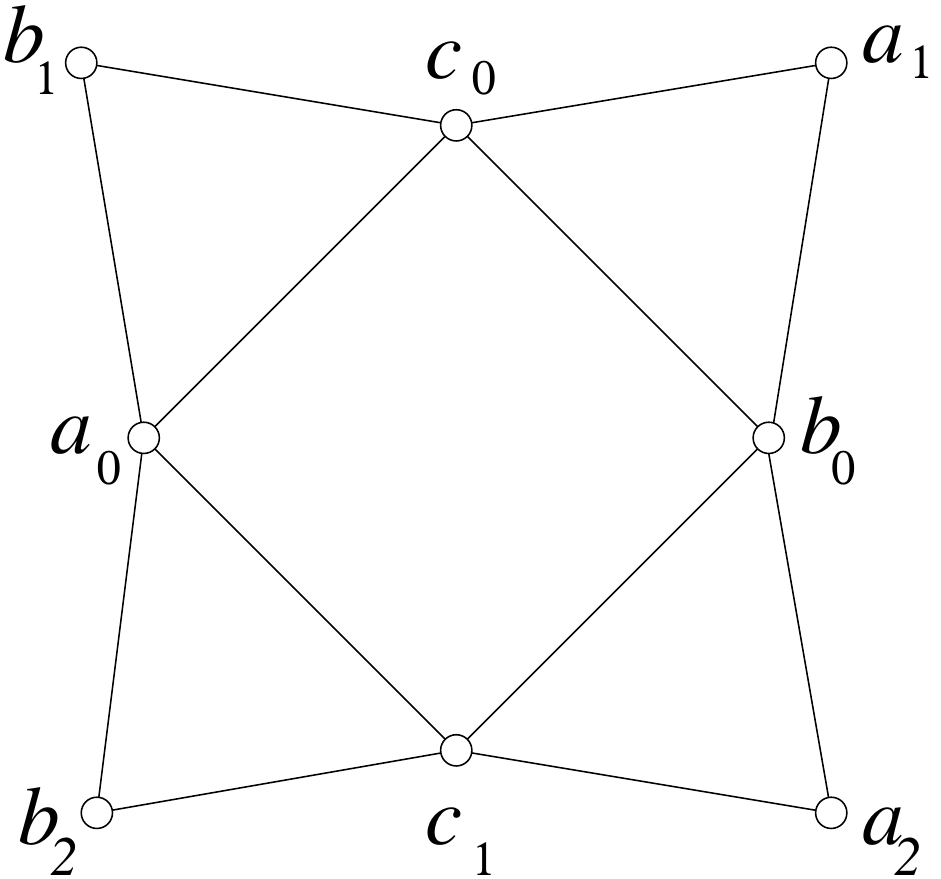}
        \caption{Forced pattern of row witnesses.\label{fig:four_cycle3}}
\end{figure}

Consider the row witnesses for edges $(a_0, c_0)$ and $(a_0, c_1)$.  These row witnesses cannot share the same state of $b$ (otherwise, there would be a cycle on two colors $b$ and $c$, a contradiction).  Similarly, row witnesses for edges $(b_0, c_0)$ and $(b_0, c_1)$ cannot share the same state of $a$.  Therefore, up to relabeling of the states, the row witnesses are forced to have the pattern shown in Figure \ref{fig:four_cycle3}.

\begin{figure}[h!]
        \centering
        \includegraphics[scale=0.35]{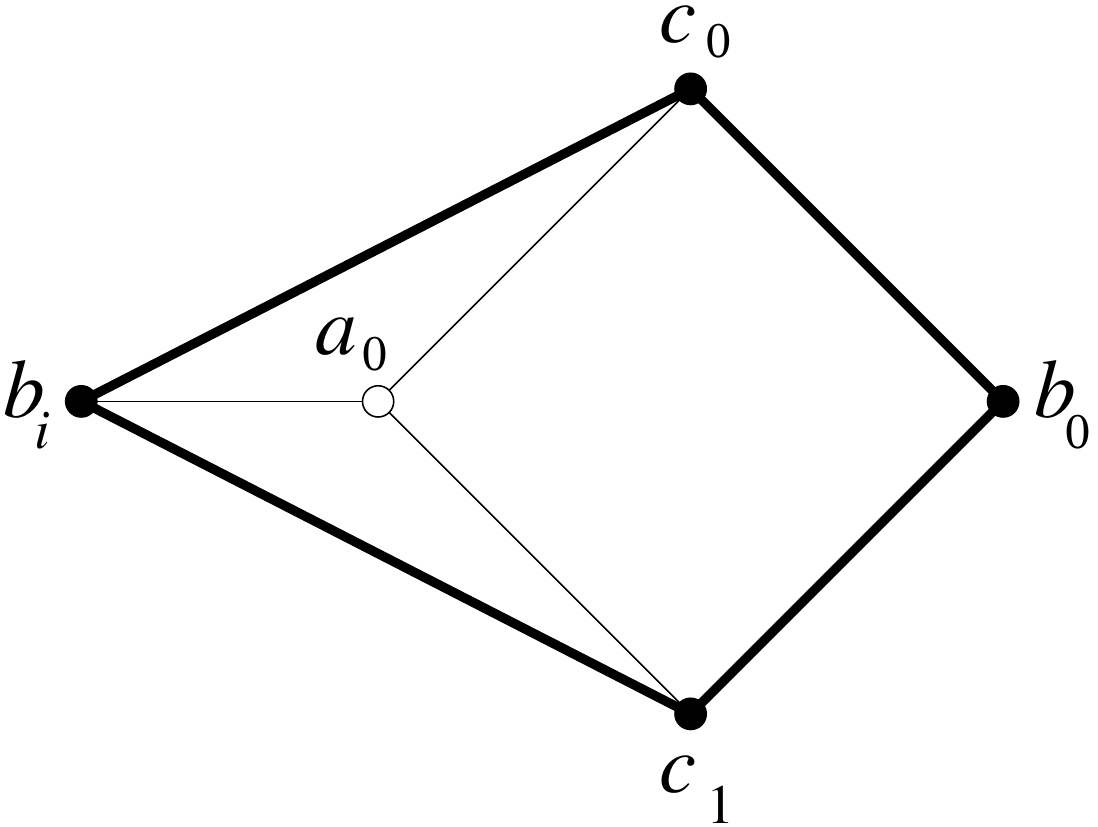}
        \caption{Row witnesses for edges $(a_0,c_0)$ and $(a_0, c_1)$ cannot share the same state of $b$.\label{fig:four_cycle2}}
\end{figure}

\begin{figure}[h!]
        \centering
        \includegraphics[scale=0.35]{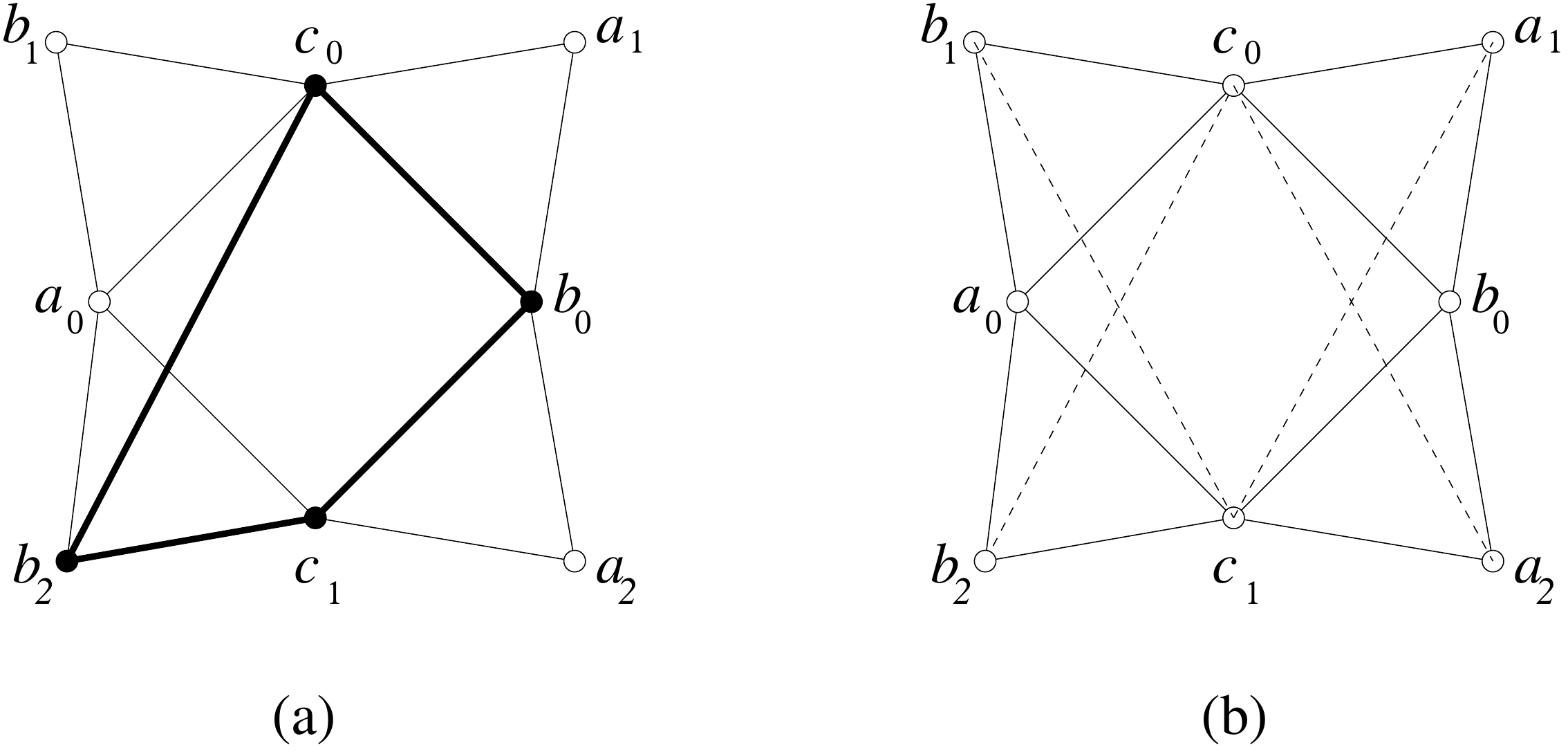}
        \caption{ $(b_2, c_0)$, $(a_2, c_0)$, $(b_1, c_1)$, $(a_1, c_1)$ induce cycles on two colors $b$ and $c$.\label{fig:four_cycle4}}
\end{figure}

Note that $b_2$ and $c_0$ cannot be adjacent in $G(S)$; otherwise, there is a cycle on two colors $b$ and $c$ (see Figure \ref{fig:four_cycle4}(a)).  By symmetry, we can argue (see Figure \ref{fig:four_cycle4})

\begin{enumerate}
\item[] pairs $(b_2, c_0)$, $(a_2, c_0)$, $(b_1, c_1)$, and $(a_1, c_1)$ are nonadjacent in $G(S)$ \hspace{5ex} (*) 
\end{enumerate}

%\begin{eqnarray}
%\text{The pairs $(b_2, c_0)$, $(a_2, c_0)$, $(b_1, c_1)$, and $(a_1, c_1)$ are %nonadjacent in $G(S)$} 
%\end{eqnarray}
%

\begin{figure}[h!]
        \centering
        \includegraphics[scale=0.35]{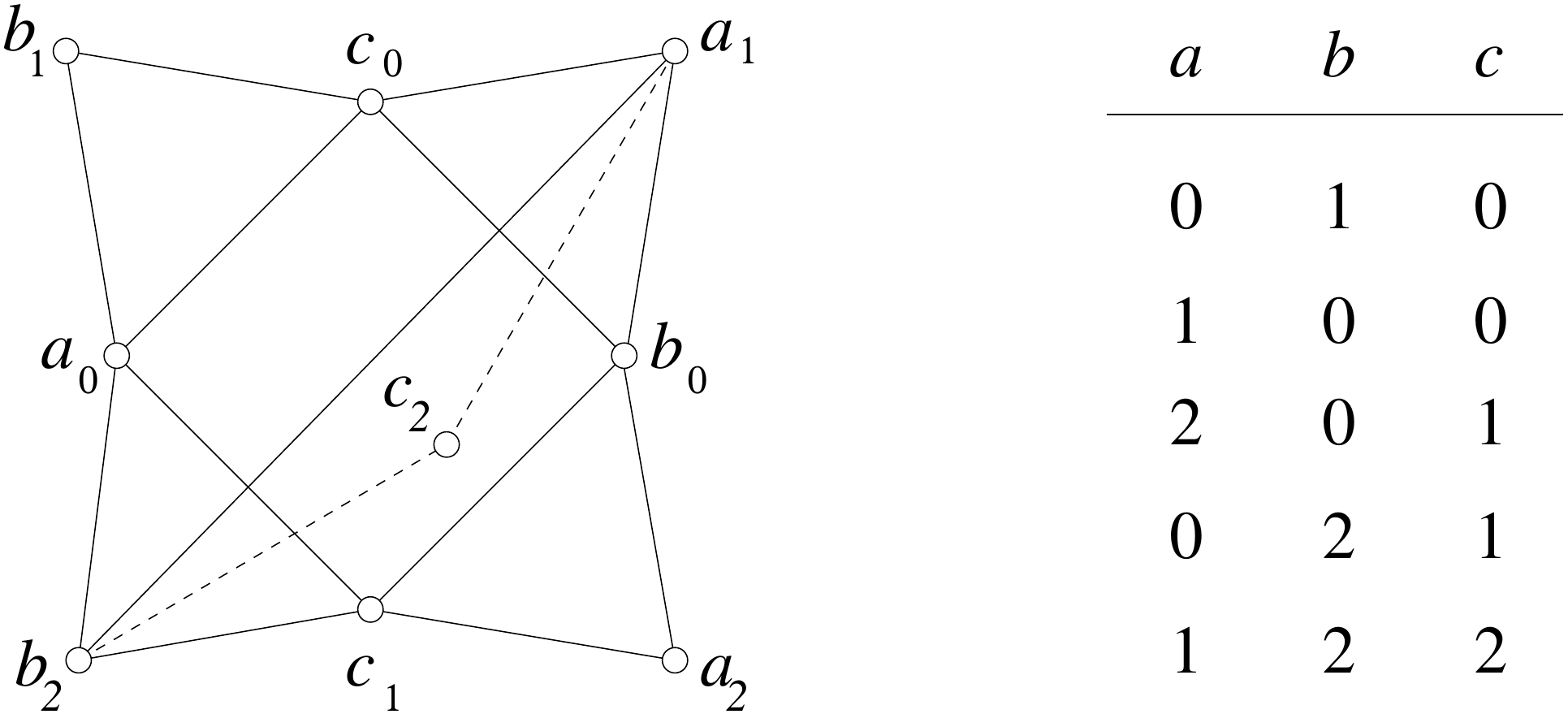}
        \caption{Input sequences $S$ and partition intersection graph $G(S)$ with a chordless cycle of length four.\label{fig:four_cycle5}}
\end{figure}

Now, suppose $b_2$ and $a_1$ are adjacent in $G(S)$.  Then the row witness for $(b_2, a_1)$ cannot be $c_0$ and cannot be $c_1$ by (*).  Therefore, the row witness for this edge must be the third state $c_2$ of character $c$ (see Figure \ref{fig:four_cycle5}).  The partition intersection graph $G(S)$ and corresponding input sequences $S$ are shown in Figure \ref{fig:four_cycle5}.  Note that $G(S)$ is not properly triangulatable and condition (iii) is satisfied, since the edge $(a_0, b_0)$ is a forced edge to triangulate cycle $C$, creating a cycle on two colors $(a_0, b_0)$, $(b_0, a_1)$, $(a_1, b_2)$, $(b_2, a_0)$ which cannot be properly triangulated. 

If $b_2$ and $a_2$ are adjacent in $G(S)$, then this induces a chordless cycle $D$ of length five $(b_2, a_0)$, $(a_0, c_0)$, $(c_0, b_0)$, $(b_0, a_2)$, $(a_2, b_2)$ (the pairs $(b_2, c_0)$ and $(a_2, c_0)$ are nonadjacent by (*) and $(a_0, b_0)$ are nonadjacent since they are nonadjacent vertices in chordless cycle $C$).  This is a contradiction since $C$ is chosen to be the largest chordless cycle in $G(S)$.

\begin{figure}[h!]
        \centering
        \includegraphics[scale=0.35]{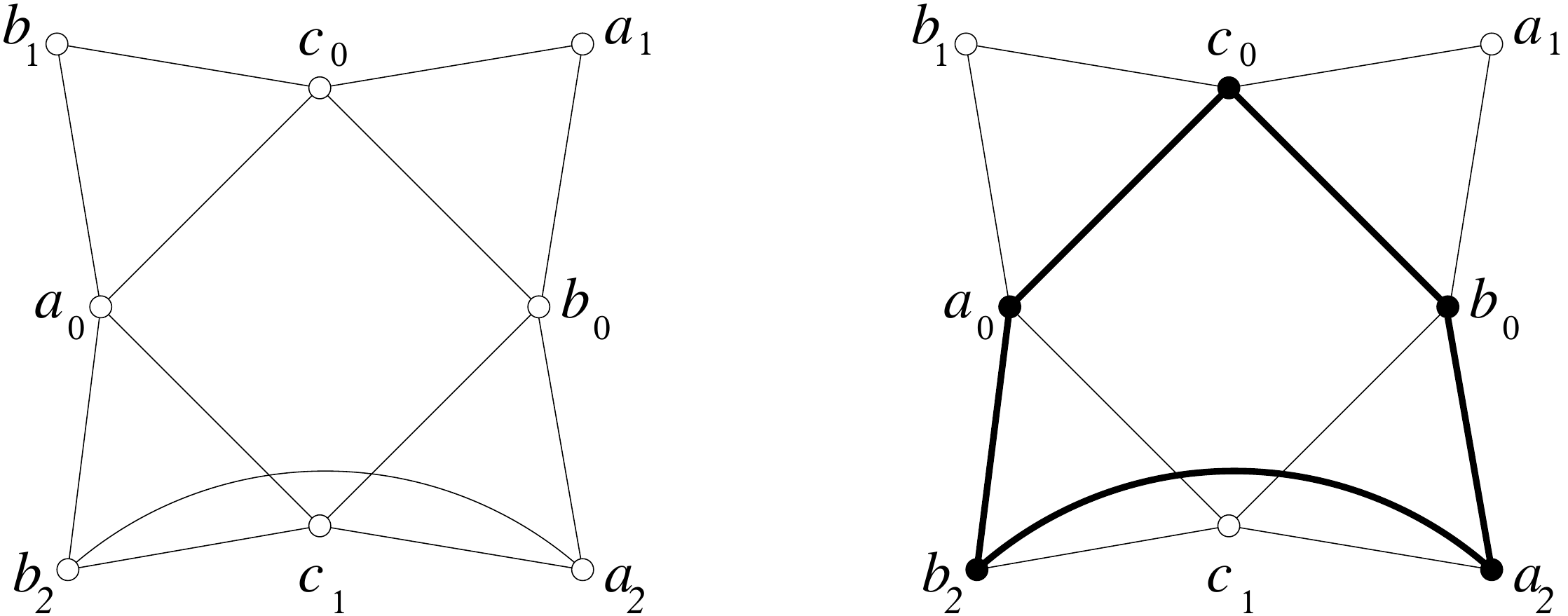}
        \caption{If $b_2$ and $a_2$ are adjacent in $G(S)$, this creates a chordless cycle of length five.\label{fig:four_cycle6}}
\end{figure}

\begin{figure}[h!]
        \centering
        \includegraphics[scale=0.35]{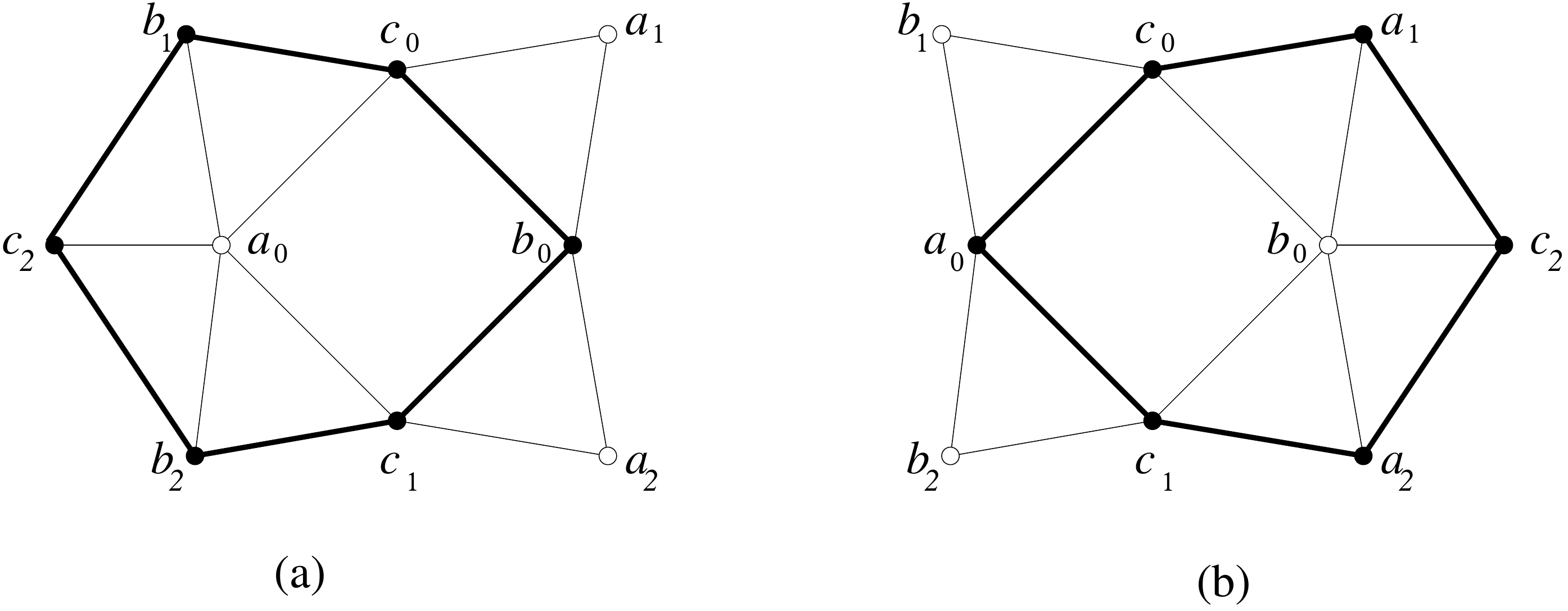}
        \caption{If $c_2$ witnesses two adjacent edges in $G(S)$, this creates a chordless cycle on two colors.\label{fig:four_cycle7}}
\end{figure}

Suppose there are no further adjacencies between vertices in Figure \ref{fig:four_cycle3}.  Then there must be additional edges formed by the final state $c_2$ of character $c$ in order for $G[a,b,c]$ to be nontriangulatable (condition (iii)).  Now, state $c_2$ is adjacent to one or more of the edges with color pattern $(a, b)$.  If $c_2$ is adjacent to exactly one such edge, then the resulting graph $G[a,b,c]$ can be properly triangulated by adding the edge $(a_0, b_0)$.  Otherwise, state $c_2$ is adjacent to two or more edges.  If the two edges share a vertex (i.e., the two edges are either $(a_1, b_0)$ and $(a_2, b_0)$ or $(b_1, a_0)$ and $(b_2, a_0)$), then there is a cycle on two colors (as shown in Figure \ref{fig:four_cycle7}(a) and \ref{fig:four_cycle7}(b)), contradicting condition (ii).  

\begin{figure}[h!]
        \centering
        \includegraphics[scale=0.35]{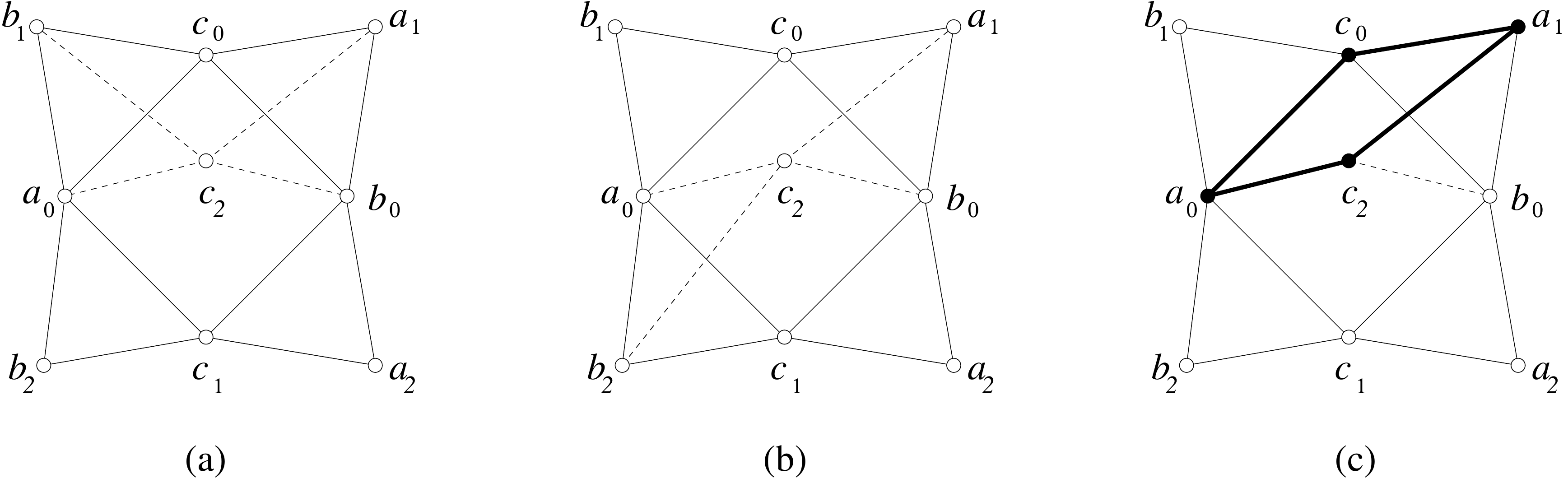}
        \caption{If $c_2$ witnesses two nonadjacent edges in $G(S)$, this creates a chordless cycle on two colors.\label{fig:four_cycle8}}
\end{figure}

Else if state $c_2$ is adjacent to two nonadjacent edges in $G(S)$ (Figure \ref{fig:four_cycle8}(a) and \ref{fig:four_cycle8}(b)), then this again creates a chordless cycle on two colors as shown in Figure \ref{fig:four_cycle8}(c), contradicting condition (ii).

\begin{figure}[h!]
        \centering
        \includegraphics[scale=0.35]{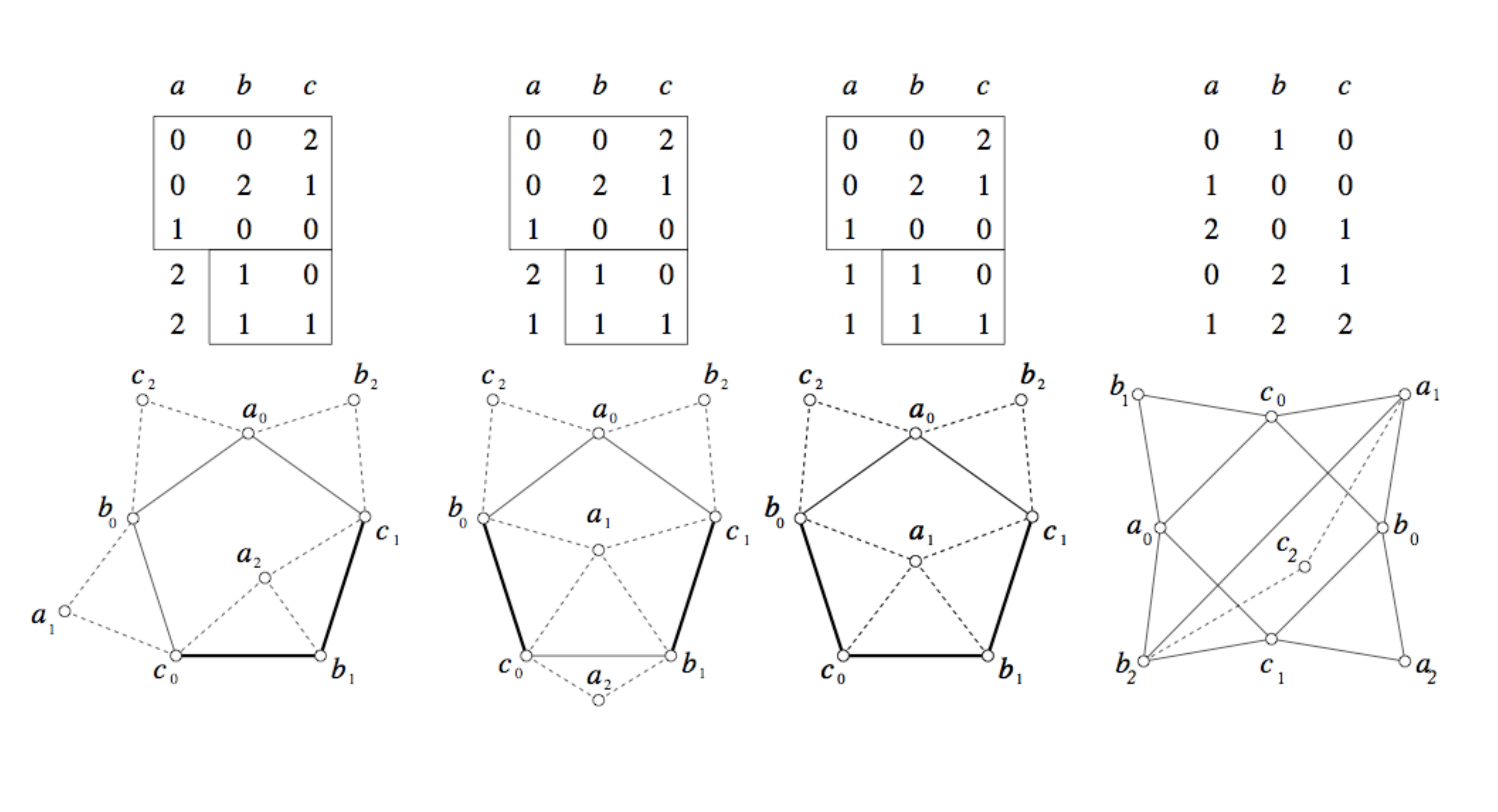}
        \caption{Minimal obstruction sets for three-state characters up to relabeling.\label{fig:obstructionsets}}
\end{figure}

In summary, Figure \ref{fig:obstructionsets} shows the minimal obstruction sets to the existence of perfect phylogenies for three-state characters up to relabeling of the character states.

\section{Structure of Proper Triangulations for Partition Intersection Graphs on Three-State Characters}

The complete description of minimal obstruction sets for three-state characters allows us to expand upon recent work of Gusfield which uses properties of legal triangulations and minimal separators of partition intersection graphs to solve several problems related to multi-state perfect phylogenies \cite{Gu09}.  In particular, the following is a necessary and sufficient condition for the existence of a perfect phylogeny for multi-state data.  We refer the reader to \cite{Gu09} for the necessary definitions and the proof.

\begin{theorem}[\cite{Gu09} Theorem 3 MSPN] Input $S$ allows a perfect phylogeny if and only if there is a set $Q$ of pairwise parallel legal minimal separators in partition intersection graph $G(S)$ such that every mono-chromatic pair of nodes in $G(S)$ is separated by some separator in $Q$. \label{thm:MSPN} \end{theorem}

For the special case of input $S$ with characters over three states, the construction of minimal obstruction sets in Section \ref{sec:forbidden} allows us to simplify Theorem \ref{thm:MSPN} to the following.

\begin{theorem} \label{thm:3statemono} For input $S$ on at most three states per character $(r \leq 3)$, there is a three-state perfect phylogeny for $S$ if and only if the partition intersection graph for every pair of characters is acyclic and every mono-chromatic pair of nodes in $G(S)$ is separated by a legal minimal separator.
\end{theorem}

% For input $S$ on at most three states per character $(r \leq 3)$, there is a three-state perfect phylogeny for $S$ if 
%and only if every pair of characters allows a perfect phylogeny and every mono-chromatic pair of nodes in $G(S)$ 
%is separated by a legal minimal separator. 

\Proof  By Theorem \ref{thm:main}, if three-state input $S$ does not allow a perfect phylogeny, then there is a triple of characters $a,b,c$ in $S$ that does not allow a perfect phylogeny.  By Section \ref{sec:forbidden}, if every pair of characters $(a,b), (a,c)$ and $(b,c)$ is acyclic (and therefore allow a perfect phylogeny) while $a,b$, and $c$ together do not allow a perfect phylogeny, then one of the graphs in Figure \ref{fig:obstructionsets} appears as an induced graph in partition intersection graph $G(S)$ (up to relabeling).  In each of these graphs, consider the mono-chromatic pair of vertices $c_0$ and $c_1$; we will show that in each graph, no legal separator can separate these two vertices.   

\begin{figure}[h!]
      	\centering
	\includegraphics[scale=0.35]{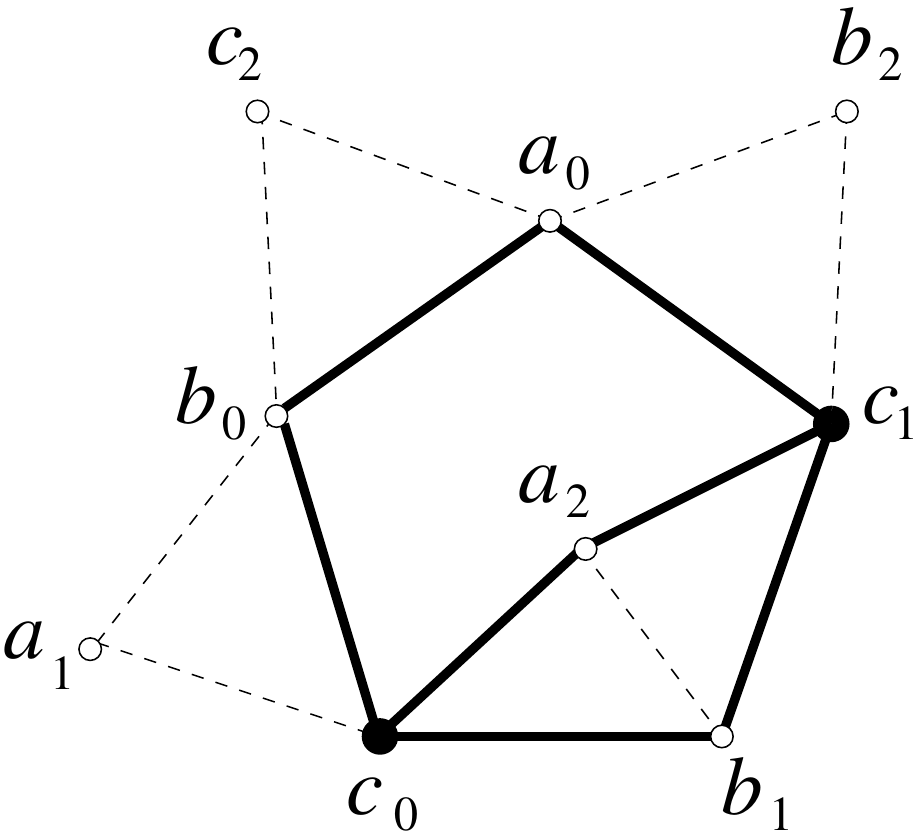}	
	\caption{Disjoint paths from $c_0$ to $c_1$\label{fig:threestate_mono1}}
\end{figure}

In the first graph, consider the following three vertex disjoint paths from $c_0$ to $c_1$ (shown in Figure \ref{fig:threestate_mono1})

\begin{enumerate}
\item $c_0 \rightarrow a_2 \rightarrow c_1$
\item $c_0 \rightarrow b_1 \rightarrow c_1$.  
\item $c_0 \rightarrow b_0 \rightarrow a_0 \rightarrow c_1$
\end{enumerate}

From these disjoint paths, we see that any separator $Q$ for $c_0$ and $c_1$ must include vertices $a_2$ and $b_1$ (to destroy paths (1) and (2)).  Furthermore, $Q$ must also contain one of $a_0$ or $b_0$ (to destroy path (3)), implying that $Q$ must contain at least two states of the same color and therefore cannot be a legal separator. 

\begin{figure}[h!]
      	\centering
	\includegraphics[scale=0.35]{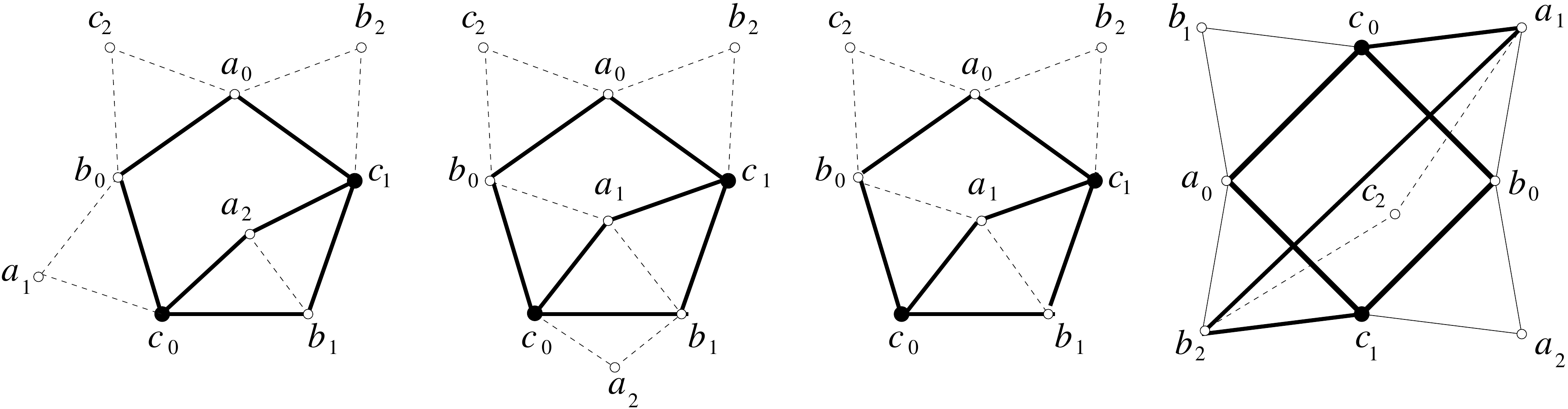}	
	\caption{Disjoint paths showing all separators of $c_0$ and $c_1$ are illegal. \label{fig:threestate_mono}}
\end{figure}

A similar argument using the disjoint paths in Figure \ref{fig:threestate_mono} shows that any separator $Q$ for vertices $c_0$ and $c_1$ cannot be a legal separator.  Therefore, if every monochromatic pair is separated by a legal minimal separator, none of these graphs in Figure \ref{fig:obstructionsets} appears as a subgraph of the partition intersection graph, implying that every subset of three characters allows a perfect phylogeny.  Theorem \ref{thm:main} then shows the entire set of characters allows a perfect phylogeny.

The other direction of the theorem follows from Theorem \ref{thm:MSPN}.   \qed

Theorem \ref{thm:3statemono} implies that the requirement of Theorem \ref{thm:MSPN} that the legal minimal separators in $Q$ be pairwise parallel, can be removed for the case of input data over three-state characters.  

\comment{Theorem \ref{thm:3statemono} shows that the case of three state input is special, as there are examples showing that the theorem does not extend to input data with four state characters.}

\section{Construction of Fitch-Meacham Examples}\label{sec:fitchmeacham}

In this section, we examine in detail the class of Fitch-Meacham examples, which were first introduced by Fitch \cite{F75,F77} and later generalized by Meacham \cite{Me83}.  The goal of these examples is to demonstrate a lower bound on the number of characters that must be simultaneously examined in any test for perfect phylogeny.  The natural conjecture generalizing our main result is that for any integer $r \geq 3$, there is a perfect phylogeny on $r$-state characters if and only if there is a perfect phylogeny for every subset of $r$ characters.  We show here that such a result would be the best possible, for any $r$.  While the general construction of these examples and the resulting lower bounds were stated by Meacham~\cite{Me83}, to the best of our knowledge, the proof of correctness for these lower bounds has not been established.  We fill this gap by explicitly describing the complete construction for the entire class of Fitch-Meacham examples and providing a proof for the lower bound claimed in \cite{Me83}.  

For each integer $r$ ($r \geq 2$), the Fitch-Meacham construction $F_r$ is a set of $r+2$ sequences over $r$ characters, where each character takes $r$ states.  We describe the construction of the partition intersection graph $G(F_r)$; the set of sequences $F_r$ can be obtained from $G(F_r)$ in a straightforward manner, with each taxon corresponding to an $r$-clique in $G(F_r)$.   

Label the $r$ characters in $F_r$ by $0,1, \ldots r-1$; each vertex labeled by $i$ will correspond to a state in character $i$.  The construction starts with two cliques EC$_1$ and EC$_2$ of size $r$, called end-cliques, with the vertices of each clique labeled by $0,1, \ldots r-1$.  The vertex labeled $i$ in EC$_1$ is adjacent to the vertex labeled $(i+1)\mod r$ in EC$_2$.  For each such edge $(i, (i+1)\mod r)$ between the two end-cliques, we create a clique of size $r-2$ with vertices labeled by $\{0,1, \ldots r-1\} \backslash \{i, (i+1)\mod r \}$.  Every vertex in this $(r-2)$-clique is then attached to both $i$ (in end-clique 1) and $(i+1)\mod r$ (in end-clique 2), creating an $r$-clique whose vertices are labeled with integers $0,1, \ldots r-1$.  There are a total of $r$ such cliques,
called \emph{tower-cliques}, and denoted by TC$_1$, TC$_2$, \ldots TC$_r$.  Note that for each $i$ ($0 \leq i \leq r-1$), there are exactly $r$ vertices labeled by $i$; we give each such vertex a distinct state, resulting in $r$ states for each character.

\comment{
This describes the construction of $F_r$ for all $r$, but we can also view this construction inductively by creating $F_{r+1}$ from $F_r$ as follows.  Starting from $F_r$, add a new vertex labeled $r+1$ to each end-clique and each tower-clique, adding edges to make $r+1$ adjacent to all vertices in the clique to obtain cliques of size $r+1$.  In the resulting graph, there are still $r$ tower-cliques, and the two vertices labeled $r+1$ in the end-cliques are not part of any tower.  Now, move the tower with end vertices labeled $r-1$ and $0$ in the first and second end-cliques respectively, to become a tower whose end vertices are $r-1$ and $r$.  Then create an edge between the vertices $r$ and $0$ in the first and second end-cliques respectively and create a tower-clique of size $r+1$ using this edge.  Label the new vertices in this tower-clique with the integers $1,2, \ldots r-1$. The resulting graph is $F_{r+1}$.
}

Note that the graph corresponding to the four gamete obstruction set is an instance of the Fitch-Meacham construction with $r=2$.  In this case, the four binary sequences $00, 01, 10, 11$ have two states, two colors and four taxa and the partition intersection graph for these sequences is precisely the graph $G(F_2)$.  Note that in this case, every subset of $r-1=1$ characters has a perfect phylogeny, while the entire set of characters does not.  Similarly, the fourth graph shown in Figure \ref{fig:obstructionsets} illustrating the obstruction set for 3-state input is the graph $G(F_3)$ corresponding to the Fitch-Meacham construction for $r=3$ (in the figure, EC$_1$= $\{ a_0,b_2,c_1\}$ and EC$_2$ = $\{ a_1,b_0,c_0\}$).  As shown in Section \ref{sec:forbidden}, every $r-1 = 2$ set of characters in the corresponding input set allows a perfect phylogeny while the entire set of characters does not.  The following theorem generalizes this property to the entire class of Fitch-Meacham examples.  Because the theorem was stated without proof in \cite{Me83}, we provide a proof of the result here.

\begin{theorem}\cite{Me83} For every $r \geq 2$, $F_r$ is a set of input sequences over $r$ state characters such that every $r-1$ subset of characters allows a perfect phylogeny while the entire set $F_r$ does not allow a perfect phylogeny.
\end{theorem}

\Proof We first show that $G(F_{r})$ does not allow a proper triangulation for any $r$.  As observed above, $G(F_2)$ is a four cycle corresponding to two colors and therefore, does not allow a proper triangulation (since any proper triangulation for a graph containing cycles must have at least three colors).  Suppose $G(F_r)$ is properly triangulatable for some $r \geq 3$, let $s$ be the smallest integer such that $G(F_s)$ has a proper triangulation, and let $G'(F_s)$ be a minimal proper triangulation of $G(F_s)$.  

For each tower-clique TC$_i$ in $G(F_s)$, consider the set of vertices in TC$_i$ that are not contained in either end-clique; call these vertices \emph{internal tower-clique vertices} and the remaining two tower vertices  \emph{end tower-clique vertices}.  Note that the removal of the two end tower-clique vertices disconnects the internal tower-clique vertices from the rest of the graph.  This implies that the internal tower-clique vertices cannot be part of any chordless cycle: otherwise, such a chordless cycle $C$ must contain \emph{both} end tower-clique vertices $i$ and $(i+1)\mod s$.  However, the two end tower-clique vertices are connected by an edge and therefore induce a chord in $C$,  a contradiction since $C$ is a chordless cycle.  

\comment{
For any subgraph $H$ of $G(F_s)$, suppose $H$ allows a proper triangulation and let $H'$ be a minimal proper triangulation for $H$.  Consider the added fill-in edges $A = E(H') \backslash E(H)$.  We first claim no internal tower-clique vertex can be adjacent to any edge in $A$.  Otherwise, consider the graph $H''$ obtained from $H'$ by removing all edges in $A$ that are adjacent to internal tower-clique vertices.  By the observation above, the internal tower-clique vertices cannot be contained in any chordless cycle in $H'$ and therefore, any chordless cycle in $H'$ remains a chordless cycle in $H''$.  Therefore, graph $H''$ is chordal, contradicting the minimality of $H'$.  This shows that the graph graph $G(F_s)$ allows a proper triangulation if and only if the graph with all internal tower-clique vertices removed allows a proper triangulation.  In this graph, we have two end-cliques on $0,1, \ldots s$ with edges (rather than tower cliques) connecting vertex $i$ (in end-clique 1) to  vertex $(i+1)\mod s$ (in end-clique 2).  
}

In the graph $G(F_s)$, onsider the following cycle of length four: $s-2$ (in EC$_1$) $\rightarrow s-1$ (in EC$_1$) $\rightarrow 0$ (in EC$_2$) $\rightarrow s-1$ (in EC$_2$) $\rightarrow s-2$ (in EC$_1$).  This four-cycle has a unique proper triangulation, which forces the edge $e$ between vertex $s-2$ in EC$_1$ and vertex $0$ in EC$_2$ to be included in $G'(F_s)$.  Consider adding edge $e$, removing all vertices labeled $s-1$ from $G'(F_s)$, and for the two vertices labeled $s-1$ in end-cliques EC$_1$ and EC$_2$, remove all interior tower-clique vertices (but not end tower-clique vertices) adjacent to $s-1$.  Then edge $e$ between vertices $s-2$ and $0$ is still present and we can expand $e$ into a tower-clique of size $s-1$ (by forming a clique with new vertices $1,2, \ldots s-3$ adjacent to both $s-2$ and $0$ of the two end-cliques).  

In the resulting graph, the vertices are exactly those of $G(F_{s-1})$ and all edges in $G(F_{s-1})$ are present.  Furthermore, if there is a chordless cycle in this graph, then it would create a chordless cycle in $G'(F_s)$ since no internal tower-clique vertex can be part of any chordless cycle (and in particular, the new vertices $1,2, \ldots s-3$ cannot be part of any chordless cycle).  Therefore, the resulting graph is a proper triangulation for $G(F_{s-1})$, a contradiction since $s$ was chosen to be the smallest integer such that $G(F_s)$ allows a proper triangulation.

To prove the second part of the theorem, we show that in $F_r$, any subset of $r-1$ characters does allow a perfect phylogeny by proving that the partition intersection graph on any subset of $r-1$ characters has a proper triangulation.  By the symmetry of the construction of $F_r$, we can assume without loss of generality that the $r-1$ characters under consideration are $\{0,1, \ldots r-2\}$.   Consider the graph obtained by connecting every vertex $i$ ($0 \leq i \leq r-3$) in EC$_1$ to every vertex $j$ satisfying $j > i$ in EC$_2$. Note the asymmetry between the first and second end-cliques in this construction and observe that none of the added edges are between characters with the same label.

Suppose the resulting graph contains a chordless cycle $C$.  Then $C$ cannot contain three or more vertices in either end-clique and cannot contain any internal tower-clique vertices (as noted earlier), so must have length exactly four with two vertices in each end-clique.   It cannot be the case that two nonadjacent vertices of $C$ are in the same end-clique, since these vertices would be adjacent and $C$ would not be chordless.  Therefore, cycle $C$ must be formed as follows: $i$ (in EC$_1$) $\rightarrow j$ (in EC$_2$) $\rightarrow j'$ (in EC$_2$) $\rightarrow i'$ (in EC$_1$).  Since $i$ and $j$ are adjacent, we have $i < j$ and since $i'$ and $j'$ are adjacent, we have $i' < j'$.  If $i<j'$, then $i$ and $j'$ are adjacent and the cycle $C$ is not chordless, a contradiction.  Therefore,  $i' < j' \leq i < j$,  which implies $i'$ and $j$ are adjacent and the cycle $C$ is not chordless, again a contradiction.  It follows that there are no chordless cycles and the added edges form a proper triangulation for the partition intersection graph on the subset of $r-1$ characters $\{0,1,\ldots r-2\}$. \qed

\section{Conflict Hypergraphs for Three State Characters}

For binary data, the Four Gamete Condition/Splits Equivalence Theorem implies that the existence of a perfect phylogeny can be determined by a pairwise test.  This motivates the following definition.

\begin{definition} For binary input $S$ on characters $\mathcal{C}$, two characters $\chi_i, \chi_j \in \mathcal{C}$ are \emph{in conflict} or \emph{incompatible} if $\chi_i$ and $\chi_j$ contain all four gametes.  The \emph{conflict}(or \emph{incompatibility}) graph $C(S)$ for $S$ is defined on vertices $V$ and edges $E$ with
\begin{eqnarray*}
V & = & \{ \chi_i  : \chi_i \in \mathcal{C} \} \\
E & = & \{ (\chi_i, \chi_j) : \chi_i \text{ and } \chi_j \text{ are incompatible} \}
\end{eqnarray*}
\end{definition}

By the Four Gamete Condition, binary input $S$ allows a perfect phylogeny if and only if the corresponding incompatibility graph $C(S)$ has no edges.  For input containing incompatible characters, there has been extensive literature devoted to studying the structure of incompatibility graphs and using these graphs to solve related problems.  For example, Gusfield et al.~\cite{GB05,GBBS07} and Huson et al.~\cite{H+05} use incompatiblity graphs to achieve decomposition theorems for phylogenies, and Gusfield, Hickerson, and Eddhu~\cite{GHE07}, Bafna and Bansal~\cite{BB05,BB04}, and Hudson and Kaplan~\cite{HK85} use it to achieve lower bounds on the number of recombination events needed to explain a set of sequences.  

The incompatibility graph is also the basis for algorithms to solve the \emph{character removal} or \emph{maximum compatibility} problem, which asks for the minimum number of characters that must be removed from an input set such that the remaining characters allow a perfect phylogeny.  For binary data, the character removal problem is equivalent to the vertex cover problem on the corresponding incompatibility graph \cite{Fe04,SS03}.  

Until this work, the notion of incompatibility graph was defined for binary characters only, as motivated by the Splits Equivalence Theorem.  Our generalization of the Splits Equivalence Theorem therefore allows us to generalize in a natural way the notion of incompatibility for three state characters. The resulting incompatibility structure will be a hypergraph whose edges correspond to pairs and triples of characters that do not allow a perfect phylogeny.  In particular, for input data $S$ on three state characters $\mathcal{C}$, let $E_2(S)$ be the set of character pairs $(\chi_i, \chi_j)$  such that $\chi_i$ and $\chi_j$ do not allow a perfect phylogeny and let $E_3(S)$ be the set of triples $(\chi_i, \chi_j, \chi_k)$ such that $\chi_i, \chi_j, \chi_k$ together do not allow a perfect phylogeny but each of the pairs $(\chi_i, \chi_j), (\chi_i, \chi_k), (\chi_j, \chi_k)$ allows a perfect phylogeny (i.e., $(\chi_i, \chi_j), (\chi_i, \chi_k), (\chi_j, \chi_k) \not\in E_2(S)$).  Then the incompatibility hypergraph $C(S)$ is defined on vertices $V$ and hyperedges $E$ with

\begin{eqnarray*}
V & = & \{ \chi_i  : \chi_i \in \mathcal{C} \} \\
E & = & E_2(S) \cup E_3(S)
\end{eqnarray*}

Note that by definition, no hyperedge of $C(S)$ is contained in another hyperedge.  The extension of incompatibility to three state characters can be used to solve algorithmic and theoretical problems for three state characters analogous to those for binary characters.  In particular, the character removal problem on three state characters can be solved using the following generalization of the vertex cover problem.

\begin{itemize}
\item[] {\bf 3-Hitting Set Problem}
\item[] {\bf Input:} A collection $M$ of subsets of size at most three from a finite ground set $\Omega$ and a positive integer $k$
\item[] {\bf Problem:} Determine if there is a set $L \subseteq \Omega$ with $\vert L \vert \leq k$ such that $L$ contains at least one element from each subset in $M$.
\end{itemize}

The 3-hitting set problem is {\tt NP}-complete \cite{GJ79} and the best known approximation algorithms for the problem have approximation ratio equal to three \cite{Ho97}.  More recently, it has been shown that the 3-hitting set problem is fixed parameter tractable in the parameter $k$.  In a series of papers, algorithms for solving the 3-hitting set problem have been given with running times $O(2.270^k + n)$ \cite{NR03}, $O(2.179^k + n)$ \cite{Fe04}, and $O(2.076^k+n)$ \cite{Wa07}.  It has also been shown that the 3-hitting set problem allows a linear-time \emph{kernelization}, a preprocessing step typical for parameterized algorithms that converts a problem of input size $n$ to an instance whose size depends only on $k$.  In \cite{Ab07}, it is shown that any 3-Hitting-Set instance can be reduced into an equivalent instance of size at most $5k^2 + k$ elements. 

Using these results for the 3-hitting set problem, we obtain the following.

\begin{corollary} The Character Removal Problem for three-state input is fixed parameter tractable.  
\end{corollary}

\section{Conclusion}

We have studied the structure of the three state perfect phylogeny problem 
and shown that there is a necessary and sufficient condition for the existence of 
a perfect phylogeny for three state characters using triples of characters. This 
extends the extremely useful Splits Equivalence Theorem and Four Gamete Condition. The obvious extension of our work would be to discover similar results 
for $r$-state characters for $r \geq 4$. 

Until this work, the notion of a conflict, or incompatibility, graph has been 
defined for two state characters only. Our 
generalization of the four gamete condition allows us to generalize this 
notion to incompatibility on three state characters. The resulting incompatibility 
structure is a hypergraph, which can be used to solve algorithmic and theoretical problems for three state characters analogous to those for binary characters. 

In addition, there are several theoretical and practical results known for two 
state characters that are still open for characters on three or more states. For 
instance, it is known that the problem of constructing near-perfect phylogenies 
for two state characters is fixed parameter tractable; the analogous problem is 
open for characters on three or more states. Similarly, the question of whether incompatibility hypergraphs can be used to find decomposition theorems for lower bounding recombination events remains open for three or more states. With the recent increase in collection of polymorphism data such as micro/mini-satellites, there is a need for the analysis of perfect 
phylogenies to be extended to multiple state characters. Our work lays a solid 
theoretical foundation we hope will help with this effort.

\vspace{2ex}

\noindent {\bf Acknowledgments} The authors gratefully acknowledge G. Blelloch, I. Coskun, R. Gysel, R. Ravi, R. Schwartz, and T. Warnow for stimulating discussions and suggestions.  This research was partially supported by NSF grants SEI-BIO 0513910, CCF-0515378, and IIS-0803564.

\end{document}